\documentclass[fleqn,11pt,a4paper,leqno]{amsart}

\usepackage{verbatim}
\usepackage{graphicx}
\usepackage{graphics}
\usepackage{latexsym}
\usepackage{amssymb,amsmath, MnSymbol}
\usepackage{multicol}
\usepackage{bussproofs}
\usepackage{enumerate}
\usepackage{xcolor}
\usepackage[latin1]{inputenc}

\def\ismodeledby{=\joinrel\mathrel|}
\date{}

\numberwithin{equation}{section}

\newtheorem{definicion}{Definition}[section]

\newtheorem{definition}[definicion]{Definition}
\newtheorem{Lemma}[definicion]{Lemma}

\newtheorem{Theorem}[definicion]{Theorem}
\newtheorem{theorem}[definicion]{Theorem}
\newtheorem{Corollary}[definicion]{Corollary}

\newtheorem{lema}[definicion]{Lemma}

\newtheorem{Remark}[definicion]{Remark}

\newenvironment{Proof}{\noindent\bf Proof \rm}{$\hfill
\square$}

\newcommand{\debaj}[2]{ #1 \to_H #2}

\newcounter{ecuacionDef} 
\renewcommand{\theecuacionDef}{\arabic{ecuacionDef}}

\newenvironment{ecuacionDef}[1]%
{
	\refstepcounter{ecuacionDef}
	\vspace{0.12cm}
	{\rm (E\theecuacionDef)}:  #1
	\vspace{0.12cm}
}%

\newcounter{numeroAxioma} 
\renewcommand{\thenumeroAxioma}{\arabic{numeroAxioma}}

\newenvironment{numeroAxioma}[1]%
{
	\refstepcounter{numeroAxioma}
	\vspace{0.12cm}
	(A\thenumeroAxioma):  #1
	\vspace{0.12cm}
}%

\begin{document}

 \newcounter{thlistctr}
 \newenvironment{thlist}{\
 \begin{list}%
 {\alph{thlistctr}}%
 {\setlength{\labelwidth}{2ex}%
 \setlength{\labelsep}{1ex}%
 \setlength{\leftmargin}{6ex}%
 \renewcommand{\makelabel}[1]{\makebox[\labelwidth][r]{\rm (##1)}}%
 \usecounter{thlistctr}}}%
 {\end{list}}

\title[A Logic for Dually Hemimorphic Semi-Heyting Algebras and Axiomatic Extensions]{A Logic for Dually Hemimorphic Semi-Heyting Algebras \\ and\\  its  Axiomatic Extensions}  

\author[J. M. CORNEJO \and H. P. SANKAPPANAVAR]{Juan M. CORNEJO* \and Hanamantagouda P. SANKAPPANAVAR**}

\newcommand{\acr}{\newline\indent}

\address{\llap{*\,} Departamento de Matem\'atica\acr
Universidad Nacional del Sur\acr
Alem 1253, Bah\'ia Blanca, Argentina\acr
INMABB - CONICET}

\email{jmcornejo@uns.edu.ar}

\address{\llap{**\,}Department of Mathematics\acr
                              State University of New York\acr
                              New Paltz, New York, 12561\acr
                              U.S.A.}

\email{sankapph@newpaltz.edu}




\begin{abstract}
The variety $\mathbb{DHMSH}$ of dually hemimorphic semi-Heyting algebras was introduced in 2011 by the second author 
as an expansion of semi-Heyting algebras by a  dual hemimorphism. 
In this paper, we focus on the variety $\mathbb{DHMSH}$ from a logical point of view.  Firstly, we present a Hilbert-style axiomatization of a new 
logic called ``Dually hemimorphic semi-Heyting logic'' ($\mathcal{DHMSH}$, for short), as an expansion of semi-intuitionistic logic $\mathcal{SI}$ (also called $\mathcal{SH}$) introduced by the first author  
by adding a weak negation (to be interpreted as a dual hemimorphism).  We then prove that it is implicative in the sense of Rasiowa and that it is complete with respect to the variety $\mathbb{DHMSH}$.  
It is deduced that the logic $\mathcal{DHMSH}$ is algebraizable in the sense of Blok and Pigozzi, with 
  the variety $\mathbb{DHMSH}$ as its equivalent algebraic semantics 
and that the lattice of axiomatic extensions of $\mathcal{DHMSH}$ is dually isomorphic to the lattice of subvarieties of $\mathbb{DHMSH}$.  A new
axiomatization for Moisil's logic is also obtained.
Secondly, we characterize the axiomatic extensions of  
$\mathcal{DHMSH}$ in which the ``Deduction Theorem'' holds. 
Thirdly, we present several new logics, extending the logic $\mathcal{DHMSH}$, corresponding to several important subvarieties of the variety $\mathbb{DHMSH}$.  
These include logics corresponding to the varieties generated by  
 two-element, three-element and some four-element dually quasi-De Morgan semi-Heyting algebras, as well as a new axiomatization for the 3-valued \L ukasiewicz logic.
Surprisingly, many of these logics turn out to be connexive logics, a few of which are presented in this paper.  
Fourthly, we present axiomatizations for two infinite sequences of logics namely,  De Morgan-G\"{o}del logics and dually pseudocomplemented G\"{o}del logics,  
Fifthly, axiomatizations are also provided for logics corresponding to many subvarieties of regular dually quasi-De Morgan Stone semi-Heyting algebras, of regular De Morgan semi-Heyting algebras of level 1, and of JI-distributive semi-Heyting algebras of level 1.
We conclude the paper with some open problems.  
Most of the logics considered in this paper are discriminator logics in the sense that they correspond to discriminator varieties.  Some of them, just like the classical logic, are even primal in the sense that their corresponding varieties are generated by primal algebras. 
\end{abstract}

\maketitle


\section{Introduction} \label{SecOne}

\

\indent  Semi-Heyting algebras were introduced by the second author\footnote {Parts of this paper were presented by the second author in invited talks at 8th International Conference on
 Non-Classical Logics: Theory and Applications, 
$\L$odz (2016), at Maltsev Meeting, Novosibirsk (2017), and at  Asubl (Algebra and Substructural Logics-Take 6) workshop, Cagliari (2018).},  
 during 1983-84, as a result of his research that went into~\cite{sankappanavar1985Heyting} (still a preprint at the time). 
Some of the early results were announced in~\cite{San87a}.  The first results on 
these algebras with their proofs, however, were published much later in 2008 (see \cite{sankappanavar2008semi}).

An algebra $\mathbf L = \langle L, \lor, \land, \to, 0, 1\rangle$ is a semi-Heyting algebra if the following conditions hold:
\begin{enumerate}[\qquad (SH1)]
	\item  $\langle L, \lor, \land, 0,1\rangle$ is a lattice with $0$ and $1$,
	\item $x \land (x \to y) \approx x \land y$,
	\item $x \land (y \to z) \approx x \land [(x \land y) \to (x \land z)]$,
	\item $x \to x \approx 1$.
\end{enumerate}
\vspace{-.1cm}

A semi-Heyting algebra is a Heyting algebra if it satisfies the identity: 

\indent \quad (H) \quad $(x \land y) \to x \approx 1$.\\

We will denote the variety of semi-Heyting algebras by $\mathbb{SH}$ and that of Heyting algebras by $\mathbb{H}.$
Semi-Heyting algebras share some important properties with Heyting algebras; for instance, semi-Heyting algebras are distributive and pseudocomplemented, with the pseudocomplement $x^* := x \to 0$; the congruences on them are determined by filters and the variety of semi-Heyting algebras is arithmetical.  For further results on $\mathbb{SH}$, see  \cite{Abad2011a}, \cite{Abad2011},  \cite{Abad2013semi}, \cite{cornejo20semiHeyting} and \cite{sankappanavar2008semi}.

It is well known that the variety of Heyting algebras is the equivalent algebraic semantics (in the sense of Blok and Pigozzi) of the intuitionistic propositional logic.  In 2011, the first author of this paper  defined, in \cite{cornejo2011semi}, a new logic called ``semi-intuitionistic logic'' ($\mathcal{SI}$, for short, also called $\mathcal{SH}$) and showed, essentially, that the variety of semi-Heyting algebras is an algebraic semantics for this logic and that the intuitionistic logic is an axiomatic extension of it.  The axioms of this logic, however, were  
expressed in a language that was not the same as 
that of semi-Heyting algebras.
In \cite{cornejo15semiLogics}, a much simpler, but equivalent, set of axioms for $\mathcal{SI}$ (or $\mathcal{SH}$), was presented in the same language as that of semi-Heyting algebras.  The logic $\mathcal{SI}$ as presented in \cite{cornejo15semiLogics} will play a fundamental role in this paper.  

In 1942, Moisil \cite{moisil1942} defined a logic called ``Logique modale'' ($\mathcal{LM}$), an expansion of intuitionistic propositional calculus by a De Morgan negation.   He also introduced Heyting algebras endowed with an involution, in \cite{moisil1942a}, as the algebraic models of the logic $\mathcal{LM}$.  
These algebras were further investigated by
 Monteiro \cite{monteiro1980algebres} under the name of symmetric Heyting algebras.  In particular, he presented a proof of an algebraic completeness theorem  for Moisil's calculus  
 by showing that $\mathcal{LM}$ is complete for the variety of symmetric Heyting algebras.  
 
Independently of the previous work,  
 motivated purely by (universal) algebraic considerations, the second author defined and studied De Morgan Heyting algebras, in \cite{sankappanavar1987dualendo}, by expanding Heyting algebras by a De Morgan operation.  
Earlier in 1985, he had also introduced (see \cite{sankappanavar1985Heyting}) the variety of Heyting algebras with a dual pseudocomplementation.   Also, in 1987, the concepts of hemimorphism (without name), semi-De Morgan algebra and (lower) quasi-De Morgan algebra were introduced  in \cite{sankappanavar1985semiDeMorgan}, unifying (and generalizing) the notions of De Morgan operation and pseudocomplementation.

In 2011, motivated by the similarities of the results and proofs in \cite{sankappanavar1985Heyting} and \cite{sankappanavar1987dualendo}, he introduced in \cite{sankappanavar2011expansions} a more general variety of algebras called ``dually hemimorphic semi-Heyting algebras''--
 an expansion of semi-Heyting algebras by a dual hemimorphism, as a common 
generalization of De Morgan Heyting algebras and dually psedocomplemented Heyting algebras.   

\begin{definicion}  {\rm \cite{sankappanavar2011expansions}} \label{definition_DHMSH}
	An algebra $\mathbf A=\langle A, \lor, \land, \to, ', 0,1\rangle$ is a dually hemimorphic semi-Heyting algebra {\rm(or,} semi-Heyting algebra with a
	dual hemimorphism\rm) if $\mathbf A$ satisfies the following conditions:

 \ecuacionDef{$\langle A, \lor, \land, \to, 0,1\rangle $ is a semi-Heyting algebra,}

\ecuacionDef{$0' \approx 1$,} \label{conditionSD2d1}

\ecuacionDef{$1' \approx 0$,} \label{conditionSD2d2}

\ecuacionDef{$(x \land y)' \approx x' \lor y'$ \quad {\rm ($\land$-De Morgan law)}.} \label{infDMlaw} \label{conditionSD3d}

The unary operation $'$ satisfying (E2)-(E4) is called a dual hemimorphism.
The variety of dually hemimorphic semi-Heyting algebras will be denoted by $\mathbb{DHMSH}$.
\end{definicion} 

Several important subvarieties of the variety 
$\mathbb{DHMSH}$, by adding the duals of those given in \cite{sankappanavar1985semiDeMorgan}, were introduced in {\rm \cite{sankappanavar2011expansions}}, some of which will be recalled in Section \ref{SecFive}.  
  
\smallskip
The following problem presents itself naturally.\\

PROBLEM 1: Find a propositional logic in the language $\langle \lor, \land, \to, \sim, \bot, \top \rangle$ 
with the following properties:

(1) It has the variety $\mathbb{DHMSH}$ of dually hemimorphic semi-Heyting algebras   
as its equivalent algebraic semantics, and 

(2) It has Moisil's logic as one of its (axiomatic) extensions.\\ 

The subvariety $\mathbb{DQDSH}$ of $\mathbb{DHMSH}$, consisting of dually quasi-De Morgan semi-Heyting algebras (see {\bf LIST 1} in Section \ref{SecFive} 
for definition), has been intensively investigated in \cite{sankappanavar2011expansions, sankappanavar2014a, sankappanavar2014b, sankappanavar2016, sankappanavar2017, sankappanavar2019}.  In Section 8 of \cite{sankappanavar2011expansions} (see also \cite{sankappanavar2014b} and \cite{sankappanavar2016})
 the following problem was raised:\\  

PROBLEM 2: Find Hilbert-type axiomatization for logics corresponding to  
 two-valued, three-valued and four-valued dually quasi-De Morgan semi-Heyting algebras, viewed as logical matrices with $\{1\}$ as the distinguished subset. \\

In this paper, we focus on the logical aspects of the variety $\mathbb{DHMSH}$ of dually hemimorphic semi-Heyting algebras and many of its subvarieties.  Firstly, we give a solution to PROBLEM 1.  More specifically, we present a Hilbert-style presentation of a new logic called ``Dually hemimorphic semi-Heyting logic'' ($\mathcal{DHMSH}$, for short), as an expansion of semi-intuitionistic logic 
presented in \cite{cornejo15semiLogics}.  We then
prove that it is implicative in the sense of Rasiowa and that it is complete with respect to the variety $\mathbb{DHMSH}$ of dually hemimorphic semi-Heyting algebras.  Using the well-known results of Abstract Algebraic Logic we deduce that 
the logic $\mathcal{DHMSH}$ is algebraizable in the sense of Blok and Pigozzi, with 
the variety $\mathbb{DHMSH}$ as its equivalent algebraic semantics.
It then follows that the lattice of axiomatic extensions of $\mathcal{DHMSH}$ is dually isomorphic to the lattice of subvarieties of $\mathbb{DHMSH}$. 
As applications of these results, we present several new logics, extending the logic $\mathcal{DHMSH}$, corresponding to some interesting subvarieties (studied in \cite{sankappanavar2011expansions}) 
 of the variety of hemimorphic semi-Heyting and Heyting algebras.  A new
axiomatization for Moisil's logic is also obtained.
Secondly, we characterize the axiomatic extensions of 
$\mathcal{DHMSH}$ in which the ``Deduction Theorem'' holds. 
Thirdly, we introduce several new logics, extending the logic $\mathcal{DHMSH}$, corresponding to  important subvarieties of the variety $\mathbb{DHMSH}$,  
including some logics corresponding to the varieties generated by  
 two-element, three-element and some four-element dually quasi-De Morgan semi-Heyting algebras, as well as a new axiomatization for the 3-valued \L ukasiewicz logic.
Many of these logics turn out to be connexive logics, a few of which are presented in this paper.  
Fourthly, we present axiomatizations for two infinite sequences of logics, namely  De Morgan-G\"{o}del logics and dually pseudocomplemented G\"{o}del logics,  
Fifthly, axiomatizations are also provided for logics corresponding to many subvarieties of regular dually quasi-De Morgan Stone semi-Heyting algebras of level 1, of Regular De Morgan Semi-Heyting Algebras of level 1 and of JI-distributive semi-Heyting algebras of level 1.  
Many of the logics considered in this paper are discriminator logics in the sense that they correspond to discriminator varieties.  Some of them, just like the classical logic, are even primal in the sense that their corresponding varieties are generated by primal algebras.  

The paper is organized as follows:  Section \ref{section_Preliminaries} contains definitions, notation and some preliminary results that are needed later in the paper.  It includes the axiomatization 
for semi-intuitionistic logic as presented in \cite{cornejo15semiLogics} which is crucial for the rest of the paper.  In Section \ref{section_logicDHMSH}, we present a Hilbert-style axiomatization for the new logic called ``Dually hemimorphic semi-Heyting logic'' ($\mathcal{DHMSH}$, for short)'' 
by expanding the language of semi-intuitionistic logic $\mathcal{SI}$ of \cite{cornejo15semiLogics} by a (weak) negation called {\em dually hemimorphic negation} and 
by adding new axioms and a new inference rule to the semi-intitionistic logic $\mathcal{SI}$.  
We then prove that the logic $\mathcal{DHMSH}$ is an implicative logic with respect to the defined connective $\to_H$, where $x \to y := x \to (x \land y)$.    
 In Section \ref{SecFour}, we prove 
the completeness theorem for the logic $\mathcal{DHMSH}$: The logic $\mathcal{DHMSH}$ is complete with respect to   
the variety $\mathbb{DHMSH}$ of dually hemimorphic semi-Heyting algebras. In Section  \ref{SecFive}, we deduce from Abstract Algebraic Logic that the logic $\mathcal{DHMSH}$ is algebraizable, in the sense of Blok and Pigozzi,  with the variety $\mathbb{DHMSH}$ as its equivalent algebraic semantics, 
 from which it follows that the lattice of axiomatic extensions of $\mathcal{DHMSH}$ is dually isomorphic to the lattice of subvarieties of $\mathbb{DHMSH}$.  These results enable us to present axiomatizations of several extensions of   
$\mathcal{DHMSH}$ by translating the (equational) axioms  of various (known) subvarieties of $\mathbb{DHMSH}$ from 
\cite{sankappanavar2011expansions, sankappanavar2014a, sankappanavar2014b, sankappanavar2016, sankappanavar2017} (see sections \ref{SecFive} and  \ref{S8}-\ref{Sec_JID}) 
into (propositional) axioms of the corresponding extensions. 
We also show that Moisil's ``logique modale''  ${\mathcal{LM}}$  
is equivalent to the logic $\mathcal{DMH}$ corresponding to the variety $\mathbb{DMH}$ of De Morgan Heyting algebras.  
 
 In Section  \ref{deducTheoExtensions}, we characterize the (axiomatic) extensions of 
$\mathcal{DHMSH}$ in which the ``Deduction Theorem'' holds. 
Sections  \ref{S7}-\ref{Sec_JID}     
deal with 
 applications of the results of Section \ref{SecFive} together with the algebraic results proved in \cite{sankappanavar2011expansions, sankappanavar2014a, sankappanavar2014b, sankappanavar2016, sankappanavar2017}. 
 More specifically, in Section  \ref{S7}, we present 
 axiomatizations for some extensions of the logic $\mathcal{DQDSH}$ 
 whose equivalent algebraic semantics are subvarieties of $\mathbb{DQDSH}$ generated by finitely many finite algebras, including two $2$-valued logics and  
twenty  $3$-valued logics and three 4-valued logics.  
Then we revisit the 3-valued {\L}ukasiewicz logic and give an alternative axiomatization for it.  In fact, we show that
the logic corresponding to the 3-element De Morgan Heyting algebra is equivalent to the 3-valued {\rm{\L}ukasiewicz} logic. 
Thereafter, we give axiomatizations for extensions of $\mathcal{DQDSH}$ corresponding to the subvarieties of the variety $\mathbb{DQDBSH}$ generated by dually quasi-De Morgan Boolean semi-Heyting algebras,  
completing the solution to PROBLEM 2 mentioned earlier.  We also give some extensions of the logic $\mathcal{DHMSH}$ which fail to possess the disjunction property.
Section \ref{S8}  describes some connections to Connexive Logic by showing that
some of these 2-valued, 3-valued and 4-valued logics are, in fact, connexive logics. 
Section \ref{section_infinite_chains} gives axiomatizations for De Morgan G\"odel logic 
and dually pseudocomplemented G\"odel logic corresponding to the varieties generated by the De Morgan Heyting chains and  
 the dually pseudocomplemented Heyting chains, respectively.  It also provides axiomatizations for the logics corresponding to their subvarieties.  
In Section \ref{Sec_Rdqdstsh}, we present axiomatizations for new logics corresponding to several subvarieties
of the variety RDQDStSH$_1$ of regular dually quasi-De Morgan Stone semi-Heyting algebras of level 1.
Section \ref{Sec_Rdmsh} presents axiomatizations for logics corresponding to a number of subvarieties of
RDMSH1 of regular De Morgan Stone semi-Heyting algebras of level 1, while Section \ref{Sec_JID} presents axiomatizations for logics corresponding to many subvarieties of
 JI-distributive linear semi-Heyting algebras of level 1.  
 Section \ref{section_concluding_remarks} concludes the paper with several open problems for future research.

\section{Preliminaries} \label{section_Preliminaries}

\indent A language $\mathbf L$ is a set of  finitary operations (or connectives), each with a fixed arity $n \geq 0$.   In this paper, we identify $\bot$ and $\top$ with $0$ and $1$ respectively and thus consider the languages $\langle \lor, \land, \to, \sim, 
 \bot, \top\rangle$ and $\langle \lor, \land, \to, ', 0, 1 \rangle$ as the same; however, we use the former in the context of logics and the latter in the context of algebras.    
For a countably infinite set $\mbox{\it Var}$ of propositional variables, the {\it formulas} of the language $\mathbf L$ are inductively defined as usual.  
A {\it logic}  (or, {\it a deductive system}) in the language $\mathbf L$ is a pair $\mathcal L =   \langle \bf{L}, \vdash_\mathcal{L} \rangle$, where 
$\vdash_\mathcal{L}$ is a substitution-invariant consequence relation on $\mbox{\it Fm}_{\,\bf L}$.   
We will present logics by means of their ``Hilbert style'' axioms and inference rules.

The set of formulas $\mbox{\it Fm}_{\,\bf L}$  can be turned into an algebra in the usual way. Throughout the paper, $\Gamma$ denotes a set of formulas and lower case Greek letters denote formulas.  The homomorphisms from the formula algebra ${\bf Fm}_{\bf L}$ into an $\mathbf L$-algebra (i.e, an algebra of type $\mathbf L$) $\mathbf A$ are called {\em interpretations} (or {\em valuations}) in A. The set of all such interpretations is denoted by  $Hom({\bf Fm}_{\bf L},\mathbf A)$. If $h \in Hom({\bf Fm}_{\bf L},\mathbf A)$ then the {\em interpretation of a formula} $\alpha$ under $h$ is its image $h \alpha \in A$, while  $h\Gamma$ denotes the set $\{h\phi \ | \ \phi \in \Gamma\}$.

\medskip

As mentioned earlier, Moisil presented in \cite{moisil1942} \rm(see also \cite{monteiro1980algebres}\rm), a propositional logic called, ``Logique modale''. 
We will refer to it as ``$\mathcal{LM}$''.

Moisil also introduced the variety of Heyting algebras endowed with an involution, in \cite{moisil1942a}, as the algebraic semantics for the logic $\mathcal{LM}$.   Monteiro \cite{monteiro1980algebres} investigated these algebras under the new name of symmetric Heyting algebras.  Among other things, he presented a proof of an algebraic completeness theorem for Moisil's calculus $\mathcal{LM}$ by showing that the logic $\mathcal{LM}$ is complete with respect to the variety of symmetric Heyting algebras. 

In the next section we will generalize Moisil's logic to a new logic called ``dually hemimorphic semi-Heyting logic''.   As a first step to achieve this goal, we need to present a generalization of intutionistic logic called ``Semi-intutionistic logic'' which was first introduced by the first author  in \cite{cornejo2011semi} in the language $\langle \vee, \land,  \to, \sim\rangle$. 
We will actually present below the more streamlined version of semi-Intuitionistic logic $\mathcal{SI}$ in the usual language 
$\langle \vee, \land, \to,  \bot, \top\rangle$, as first presented in \cite{cornejo15semiLogics} with the intuitionistic logic as an axiomatic extension.  To facilitate this presentation, it will be convenient to use $\alpha\to_H \beta$ as an abbreviation for $\alpha\to(\alpha\land \beta)$ so that the axioms given are easier to read.  Moreover, the connective $\to_H$ plays a crucial role in this section and in the sections that follow.

\begin{definition} {\rm  \cite{cornejo15semiLogics}}
The semi-intuitionistic logic $\mathcal{SI}$ (also called $\mathcal{SH}$) is defined in the language \\
$\{\vee, \land,  \to,  \bot, \top \}$ and it has the following   
axioms and the inference rule:\\

{\bf  \noindent  AXIOMS:}  

\noindent \rm \numeroAxioma{$\debaj{\alpha}{(\alpha \vee \beta)}$,} \label{axioma_supremo_izq}

\noindent \numeroAxioma{$\debaj{\beta}{(\alpha \vee \beta)}$,} \label{axioma_supremo_der}

\noindent \numeroAxioma{$\debaj{(\debaj{\alpha}{\gamma})} {  [\debaj{(\debaj{\beta}{\gamma})} {(\debaj{(\alpha \vee \beta)}{\gamma})}]}$,} \label{axioma_supremo_cota_inferior}

\noindent \numeroAxioma{$\debaj{(\alpha \wedge \beta)}{\alpha}$,} \label{axioma_infimo_izq}

\noindent \numeroAxioma{$\debaj{(\debaj{\gamma}{\alpha})} { [\debaj{(\debaj{\gamma}{\beta})}{(\debaj{\gamma}{(\alpha \wedge \beta)})}]}$,} \label{axioma_infimo_cota_superior}

\noindent \numeroAxioma{$\top$,} \label{axioma_top}

\noindent \numeroAxioma{$\debaj{\bot}{\alpha}$,} \label{axioma_bot}

\noindent \numeroAxioma{$\debaj{(\debaj{(\alpha \wedge \beta)}{\gamma})}{(\debaj{\alpha}{(\debaj{\beta}{\gamma})})}$,} \label{axioma_condicRes_InfAImplic}

\noindent \numeroAxioma{$\debaj{(\debaj{\alpha}{(\debaj{\beta}{\gamma})})}{(\debaj{(\alpha \wedge \beta)}{\gamma})}$,} \label{axioma_condicRes_ImplicAInf}

\noindent \numeroAxioma{$\debaj{(\debaj{\alpha}{\beta})}{(\debaj{(\debaj{\beta}{\alpha})}{(\debaj{(\alpha \to \gamma)}{(\beta \to \gamma)})})}$,} \label{axioma_BuenaDefImplic2} \label{axioma_ImplicADerecha}

\noindent \numeroAxioma{$\debaj{(\debaj{\alpha}{\beta})}{(\debaj{(\debaj{\beta}{\alpha})}{(\debaj{(\gamma \to \beta)}{(\gamma \to \alpha)})})}$.} \label{axioma_BuenaDefImplic1} \label{axioma_ImplicAIzquierda} \\

  {\bf RULE OF INFERENCE:}\\

(SMP):   From $\phi$ and $\phi \to_H {\gamma}$, deduce $\gamma$  ({\rm semi-Modus Ponens}).

\end{definition}

The following theorem and the lemma, proved in \cite{cornejo15semiLogics}, are useful in later sections.  

\begin{theorem} {\rm\cite{cornejo15semiLogics} (Completeness Theorem) \label{SIcompleteness}}
$$	\mbox{For all }\, \Gamma\cup\{\alpha\}\subseteq \mbox{\it Fm},\, 
 \Gamma \vdash_{\mathcal{SI}} \alpha \mbox{ if and only if }\, \Gamma \models_{\mathbb{SH}} \alpha.
 $$
\end{theorem}

\begin{lema} {\rm \cite{cornejo15semiLogics}} \label{lema_propiedades_SHB2}
	The following statements hold  in the logic $\mathcal{SI}$:
	\begin{enumerate}[{\rm (a)}]
		\item If $\Gamma \vdash_\mathcal{SI} \psi$ then $\Gamma \vdash_\mathcal{SI} \alpha \to_H \psi$,  \label{300513_15b} \label{saquemosAxioma1inciso01} 
		\item $\vdash_\mathcal{SI} \alpha \to_H \alpha$, \label{110613_12b}  \label{030613_18b} 
		\item $\vdash_\mathcal{SI} (\alpha \wedge \beta) \to_H \beta$, \label{190618_13b}  
		\item If $\vdash_\mathcal{SI} \alpha \to_H \beta$  then $\vdash_\mathcal{SI} (\alpha \wedge \gamma) \to_H (\beta \wedge \gamma)$ and $\vdash_\mathcal{SI} (\gamma \wedge \alpha) \to_H (\gamma \wedge \beta)$,   \label{030713_07} \label{saquemosAxioma1inciso05}
		\item $\vdash_\mathcal{SI} \debaj{(\debaj{\alpha}{\beta})}{[\debaj{(\debaj{\beta}{\gamma})}{(\debaj{\alpha}{\gamma})}]}$, \label{020713_42b}   
		\item  If $ \Gamma \vdash_{\mathcal{SI}} \alpha$  and $\Gamma \vdash_{\mathcal{SI}} \beta$ then $\Gamma \vdash_{\mathcal{SI}} \alpha \wedge \beta$.  \label{290413_infimos}
	\end{enumerate}
\end{lema}

We now give some useful information about semi-Heyting algebras. 

A key feature of semi-Heyting algebras is the following:\\
 Every semi-Heyting algebra $\langle A,\lor,\land,\to, 0,1 \rangle$ gives rise, in a natural way, to a Heyting algebra $\langle A,\lor,\land,\to_H, 0, 1 \rangle$, 
where $x\to_H y := x\to(x\land y)$, for $x,y \in A$ (see \cite{Abad2013semi}).\\

\begin{Lemma}{\rm\cite{Abad2013semi}}	 Let  $\mathbf{A} = \langle A, \lor, \land, \to, 0, 1\rangle$ be a semi-Heyting algebra and let
 $a,b,c \in A$.  Then, 
 \begin{itemize}
\item[{\rm(1)}] $a \land b \leq c$ if and only if  $a \leq b \to (b \land c)$,
\item[{\rm(2)}] $(a \to b) \wedge (b \to a) = 1$ if and only if $(a \to_H b) \wedge (b \to_H a) = 1$,
\item[{\rm(3)}] $a = b$ if and only if $(a \to_H b) \wedge (b \to_H a) = 1$, 
 \item[{\rm(4)}] $a \to b \leq a \to_H b$.    
\end{itemize}
\end{Lemma}

\begin{Lemma}{\rm \cite[Corollary 3.9]{cornejo2011semi}}
Let  $\mathbf{A} = \langle A, \lor, \land, \to, 0, 1\rangle$ be a semi-Heyting algebra and $a,b \in A$.  Then,  $a \to_H b = 1$ if and only if $a \leq b$. 
\end{Lemma}


\vspace{1cm}
\section{The logic $\mathcal{DHMSH}$: Axioms, Rules, and Rasiowa's Implicativeness} \label{section_logicDHMSH}
\

\indent In this section we present a new propositional logic called ``dually hemimorphic semi-Heyting logic'' 
denoted by $\mathcal{DHMSH}$ and prove, as a first step toward a completeness theorem, that 
the logic $\mathcal{DHMSH}$ is implicative in the sense of Rasiowa with respect to the implication $\to_H$.

\begin{definition} 
	The dually hemimorphic semi-Heyting logic, $\mathcal{DHMSH}$,  is defined in the language 
$\langle \vee,  \land, \to, \sim,  \bot, \top \rangle$ and it has the following axioms and rules of inference:\\
	 
{\bf \noindent AXIOMS:}\\

{\rm (A1), (A2), $\dots$, (A11)} of the logic $\mathcal{SI}$, together with the following three additional axioms:

 \rm\numeroAxioma{$\top \to_H \sim\bot,$}  \label{axiom_bot_neg}

 \rm \numeroAxioma{$\sim\top \to_H \bot,$} \label{axiom_top_neg}

 \numeroAxioma{$\sim(\alpha \wedge \beta) \to_H (\sim\alpha \ \vee \sim\beta).$} \label{axiom_DM_law1}\\


{\bf \noindent RULES OF INFERENCE:}\\

	       {\bf (SMP)}  \quad From $\phi$ and $\phi \to_H {\gamma}$, deduce $\gamma$  ({\rm semi-Modus Ponens}),\\
 
              {\bf (SCP)} \quad From  $\phi \to_H \gamma$, deduce $\sim\gamma \to_H \ \sim\phi$ \rm(semi contraposition rule\rm).
\end{definition}

Since the axioms and the inference rule of the logic $\mathcal{SI}$ are included in the logic $\mathcal{DHMSH}$, the following result is immediate.

\begin{theorem} \label{teorema_validity_SI_to_DMSH}
	Let  $\Gamma\cup\{\alpha\}\subseteq \mbox{\it Fm}$. If $\Gamma \vdash_{\mathcal{SI}} \alpha$ then $\Gamma \vdash_{\mathcal{DHMSH}} \alpha$.
\end{theorem}

The following lemma is needed later.
\begin{lema}  \label{lema_propiedades_SI}
	Let $\Gamma\cup\{\alpha, \beta, \gamma, \psi\}\subseteq \mbox{\it Fm}$. The following statements hold  in the logic $\mathcal{DHMSH}$:
	\begin{enumerate}[{\rm (a)}]
		\item If $\Gamma \vdash_\mathcal{DHMSH} \psi$, then $\Gamma \vdash_\mathcal{DHMSH} \alpha \to_H \psi$,  \label{300513_15} 
		\item $\vdash_\mathcal{DHMSH} \alpha \to_H \alpha $, \label{110613_12}  \label{030613_18}
		\item $\vdash_\mathcal{DHMSH} (\alpha \wedge \beta) \to_H \beta$, \label{190618_13} 
		\item $\Gamma \vdash_\mathcal{DHMSH} \debaj{(\debaj{\alpha}{\beta})}{[\debaj{(\debaj{\beta}{\gamma})}{(\debaj{\alpha}{\gamma})}]}$, \label{020713_42} 
		\item $\Gamma \vdash_\mathcal{DHMSH} \alpha \wedge \beta$ if and only if $\Gamma \vdash_\mathcal{DHMSH} \alpha$ and $\Gamma \vdash_\mathcal{DHMSH} \beta$, \label{300317_01}
		\item  If $\Gamma \vdash_\mathcal{DHMSH} \alpha \to_H \beta$, then $\Gamma \vdash_\mathcal{DHMSH} (\alpha \wedge \gamma) \to_H (\beta \wedge \gamma)$ and $\Gamma \vdash_\mathcal{DHMSH} (\gamma \wedge \alpha) \to_H (\gamma \wedge \beta)$. \label{190917_08}
	\end{enumerate}
\end{lema}

\begin{Proof}
	Items (\ref{300513_15}), (\ref{030613_18}), (\ref{190618_13}) and (\ref{020713_42}) follow from Theorem \ref{teorema_validity_SI_to_DMSH} and items (\ref{300513_15b}), (\ref{110613_12b}), (\ref{190618_13b}) and (\ref{020713_42b}) of Lemma \ref{lema_propiedades_SHB2}, respectively.
We have, by Theorem \ref{teorema_validity_SI_to_DMSH} and Lemma \ref{lema_propiedades_SHB2} (\ref{290413_infimos}),  that if $\Gamma \vdash_\mathcal{DMSH} \alpha$ and $\Gamma \vdash_\mathcal{DHMSH} \beta$ then $\Gamma \vdash_\mathcal{DHMSH} \alpha \wedge \beta$. The other half of the item (\ref{300317_01}) follows easily from axiom (A\ref{axioma_infimo_izq}), item (\ref{190618_13}) and (SMP).  Finally, Item (\ref{190917_08}) follows from Lemma \ref{lema_propiedades_SHB2} (\ref{030713_07}).
\end{Proof}

\begin{Lemma} \label{lemma_properties_DMSH_logic} \label{lemma_properties_DHMSH_logic}
	Let $\Gamma\cup\{\alpha, \beta, \gamma\}\subseteq \mbox{\it Fm}$. 
	\begin{enumerate}[{\rm(a)}]
		\item  If $\Gamma \vdash_\mathcal{DHMSH} \alpha \to_H \beta$ and  $\Gamma \vdash_\mathcal{DHMSH} \beta \to_H \gamma$ then $\Gamma \vdash_\mathcal{DHMSH} \alpha \to_H \gamma$. \label{190917_01}
		\item $\Gamma, \beta \to_H \alpha \vdash_\mathcal{DHMSH} \ \sim\alpha \to_H \sim\beta.$ \label{IL3_partial}
		\item $\Gamma \vdash_\mathcal{DHMSH} \sim(\alpha \vee \beta) \to_H (\sim\alpha \wedge \sim\beta).$\label{241016_08}
		\item  If $\Gamma \vdash_\mathcal{DHMSH} \alpha \to_H \beta$, then \\
		\quad $\Gamma \vdash_\mathcal{DHMSH} (\alpha \lor \gamma) \to_H (\beta \lor \gamma)$ \ and \ $\Gamma \vdash_\mathcal{DHMSH} (\gamma \lor \alpha) \to_H (\gamma \lor \beta)$. \label{190917_11}
		\item $\Gamma \vdash_\mathcal{DHMSH} (\sim\alpha \vee \sim\beta) \to_H \ \sim(\alpha \wedge \beta)$. \label{171117_01}
	\end{enumerate}
\end{Lemma}

\begin{Proof}
	\begin{enumerate}[(a)]
		\item 
		This follows from \ref{lema_propiedades_SI} (\ref{020713_42}), 
		using (SMP).
		
		\item 
		This is immediate from (SCP).

		\item 		
		Using (A1), (A2) and (SCP), we get
		\begin{equation} \label{eq1}
		\Gamma \vdash_\mathcal{DHMSH} \sim(\alpha \vee \beta) \to_H \sim\alpha, \text{and }
		\Gamma \vdash_\mathcal{DHMSH} \ \sim(\alpha \vee \beta) \to_H \ \sim\beta.
		\end{equation}
\noindent The conclusion follows from \eqref{eq1}, (A5) and (SMP).
		
		\item 
		\begin{enumerate}[1.]
			\item $\Gamma \vdash_\mathcal{DHMSH} \alpha \to_H \beta$ by hypothesis. \label{190917_12} \label{190917_19}
			\item $\Gamma \vdash_\mathcal{DHMSH} \beta \to_H (\beta \lor \gamma)$ by (A\ref{axioma_supremo_izq}) \label{190917_13}
			\item $\Gamma \vdash_\mathcal{DHMSH} \alpha \to_H (\beta \lor \gamma)$  by  (\ref{190917_01}) in \ref{190917_12} and \ref{190917_13}. \label{190917_15}
			\item $\Gamma \vdash_\mathcal{DHMSH} \gamma \to_H(\beta \lor \gamma) $ by (A\ref{axioma_supremo_der}) \label{190917_17}
			\item $\Gamma \vdash_\mathcal{DHMSH} \debaj{(\debaj{\alpha}{(\beta \lor \gamma)})} {  [\debaj{(\debaj{\gamma}{(\beta \lor \gamma)})} {(\debaj{(\alpha \vee \gamma)}{(\beta \lor \gamma)})}]}$ by (A\ref{axioma_supremo_cota_inferior}) \label{190917_16}
			\item $\Gamma \vdash_\mathcal{DHMSH}   \debaj{(\debaj{\gamma}{(\beta \lor \gamma)})} {(\debaj{(\alpha \vee \gamma)}{(\beta \lor \gamma)})}$  by (SMP) in \ref{190917_15} and \ref{190917_16}. \label{190917_18}
			\item $\Gamma \vdash_\mathcal{DHMSH} \debaj{(\alpha \vee \gamma)}{(\beta \lor \gamma)}$  by (SMP) in \ref{190917_17} and \ref{190917_18}.
			\item $\Gamma \vdash_\mathcal{DHMSH} \gamma \to_H (\gamma \lor \beta)$ by (A\ref{axioma_supremo_izq})  \label{190917_21}
			\item $\Gamma \vdash_\mathcal{DHMSH} \beta \to_H (\gamma \lor \beta)$ by (A\ref{axioma_supremo_der}) \label{190917_20}
			\item $\Gamma \vdash_\mathcal{DHMSH} \alpha \to_H (\gamma \lor \beta)$  by (\ref{190917_01}) and (SMP) in \ref{190917_19} and \ref{190917_20}. \label{190917_23}
			\item $\Gamma \vdash_\mathcal{DHMSH} \debaj{(\debaj{\gamma}{(\gamma \lor \beta)})} {  [\debaj{(\debaj{\alpha}{(\gamma \lor \beta)})} {(\debaj{(\gamma \vee \alpha)}{(\gamma \lor \beta)})}]}$ by (A\ref{axioma_supremo_cota_inferior}) \label{190917_22}
			\item $\Gamma \vdash_\mathcal{DHMSH}   \debaj{(\debaj{\alpha}{(\gamma \lor \beta)})} {(\debaj{(\gamma \vee \alpha)}{(\gamma \lor \beta)})}$  by (SMP) in \ref{190917_21} and \ref{190917_22}. \label{190917_24}
			\item $\Gamma \vdash_\mathcal{DHMSH} \debaj{(\gamma \vee \alpha)}{(\gamma \lor \beta)}$  by (SMP) in \ref{190917_23} and \ref{190917_24}.
		\end{enumerate}	
		
		\item 
		\begin{enumerate}[1.]
			\item $\Gamma \vdash_\mathcal{DHMSH} (\alpha \wedge \beta) \to_H \alpha$ by axiom (A\ref{axioma_infimo_izq}).
			\item $\Gamma \vdash_\mathcal{DHMSH} \ \sim\alpha \to_H \ \sim(\alpha \wedge \beta)$ by (SCP). \label{161117_01}
			\item $\Gamma \vdash_\mathcal{DHMSH} (\alpha \wedge \beta) \to_H \beta$ by \ref{lema_propiedades_SI} (\ref{190618_13}).
			\item $\Gamma \vdash_\mathcal{DHMSH} \ \sim\beta \to_H \ \sim(\alpha \wedge \beta)$ by (SCP). \label{161117_03}
			\item $\Gamma \vdash_\mathcal{DHMSH} \debaj{(\debaj{\sim\alpha}{\ \sim(\alpha \wedge \beta)})} {  [\debaj{(\debaj{\sim\beta}{\ \sim(\alpha \wedge \beta)})} {(\debaj{(\ \sim\alpha \vee \sim\beta)}{\ \sim(\alpha \wedge \beta)})}]}$  by axiom (A\ref{axioma_supremo_cota_inferior}). \label{161117_02}
			\item $\Gamma \vdash_\mathcal{DHMSH}   \debaj{(\debaj{ \ \sim\beta}{\ \sim(\alpha \wedge \beta)})} {(\debaj{(\ \sim\alpha \vee \sim\beta)}{\ \sim(\alpha \wedge \beta)})}$  by (SMP) in \ref{161117_01}, \ref{161117_02}. \label{161117_04}
			\item $\Gamma \vdash_\mathcal{DHMSH}   \debaj{(\sim\alpha \vee \sim\beta)}{\ \sim(\alpha \wedge \beta)}$  by (SMP) in \ref{161117_03} and \ref{161117_04},
		\end{enumerate}			
	\end{enumerate}
proving the lemma.	
\end{Proof}

\subsection{$\mathcal{DHMSH}$ as an implicative logic in the sense of Rasiowa}
\ \\

 In 1974, Rasiowa (\cite[page 179]{rasiowa1974algebraic}) introduced an important class of logics called  ``standard systems of implicative extensional propositional calculus'' and associated a class of algebras with each of them, by a generalization of the classical Lindenbaum-Tarski process.  We will refer to these logics as ``implicative logics in the sense of Rasiowa'' (``implicative logics'', for short).  These logics have played a pivotal role in the development of Abstract Algebraic Logic.  We now recall the definition  
of implicative logics.  We follow Font \cite{Fo17}.

\begin{definicion} {\rm \cite{rasiowa1974algebraic}}  
	Let $\mathcal L$  be a logic in a language $\mathbf{L}$ that includes a binary connective $\to$, either primitive or defined by a term in exactly two variables. Then $\mathcal L$ is called 
	 an implicative logic 
	 with respect to the binary connective $\to$, if the following conditions are satisfied:
	\begin{itemize}
		\item[\rm (IL1)] $\vdash_{\mathcal L} \alpha \to \alpha$, 
		\item[\rm (IL2)] $\alpha \to \beta, \beta\to \gamma \vdash_{\mathcal L} \alpha \to \gamma$,
		\item[\rm (IL3)] For each symbol $f \in \mathbf L$ of arity $n \geq 1$,\\ $\left\{\begin{array}{c}
		\alpha_1 \to \beta_1, \ldots, \alpha_n \to \beta_n, \\
		\beta_1 \to \alpha_1, \ldots, \beta_n \to \alpha_n
		\end{array}
		\right\} \vdash_{\mathcal L} f(\alpha_1, \ldots, \alpha_n) \to f(\beta_1, \ldots, \beta_n)$,
		\item[\rm (IL4)] $\alpha, \alpha \to \beta \vdash_{\mathcal L} \beta$,
		\item[\rm (IL5)] $\alpha \vdash_{\mathcal L} \beta \to \alpha$.
	\end{itemize}
\end{definicion}

The following lemma was proved in  \cite[Lemma 4.6]{cornejo15semiLogics}.

\begin{Lemma}{\rm \cite[Lemma 4.6]{cornejo15semiLogics}} \label{teorema_SI_implicativa}
	The logic $\mathcal{SI}$ is implicative with respect to the connective $\to_H$.
\end{Lemma}

The following theorem follows from Theorem \ref{teorema_validity_SI_to_DMSH}, Lemma \ref{lemma_properties_DMSH_logic} (\ref{IL3_partial}) and Lemma \ref{teorema_SI_implicativa}.

\begin{theorem} \label{teorema_DMSH_implicativa} \label{teorema_DHMSH_implicativa}
	The logic $\mathcal{DHMSH}$ is implicative with respect to the connective $\to_H$.
\end{theorem}

\bigskip
\section{Completeness of $\mathcal{DHMSH}$} \label{SecFour}
\
Let $\mathbf L$ denote a language.  Identities in $\mathbf L$ are ordered pairs of $\mathbf L$-formulas that will be written in the form $\alpha \approx \beta$.
An interpretation $h$ in $\mathbf A$ satisfies an identity $\alpha \approx \beta$ if $h \alpha = h \beta$.  
We denote this satisfaction relation by the notation: $\mathbf{A} \models_h
\alpha \approx \beta $. An algebra $\mathbf A$ {\it satisfies the equation} $\alpha \approx \beta$ if all the interpretations in $\mathbf A$ satisfy it; in symbols,
$$\mathbf{A} \models \alpha \approx \beta \mbox{ if and only if } \mathbf{A} \models_h
\alpha \approx \beta, \mbox{ for all } h \in Hom({\bf Fm}_{\bf L},\mathbf A).$$
A class $\mathbb K$ of algebras  {\it satisfies the identity} $\alpha \approx \beta$ when all the algebras in $\mathbb K$ satisfy it; i.e.
$$\mathbb{K} \models \alpha \approx \beta \mbox{ if and only if } \mathbf{A}\models \alpha \approx \beta, \mbox{ for all }\mathbf A \in \mathbb K.$$ 

If $\bar{x}$ is a sequence of variables and $h$ is an interpretation in $\mathbf{A}$, then we write $\bar{a}$ for $h(\bar{x})$.
For a class $\mathbb K$ of $\mathbf L$-algebras, 
we define the relation $\models_{\mathbb K}$ that holds between a set $\Delta$ of identities and a single identity $\alpha \approx  \beta$ as follows: \\ 
$$\Delta \models_{\mathbb K} \alpha \approx \beta \text{ if and only if }$$

for every $\mathbf A \in K$ and every interpretation $\bar{a}$ of the variables of $\Delta \cup  \{\alpha \approx \beta\}$ in $\mathbf{A}$, 

$\phi^{\mathbf A}(\bar{a} ) = \psi^{\mathbf A} (\bar{a}) $,  for every $\phi \approx \psi \in \Delta \quad \Rightarrow  \quad \alpha^{\mathbf A}(\bar{a}) = \beta^{\mathbf A}(\bar{a})$.\\

In this case, we say that $\alpha \approx \beta$ is a $\mathbb K$-consequence of $\Delta$. 
The relation $\models_K$ is called the {\it semantic equational consequence relation} determined by K.\\ 

\indent Our goal in this section is to prove that the logic $\mathcal{DHMSH}$ is complete with respect to the variety $\mathbb{DHMSH}$.   For this we need the following definition from \cite{rasiowa1974algebraic}.

\begin{definicion} {\rm \cite[Definition 6, page 181]{rasiowa1974algebraic}}, {\rm \cite[Definition 2.5]{Fo17}}  \label{definicion_lalg} 
	Let $\mathcal L$ be an implicative logic in the language $\mathbf L$  with an implication connective $\to$.  An $\mathcal L$-algebra is an algebra $\mathbf A$ in the language $\mathbf L$ that has an element $1$ with the following properties:
	\begin{quote}
		\begin{itemize}
			\item[{\rm (LALG1)}] For all $\Gamma \cup \{\phi\} \subseteq \mbox{\it Fm} $ and all $h \in Hom({\bf Fm}_{\bf L}, \mathbf A)$,\\
			 if $\Gamma \vdash_{\mathcal L} \phi$ and $h \Gamma \subseteq \{1\}$ then $h \phi = 1$, 
			\item[{\rm (LALG2)}] For all $a,b \in A$, if $a \to b = 1$ and $b \to a = 1$ then $a=b$.
		\end{itemize}
	\end{quote} 
	
\noindent	The class of $\mathcal L$-algebras is denoted by $\mathbf{Alg^*\mathcal{L}}$. 
\end{definicion}

We also need the following result from \cite{cornejo15semiLogics}.  

\begin{Theorem} {\rm \cite[Corollary 4.8]{cornejo15semiLogics}} \label{theorem_algSH_SH}
	$\mathbf{Alg^*{\mathcal{SI}}} = \mathbb{SH}$.
\end{Theorem}

Since $\mathcal{DHMSH}$ is an implicative logic with respect to the binary connective $\to_H$ by Theorem \ref{teorema_DMSH_implicativa}, we obtain the following result, in view of \cite[Theorem 7.1, pag 222]{rasiowa1974algebraic}.  

\begin{theorem}\label{lema_SI_SIalgebras}\label{lagSI_igual_SH}
	The logic $\mathcal{DHMSH}$ is complete with respect to the class 
	$\sf{A}lg^*{\mathcal{DHMSH}}$. 
	In other words,	
	$$\mbox{for all } \Gamma \cup \{\phi\} \subseteq \mbox{\it Fm, } 
	\Gamma \vdash_{\mathcal{DHMSH}}\phi \mbox{ if and only if }
	\Gamma \models_{\mathbb{DHMSH}} \phi.$$			
\end{theorem}

As a last step to complete the proof of algebraic completeness of the logic $\mathcal{DHMSH}$, we need to prove that $\mathbf{Alg^*{\mathcal{DHMSH}}} = \mathbb{DHMSH}$. 

\begin{Lemma} \label{lemma_DMSH_included_LAlg}
	$\mathbb{DHMSH} = \sf{A}lg^*{\mathcal{DHMSH}}$.
\end{Lemma}

\begin{Proof} 
First, we wish to prove that $\mathbb{DHMSH} \subseteq \sf{A}lg^*{\mathcal{DHMSH}}$.
	Let $\mathbf A \in \mathbb{DHMSH}$, $\Gamma \cup \{\phi\} \subseteq \mbox{\it Fm} $ and $h \in Hom(\mbox{\it Fm}_{\bf L}, \mathbf A)$ such that $\Gamma \vdash_{\mathcal{DMSH}} \phi$
	and $h \Gamma \subseteq \{1\}$.  We need to verify that $h \phi = 1$. We will proceed by induction on the length of the proof of $\Gamma \vdash_{\mathcal{DMSH}} \phi$.
	
	\begin{itemize}
		\item Assume that $\phi$ is an axiom.
		
		If $\phi$ is one of the axioms (A\ref{axioma_supremo_izq}) to (A\ref{axioma_ImplicAIzquierda}) then $\vdash_{\mathcal{SI}} \phi$.  Hence, by theorem \ref{SIcompleteness}, 
		$\models_{\mathbb{DHMSH}} \phi$ and so, $h(\phi) = \top$.

		If $\phi$ is the axiom (A\ref{axiom_bot_neg}) then, using (E\ref{conditionSD2d1}), we have  $h(\phi) = h(\top \to_H \ \sim\bot) = 1 \to_H 0' 
		= 1$. 
		
		If $\phi$ is the axiom (A\ref{axiom_top_neg}) then, using (E\ref{conditionSD2d2}), we get that $h(\phi) = h(\sim\top \to_H \bot) = 0 \to_H 0 = 1.$                             
		
		If $\phi$ is the axiom (A\ref{axiom_DM_law1}) then, using (E\ref{conditionSD3d}), we obtain that $h(\phi) = h(\sim(\alpha \wedge \beta) \to_H (\ \sim\alpha \vee \sim\beta)) = (h(\alpha) \wedge h(\beta))' \to_H (h(\alpha)' \vee h(\beta)') = (h(\alpha) \wedge h(\beta))' \to_H (h(\alpha) \wedge h(\beta))' = 1$.
		
		\item If $\phi \in \Gamma$ then $h(\phi) = \top$ by hypothesis.
		
		\item Assume now that $\Gamma \vdash_{\mathcal L} \phi$ is obtained from an application of (SMP). Then there exist a formula $\psi$ such that $\Gamma \vdash_{\mathcal L} \psi$ and $\Gamma \vdash_{\mathcal L} \psi \to_H \phi$. By induction, $h(\psi) = 1$ and $h(\psi \to_H \phi) = 1$. Then 
$1= h(\psi) \to_H h(\phi) = 1 \to_H h(\phi)= h(\phi).$
	
		\item Assume that $\Gamma \vdash_{\mathcal L} \phi$ is the result of an application of the rule (SCP). Then for $\alpha, \beta \in \mbox{\it Fm}$, $\phi = \sim\beta \to_H \sim\alpha$ and $\Gamma \vdash_{\mathcal L} \alpha \to_H \beta$.  By induction, $1 = h(\alpha \to_H \beta) = h(\alpha) \to_H h(\beta)$ and, consequently $h(\alpha) \leq h(\beta)$. 
		Then, using condition (E\ref{conditionSD3d}), $h(\beta)' \leq h(\alpha)'$. Hence $h(\beta)' \to_H h(\alpha)' = 1$. Therefore $h(\phi) = h(\sim\beta \to_H \sim\alpha) = h(\beta)' \to_H h(\alpha)' = 1$.
	\end{itemize}
	Hence, the induction is complete and so, we concludes that $\mathbf A$ satisfies (LALG1). It is easy to see that the condition (LALG2) also holds, implying $\mathbf A  \in \sf{A}lg^*{\mathcal{DHMSH}}$.  

Next, we prove the other inclusion.
Let $\mathbf A = \langle A, \vee, \wedge, \to, ', 0, 1\rangle \in \sf{A}lg^*{\mathcal{DHMSH}}$.  Notice that $\langle A, \lor, \wedge,  \to, 0,1\rangle \in \sf{A}lg^*{\mathcal{\mathcal{SI}}}$.  By Theorem \ref{theorem_algSH_SH}, $\langle A, \lor, \wedge ,  \to, 0,1 \rangle \in \mathbb{SH}$.  Now, it only remains to show that $\mathbf A$ satisfies the conditions (E\ref{conditionSD2d1}) to (E\ref{conditionSD3d}).
	
	In view of axiom $(A\ref{axiom_bot_neg})$ and (LALG1), we have that  $\mathbf A \models 1 \to_H 0' \approx 1$. Using (LALG1) and Lemma \ref{lema_propiedades_SI} (\ref{300513_15}), we get               
	$\mathbf A \models 0' \to_H 1 \approx 1$. Then by (LALG2), $\mathbf A \models 1 \approx 0'$.
	In view of axioms $(A\ref{axioma_bot})$ and $(A\ref{axiom_top_neg})$, together with (LALG1), we have that $\mathbf A \models 0 \to_H 1' \approx 1$ and  $\mathbf A \models 1' \to_H 0 \approx 1$. Then by (LALG2), $\mathbf A \models 1' \approx 0$.
	
	By Lemma \ref{lemma_properties_DMSH_logic} (\ref{171117_01}) and the condition (LALG1), $\mathbf A$ satisfies the identity $(x' \vee y') \to_H (x \wedge y)' \approx 1$.  Also, In view of axiom $(A\ref{axiom_DM_law1})$, and (LALG1), $\mathbf A$ satisfies the identity $(x \vee y)' \to_H (x' \wedge  y') \approx 1$. Applying (LALG2), the algebra satisfies (E\ref{conditionSD3d}).
	Consequently $\mathbf A \in \mathbb{DHMSH}$.
\end{Proof}

\medskip

\, \,
We are now ready to present the completeness theorem for the logic $\mathcal{DHMSH}$.

\begin{theorem}\label{completeness_DMSH} \label{completeness_DHMSH}
	The logic $\mathcal{DHMSH}$ is complete with respect to the variety
	$\mathbb{DHMSH}$. 
\end{theorem}
\begin{Proof}  From Lemma \ref{lemma_DMSH_included_LAlg} 
we have $\sf{A}lg^*{\mathcal{DHMSH}} = \mathbb{DHMSH}$.  The conclusion follows from Theorem \ref{lagSI_igual_SH}. 
\end{Proof}

\bigskip

\section{Equivalent Algebraic Semantics and Axiomatic Extensions of the logic $\mathcal{DHMSH}$ } \label{SecFive} 
\

 Our goal in this section is to improve Theorem \ref{completeness_DHMSH} by proving that 
the logic $\mathcal{DHMSH}$ is algebraizable and $\mathbb{DHMSH}$ is an equivalent algebraic semantics of the logic $\mathcal{DHMSH}$.

Here we first recall some relevant notions and results from Abstract Algebraic Logic 
(see \cite[Section 2.1]{BP89}, \cite{Font03survey}, or \cite{Fo17}).  

\begin{definition} {\rm \cite[Definition 2.2]{BP89} {\rm \cite[Definition 3.4]{Fo17}}}
Let $\langle {\bf L}, \vdash_{\mathcal L}\rangle$ be a logic (i.e., deductive system) and $\mathbb K$ a class of {\bf L}-algebras.   $\mathbb K$ is called an ``algebraic semantics'' for $\langle {\bf L}, \vdash_{\mathcal L}\rangle$  if \ $\vdash_{\mathcal L}$ can be interpreted in $\models_\mathbb K$ in the following sense:\\

 There exists a finite set $\delta_i(p) \approx \epsilon_i(p)$, for $i < n$, 
 of identities with a single variable $p$ such that, for all $\Gamma \cup {\phi} \subseteq  Fm$ and each $j < n$,

$${\rm(A)} \qquad \Gamma \vdash_{\mathcal L} \phi \ \Leftrightarrow \ \{\delta_i[\psi/p] \approx \epsilon_i[\psi/p] : i < n,  \psi  \in \Gamma\} \models_K  \delta_j[\phi/p] \approx \epsilon_j[\phi/p],$$
where $\delta[\psi/p] $ denotes the formula obtained by the substitution of $\psi$ at every occurrences of $p$ in $\delta$.
The identities $\delta_i \approx  \epsilon_i$, for $i < n$, are called ``defining identities'' for 
$\langle {\bf L}, \vdash_{\mathcal L}\rangle$ and $\mathbb K$.
\end{definition}

In what follows, we will use ``$\Gamma \ismodeledby\models_{\mathbb K} \Delta$'' as an abbreviation for ``$\Gamma \models_{\mathbb K} \Delta$ and $\Delta \models_{\mathbb K} \Gamma$.''   Also, $\delta(\Delta(\phi, \psi))$ denotes the formula obtained by the substitution of the formula $\Delta(\phi, \psi)$ at every occurrence of $p$ in $\delta(p)$.

\begin{definition} {\rm \cite[Definition 2.8]{BP89}}, {\rm \cite[Definition 3.11]{Fo17}} \label{def2.8}
Let $\langle {\bf L}, \vdash_{\mathcal L}\rangle$ be a logic and $\mathbb K$ an algebraic semantics for 
$\langle {\bf L}, \vdash_{\mathcal L}\rangle$  with defining identities $\delta_i =  \epsilon_i$, for $i < n$.\\ 
$\mathbb K$ is said to be ``equivalent to'' $\langle {\bf L}, \vdash_{\mathcal L}\rangle$ if there exists a finite set $\Delta_j(p,q)$, for $j < m$, of formulas with two variables $p, q$ such that 

for every identity $\phi  \approx \psi$, for $i<n$, and for $j<m$, 
 $${\rm(E)} \qquad   \phi  \approx \psi \ismodeledby \models_{\mathbb K}  \{\delta_i(\Delta_j(\phi, \psi)) \approx \epsilon_i(\Delta_j(\phi, \psi)) : i<n, j<m\}.$$
 
 The set $\Delta_j(p,q)$, $j < m$, of formulas with two variables, satisfying (E) is called a set of ``equivalence formulas'' for $\langle {\bf L}, \vdash_{\mathcal L}\rangle$ and $\mathbb K$.
 A logic $\mathcal L$ is said to be ``algebraizable'' if and only if it has an equivalent algebraic semantics $\mathbb K$.

\end{definition}

The following theorem, proved in \cite{BP89}, is crucial in what follows. 

\begin{Theorem} {\rm(\cite{BP89}, \cite[Proposition 3.15]{Fo17}}) \label{EquivalentSemanticsTheorem}
Every implicative logic  $\mathcal L$ is algebraizable with respect to the class $Alg^*\mathcal L$ and the algebraizability is witnessed by the defining identity $p \approx p \to p$ and the equivalence formulas 
$\Delta =\{p \to q, q \to p\}$. 
\end{Theorem}

 As an immediate consequence of Theorem \ref{EquivalentSemanticsTheorem}, Theorem \ref{teorema_DHMSH_implicativa} and Theorem \ref{completeness_DHMSH}, we obtain the following
crucial result. 
 
\begin{Corollary} \label{CorExt} 
 The logic $\mathcal{DHMSH}$ is algebraizable, and  the variety  
 $\mathbb{DHMSH}$ is the equivalent algebraic semantics for $\mathcal{DHMSH}$ with the defining identity $p \approx p \to_H p$ \rm(equivalently, $p \approx 1$) and the equivalence formulas $\Delta = \{p \to_H q, q \to_H p \}$.
 \end{Corollary}

Also, we just mention, in passing, the following fact  which follows from Theorem \ref{EquivalentSemanticsTheorem}, Lemma \ref{teorema_SI_implicativa} and Theorem \ref{theorem_algSH_SH} about the semi-intuitionistic logic $\mathcal{SH}$.

 \begin{Corollary} \label{CorExtSH} 
 The logic $\mathcal{SI}$ is algebraizable, and the variety 
 $\mathbb{SH}$ is the equivalent algebraic semantics for $\mathcal{SI}$ with the defining identity $p \approx p \to_H p$ and the equivalence formulas $\Delta = \{p \to_H q, q \to_H p \}.$\\
 \end{Corollary}

\subsection{Axiomatic Extensions of $\mathcal{DHMSH}$ }

\
\medskip

A logic $\mathcal L'$ is an {\it axiomatic extension} of $\mathcal L$ if $\mathcal L'$ is obtained by adjoining new axioms but keeping the rules of inference the same as in $\mathcal L$. 
\indent In the sequel, we sometimes use the term ``extension'' for ``axiomatic extension''. 	
Let $Ext(\mathcal L)$ denote the lattice of axiomatic extensions of the logic $\mathcal L$ and 
$\mathbf{L_V(\mathbb K)}$ denote the lattice of subvarieties of the variety $\mathbb K$ of algebras. 

The following theorem is one of the hallmark accomplishments of Abstract Algebraic Logic.  
	
\begin{Theorem} {\rm \cite[Theorem 3.33] {Fo17}} \label{alg-logic-theorem}
Let $\mathcal L$ be an algebraizable logic with the equivalent algebraic semantics $\mathbb K$.  Then
$Ext(\mathcal L)$ is dually isomorphic to $\mathbf{L_V(\mathbb K)}$. 
\end{Theorem}

 The following theorem is a consequence of Theorem \ref{alg-logic-theorem},   
 Theorem \ref{completeness_DHMSH} and 
 Corollary \ref{CorExt}.

\begin{Theorem} {\rm (Isomorphism Theorem for $\mathcal{DHMSH}$)} \label{Ext_subvariety_Theorem} 
$Ext(\mathcal{DHMSH})$ is dually isomorphic to $\mathbf{L_V(\mathbb{DHMSH})}$.
\end{Theorem}  
In a similar fashion, the following result is a consequence of Theorem \ref{alg-logic-theorem}, 
 Theorem \ref{SIcompleteness} and 
 Corollary \ref{CorExtSH}.  
\begin{Theorem} {\rm (Isomorphism Theorem for $\mathcal{SH}$)} \label{ExtSH_subvariety_Theorem} 
$Ext(\mathcal{SH})$ is dually isomorphic to $\mathbf{L_V(\mathbb{SH})}$.
\end{Theorem} 
The following theorem is an immediate consequence of Theorem \ref{Ext_subvariety_Theorem} and plays an important role in the rest of the paper.
Let $\mathbf{Mod}(\mathcal E) := \{\mathbf{A} \in \mathbb{DHMSH}: \mathbf{A} \models \delta \approx 1, \mbox{ for every } \delta \in \mathcal E\} $.  

 \begin{Theorem} \label{completeness_DHMSH_extension} \label{teo_040417_01}
 Let $\mathcal E$ be an extension of the logic $\mathcal{DHMSH}$.  Then\\
 \rm(a)  $\mathcal E$ is algebraizable with the same equivalence formulas and defining equations as those of the logic $\mathcal{DHMSH}$. \\ 
 (b)	$\mathbf{Mod}(\mathcal E)$ is the equivalent algebraic semantics for $\mathcal E$.
 \end{Theorem}

Note that Theorem \ref{ExtSH_subvariety_Theorem}  justifies the use of the phrase ``the logic corresponding to the subvariety $\mathbb  V$'' of $\mathbb{DHMSH}$.

\indent In the later sections of the paper, 
we give 
various applications of the results proved above and the algebraic results proved in \cite{sankappanavar2011expansions, sankappanavar2014a, sankappanavar2014b, sankappanavar2016, sankappanavar2017}.
We now give Hilbert-style axiomatizations for several important extensions of the logic $\mathcal{DHMSH}$.

To facilitate the presentation of the extensions of the logic $\mathcal{DHMSH}$, we first list several  important subvarieties of the variety $\mathbb{DHMSH}$ of dually hemimorphic semi-Heyting algebras that were introduced (or implicit) in \cite{sankappanavar2011expansions}.\\

{\bf LIST 1: SOME IMPORTANT SUBVARIETIES OF THE VARIETY $\mathbb{DHMSH}$} \label{listOfSubvarieties}

\begin{enumerate}
\item $\mathbb{DHMH}:$ {\em Dually hemimorphic Heyting algebras} are 
 $\mathbb{DHMSH}$-algebras satisfying the identity:
	
		\vspace{0.12cm}
		(H):  \texttt{$(x \wedge y) \to x \approx 1$.}
		\vspace{0.12cm}

\item $\mathbb{OCKSH}:$ {\em Ockham semi-Heyting algebras} are 
$\mathbb{DHMSH}$-algebras satisfying the identity:

\ecuacionDef{$(x \lor y)' \approx x' \land y'$.} \label{071217_02}

\item $\mathbb{OCKH}:$ {\em Ockham Heyting algebras} are 
$\mathbb{OCKSH}$-algebras satisfying the identity (H).

\item $\mathbb{D}ms\mathbb{SH}:$ {\em Dually {\rm ms} semi-Heyting algebras} are 
$\mathbb{OCKSH}$-algebras satisfying the identity:

\ecuacionDef{$x'' \leq x.$} \label{111217_05}

\item $\mathbb{D}ms\mathbb{H}:$ {\em Dually {\rm ms} Heyting algebras} are $\mathbb{D}ms\mathbb{SH}$-algebras satisfying the identity (H).

\item $\mathbb{DMSH}:$ {\em De Morgan semi-Heyting  algebras} are 
$\mathbb{OCKSH}$-algebras (or $\mathbb{DHMSH}$-algebras) satisfying the identity 

 \ecuacionDef{$x'' \approx x.$} \label{111217_07}

\item $\mathbb{DMH}:$ {\em De Morgan Heyting algebras} are $\mathbb{DMSH}$-algebras satisfying the identity (H).

\item $\mathbb{DSDSH}:$ {\em Dually semi-De Morgan semi-Heyting algebras 
{\rm (\cite{sankappanavar1985semiDeMorgan})} are $\mathbb{DHMSH}$-algebras satisfying the identities:

\ecuacionDef{$(x \vee y)'' \approx x'' \vee y''$,} \label{131217_09}

\ecuacionDef{$x''' \approx x'$.} \label{131217_10}

\item $\mathbb{DSDH}:$ {\em Dually semi-De Morgan Heyting algebras} are $\mathbb{DSDSH}$-algebras satisfying the identity \rm(H).

\item $\mathbb{DQDSH}:$ {\em Dually quasi-De Morgan semi-Heyting algebras} {\rm (\cite{sankappanavar1985semiDeMorgan})} are 
$\mathbb{DSDSH}$-algebras satisfying the identity {\rm (E\ref{111217_05})}.

\item $\mathbb{DQDH}:$ {\em Dually quasi-De Morgan Heyting algebras} are $\mathbb{DQDSH}$-algebras satisfying the identity (H).




\item $\mathbb{DPCSH}:$ {\em Dually pseudocomplemented semi-Heyting algebras} are 
$\mathbb{DQDSH}$-algebras satisfying the identity: 

\ecuacionDef{$x \vee x' \approx 1$.} \label{131217_14}

\item $\mathbb{DPCH}:$ {\em Dually pseudocomplemented Heyting algebras} are 
$\mathbb{DPCSH}$-algebras satisfying the identity \rm(H\rm).

\item $\mathbb{BDQDSH}:$ {\em Blended dually quasi-De Morgan semi-Heyting algebras } are 
$\mathbb{DQDSH}$-algebras satisfying the identity: 

\ecuacionDef{$(x \vee x^*)' \approx x' \wedge x{^*}'$ \rm(Blended $\vee$-De Morgan law).} \label{181217_01} 

\item $\mathbb{BDQDH}:$ {\em Blended dually quasi-De Morgan Heyting algebras } are 
$\mathbb{BDQDSH}$-algebras satisfying the identity (H).

\item $\mathbb{SBDQDSH}:$ {\em Strongly blended dually quasi-De Morgan semi-Heyting algebras} are $\mathbb{DQDSH}$-algebras satisfying the identity: 

\ecuacionDef{$(x \vee y^*)' \approx x' \wedge y{^*}'$} } \label{181217_03}  \rm(Strongly blended $\vee$-De Morgan law). 

\item $\mathbb{SBDQDH}:$ {\em Strongly blended dually quasi-De Morgan Heyting algebras } are 
$\mathbb{SBDQDSH}$-algebras satisfying the identity (H).

 \item $\mathbb{DQDBSH}:$ {\em Dually quasi-De Morgan Boolean semi-Heyting algebras} are 
$\mathbb{DQDSH}$-algebras satisfying the identity: 

\ecuacionDef{$x \vee x^* \approx 1$.}  \label{181217_05}

\item $\mathbb{DQDBH}:$ {\em dually quasi-De Morgan Boolean Heyting  algebras } are 
$\mathbb{DQDBSH}$-algebras satisfying the identity (H).

\item $\mathbb{DQS}\rm t\mathbb{SH}:$ {\em Dually quasi-Stone semi-Heyting algebras} are 
$\mathbb{DHMSH}$-algebras satisfying the identities: (E\ref{111217_05}),

\ecuacionDef{$(x \vee y')' \approx x' \wedge y''$ (weak $\vee$-De Morgan law),} \label{181217_07}

\ecuacionDef{$x' \wedge x'' \approx 0$ (Dual Stone identity).} \label{181217_09}

\item $\mathbb{DQS}\rm t\mathbb{H:} $ {\em Dually quasi-Stone Heyting algebras } are 
$\mathbb{DQS}\rm t\mathbb{SH}$-algebras satisfying the identity (H).

\item $\mathbb{BDQS}\rm t\mathbb{SH}:$ {\em Blended dually quasi-Stone semi-Heyting algebras} are 
$\mathbb{DQS}\rm t\mathbb{SH}$-algebras satisfying  the identity 
(E\ref{181217_01}).  

\item $\mathbb{BDQS}\rm t\mathbb{H:} $ {\em Blended dually quasi-Stone Heyting algebras } are 
$\mathbb{BDQS}\rm t\mathbb{SH}$-algebras satisfying the identity (H).

\item $\mathbb{SBDQS}\rm t\mathbb{SH}:$ {\em Strongly blended dually quasi-Stone semi-Heyting algebras} are $\mathbb{DQS}\rm t\mathbb{SH}$-algebras satisfying the identity (E12). 

\item $\mathbb{SBDQS}\rm t\mathbb{H:} $ {\em 
Strongly blended dually quasi-Stone Heyting algebras } are 
$\mathbb{SBDQS}\rm t\mathbb{SH}$-algebras satisfying the identity (H).

\item $\mathbb{DS}\rm t\mathbb{SH}:${\em Dually Stone semi-Heyting algebras} are $\mathbb{DPCSH}$-algebras satisfying the identity (E15). 

\item $\mathbb{DS}\rm t\mathbb{H}:${\em Dually Stone Heyting algebras} are $\mathbb{DS}\rm t\mathbb{SH}$-algebras satisfying the identity (H). 

\item $\mathbb{DSCSH}:$ {\em Dually semi-complemented semi-Heyting algebras } are 
$\mathbb{DHMSH}$-algebras satisfying the identity (E\ref{131217_14}). 

 \item $\mathbb{DSCH:} $ {\em Dually semi-complemented Heyting algebras } are 
$\mathbb{DSCSH}$-algebras satisfying the identity (H).

\item $\mathbb{DDPCSH}:$ {\em Dually demi-pseudocomplemented semi-Heyting algebras} are $\mathbb{DSDSH}$-algebras satisfying the identity: 

\ecuacionDef{$x' \vee x'' \approx 1.$} \label{211118_01} 

 \item $\mathbb{DDPCH:} $ {\em Dually demi-pseudocomplemented Heyting algebras } are 
$\mathbb{DDPCSH}$-algebras satisfying the identity (H).

\item $\mathbb{DAPCSH}:$ {\em Dually almost pseudocomplemented semi-Heyting algebras} are \\
$\mathbb{DDPCSH}$-algebras in which $'$ satisfies the identity dual to (E\ref{181217_05}). 
\label{ultimoNumeroIdentidad}

 \item $\mathbb{DAPCH:} $ {\em Dually almost pseudocomplemented Heyting algebras } are 
$\mathbb{DAPCSH}$-algebras in which $'$ satisfies the identity (H).


\end{enumerate}
	



 Next, we present Hilbert-type axiomatizations for the (new) logics that are extensions of $\mathcal{DHMSH}$ and that correspond to the subvarieties of $\mathcal{DHMSH}$ mentioned in LIST 1.  For the relationships of these logics to the varieties in LIST 1, the reader is referred to Theorem \ref{completenessAlgrLogics} below.  \\

{\bf LIST 2: SOME IMPORTANT EXTENSIONS OF $\mathcal{DHMSH}$} 

\begin{enumerate}
	\item $\mathcal{DHMH}:$ The dually hemimorphic Heyting logic is the extension of $\mathcal{DHMSH}$ given by

\numeroAxioma{$(\alpha \wedge \beta) \to \alpha$.} \label{axiom_Heyting}

	\item $\mathcal{OCKSH}:$ The Ockham semi-Heyting logic is the extension of $\mathcal{DHMSH}$ given by

	\numeroAxioma{$\sim(\alpha \vee \beta) \to_H (\sim\alpha \wedge \sim\beta)$,} \label{axiom_OCKSH_1}	
	
	\numeroAxioma{$(\sim\alpha \wedge \sim\beta) \to_H \ \sim(\alpha \vee \beta)$.} \label{axiom_OCKSH_2}

\item $\mathcal{OCKH}:$ The Ockham Heyting logic is the extension of $\mathcal{OCKSH}$ given by
(A\ref{axiom_Heyting}).

	\item $\mathcal{D}ms\mathcal{SH}:$ The dually {\rm ms} semi-Heyting logic is the extension of $\mathcal{OCKSH}$ given by

	\numeroAxioma{$\sim\sim\alpha \to_H \alpha$.} \label{axioma_involucion_neg2} \label{axiom_DmsSH}

\item $\mathcal{D}ms\mathcal{H}:$ The dually {\rm ms} Heyting logic is the extension of $\mathcal{D}ms\mathcal{SH}$ given by (A\ref{axiom_Heyting}).

	\item $\mathcal{DMSH}:$ The De Morgan semi-Heyting logic is the extension of $\mathcal{OCKSH}$ given by (A\ref{axioma_involucion_neg2}) and

	\numeroAxioma{$\alpha \to_H \ \sim\sim\alpha$.} \label{axioma_involucion_neg}

	\item $\mathcal{DMH}:$ The De Morgan Heyting logic is the extension of $\mathcal{DMSH}$ given by (A\ref{axiom_Heyting}).

	\item $\mathcal{DSDSH}:$ The dually semi-De Morgan semi-Heyting logic is the extension of $\mathcal{DHMSH}$ given by

	\numeroAxioma{$\sim\sim(\alpha \vee \beta) \to_H (\sim\sim\alpha \ \vee \ \sim\sim\beta)$,} \label{axiom_DSDSH_1}

	\numeroAxioma{$(\sim\sim\alpha \ \vee \ \sim\sim \beta) \to_H \ \sim\sim(\alpha \vee \beta)$,}  \label{axiom_DSDSH_2}

	\numeroAxioma{$\sim\sim\sim\alpha \to_H \ \sim\alpha$,} \label{axiom_DSDSH_3}

	\numeroAxioma{$\sim\alpha \to_H \  \sim\sim\sim\alpha$.} \label{axiom_DSDSH_4}
	
\item $\mathcal{DSDH}:$ The dually semi-De Morgan Heyting logic is the extension of $\mathcal{DSDSH}$ given by (A\ref{axiom_Heyting}).

	\item $\mathcal{DQDSH}:$ The dually quasi-De Morgan semi-Heyting logic is the extension of $\mathcal{DSDSH}$ given by (A\ref{axiom_DmsSH}). 
	
\item $\mathcal{DQDH}:$ The dually quasi-De Morgan Heyting logic is the extension of $\mathcal{DQDSH}$ given by (A\ref{axiom_Heyting}).

	\item $\mathcal{DPCSH}:$ The dually pseudocomplemented semi-Heyting logic is the extension of $\mathcal{DQDSH}$ given by

	\numeroAxioma{$\alpha \ \vee \ \sim\alpha$.} \label{axiom_DPCSH}

	\item $\mathcal{DPCH}:$ The dually pseudocomplemented Heyting logic is the extension of $\mathcal{DPCSH}$ given by (A\ref{axiom_Heyting}).

	\item $\mathcal{BDQDSH}:$ The blended dually quasi-De Morgan semi-Heyting logic is the extension of $\mathcal{DQDSH}$ given by

	\numeroAxioma{$\sim(\alpha \vee (\alpha \to \bot)) \to_H (\sim\alpha \ \wedge \ \sim(\alpha \to \bot))$,} \label{axiom_BDQSSH_1}

	\numeroAxioma{$(\sim\alpha \ \wedge \ \sim(\alpha \to \bot)) \to_ H \ \sim (\alpha \ \vee \ (\alpha \to \bot))$.} \label{axiom_BDQSSH_2}

\item $\mathcal{BDQDH}:$  The blended dually quasi-De Morgan Heyting logic is the extension of $\mathcal{BDQDSH}$ given by
(A\ref{axiom_Heyting}).

	\item $\mathcal{SBDQDSH}:$ The strongly blended quasi-De Morgan semi-Heyting logic is the extension of $\mathcal{DQDSH}$ given by

	\numeroAxioma{$\sim(\alpha \vee (\beta \to \bot)) \to_H (\sim \alpha \ \wedge \ \sim(\beta \to \bot))$,} \label{axiom_SBDQDSH_1}

	\numeroAxioma{$(\sim \alpha \ \wedge \ \sim(\beta \to \bot)) \to_ H \ \sim(\alpha \vee (\beta \to \bot))$.} \label{axiom_SBDQDSH_2}

 \item $\mathcal{SBDQDH}:$  The strongly blended dually quasi-De Morgan Heyting logic is the extension of $\mathcal{SBDQDSH}$ given by
(A\ref{axiom_Heyting}).
	
 \item $\mathcal{DQDBSH}:$ The dually quasi-De Morgan Boolean semi-Heyting logic is the extension of $\mathcal{DQDSH}$ given by (A15).

	\numeroAxioma{$\alpha \vee (\alpha \to \bot)$.}  \label{axiom_DQDBSH}
	
\item $\mathcal{DQDBH}:$ The dually quasi-De Morgan Boolean Heyting logic is the extension of $\mathcal{DQDBSH}$ given by (A\ref{axiom_Heyting}).

	\item $\mathcal{DQS}\rm t\mathcal{SH}:$ The dually quasi-Stone semi-Heyting logic is the extension of $\mathcal{DHMSH}$ given by 
	(A\ref{axiom_DmsSH}) and the following axioms:

	\numeroAxioma{$\sim(\alpha \ \vee \ \sim\beta) \to_H (\sim\alpha \ \wedge \ \sim\sim\beta)$,} \label{axiom_DQSSH_1}	
	
	\numeroAxioma{$(\sim\alpha \ \wedge \ \sim\sim\beta) \to_H \ \sim(\alpha \ \vee \sim\beta)$,} \label{axiom_DQSSH_2}

	\numeroAxioma{$(\sim\alpha \ \wedge \sim\sim\alpha) \to_H \bot$.} \label{axiom_DQSSH_3}
	
\item $\mathcal{DQS}\rm t\mathcal{H}:$  The dually quasi-Stone Heyting logic is the extension of $\mathcal{DQS}\rm t\mathcal{SH}$ given by
(A\ref{axiom_Heyting}).

	\item $\mathcal{BDQS}\rm t\mathcal{SH}:$ The blended dually quasi-Stone semi-Heyting logic is the extension of $\mathcal{DHMSH}$ given by (A\ref{axiom_DmsSH}), (A\ref{axiom_BDQSSH_1}), (A\ref{axiom_BDQSSH_2}), (A\ref{axiom_DQSSH_1}), (A\ref{axiom_DQSSH_2}) and (A\ref{axiom_DQSSH_3}).
	
\item $\mathcal{BDQS}\rm t\mathcal{H}:$  The blended dually quasi-Stone Heyting logic is the extension of $\mathcal{BDQS}\rm t\mathcal{SH}$ given by
(A\ref{axiom_Heyting}).

	\item $\mathcal{SBDQS}\rm t\mathcal{SH}:$ The strongly blended dually quasi-Stone semi-Heyting logic is the extension of $\mathcal{DHMSH}$ given by (A\ref{axiom_DmsSH}),   (A\ref{axiom_SBDQDSH_1}), (A\ref{axiom_SBDQDSH_2}),  (A\ref{axiom_DQSSH_1}), (A\ref{axiom_DQSSH_2}) and (A\ref{axiom_DQSSH_3}).
	
 \item $\mathcal{SBDQS}\rm t\mathcal{H}:$  The strongly blended dually quasi-Stone Heyting logic is the extension of $\mathcal{SBDQS}\rm t\mathcal{SH}$ given by
(A\ref{axiom_Heyting}).

\item $\mathcal{DS}\rm t \mathcal{SH}:$ The dually Stone semi-Heyting logic is the extension of $\mathcal{DPCSH}$ given by (A32).          

\item $\mathcal{DS}\rm t \mathcal{H}:$ The dually Stone Heyting logic is the extension of $\mathcal{DS}\rm t \mathcal{SH}$ given by (A15).                             
	
	\item $\mathcal{DSCSH}:$ The dually semi-complemented semi-Heyting logic is the extension of $\mathcal{DHMSH}$ given by (A\ref{axiom_DPCSH}). 
	
\item $\mathcal{DSCH}:$ The dually semi-complemented Heyting logic is the extension of $\mathcal{DSCSH}$ given by (A\ref{axiom_Heyting}).	
	
	\item $\mathcal{DDPCSH}:$ The dually demi-pseudocomplemented semi-Heyting logic is the extension of 
	$\mathcal{DSDSH}$  given by 
	
		\numeroAxioma{$\sim\alpha \ \vee \ \sim\sim\alpha$.} \label{axiom_DDPCSH}

	\label{B_ultimoNumeroIdentidad}

\item $\mathcal{DDPCH}:$ The dually demi-pseudocomplemented Heyting logic is the extension of $\mathcal{DDPCSH}$ given by (A\ref{axiom_Heyting}).	

\item $\mathcal{DAPCSH}:$ {\em The dually almost pseudocomplemented semi-Heyting logic} is the extension of $\mathcal{DDPCSH}$ given by
  
  $\sim\sim\alpha \to_H \alpha$.
  
\item $\mathcal{DAPCH}:$ The dually almost pseudocomplemented Heyting logic is the extension of $\mathcal{DAPCSH}$ given by (A\ref{axiom_Heyting}). 
	
\end{enumerate}	

\medskip	

The following theorem which is immediate from Theorem \ref{completeness_DHMSH_extension} describes the correspondence between the logics in LIST 2 and the varieties in LIST 1.

 \begin{Theorem} \label{completenessAlgrLogics} \label{Extensions}
Let $\mathbb V_i$ be the variety of algebras mentioned in the $i$-th item of {\bf LIST 1} and $\mathcal V_i$ be the logic  
appearing in the $i$-th item of  {\bf LIST 2}. 	Then, the logic $\mathcal V_i$ corresponds to the variety $\mathbb V_i$  in the sense that $\mathbb V_i$ is its equivalent algebraic semantics for $\mathcal V_i$.
 \end{Theorem}

In the figure below, we present a (partial) poset describing the mutual relationships among the varieties, whose names end with ``SH'', mentioned in PART 1 above.  The dual of this poset will show the relations among the logics, whose names end with ``SH'', presented in PART 2, $\mathbb{T}$ being the trivial variety. Note that the links do not necessarily represent covers.  

\newpage

{\small
\setlength{\unitlength}{1cm}
\begin{picture}(16,22)(0,0)
\put(8,1){\circle*{0.15}}
\put(7.3,0.9){$\mathbb T$}

\put(6,3){\circle*{0.15}}
\put(4.8,3){$\mathbb V(\mathbf{2}^e)$}
\put(10,3){\circle*{0.15}}
\put(8.6,3){$\mathbb V(\mathbf{\bar{2}^e})$}

\put(8,5){\circle*{0.15}}
\put(6,5){$\mathbb V(\mathbf{2^e}, \mathbf{\bar{2}^e})$}

\put(4,7){\circle*{0.15}}
\put(2.6,7){$\mathbb{DS}\rm t \mathbb{SH}$}
\put(12,7){\circle*{0.15}}
\put(12.2,7){$\mathbb{DQDBSH}$}

\put(6,9){\circle*{0.15}}
\put(4.6,9){$\mathbb{DPCSH}$}
\put(10,9){\circle*{0.15}}
\put(8.6,9){$\mathbb{DMSH}$}

\put(8,11){\circle*{0.15}}
                                         \put(8.1,11){$\mathbb{SBDQDSH}$}

\put(4,13){\circle*{0.15}}
\put(2,13){$\mathbb{SBDQS}\rm t \mathbb{SH}$}
\put(6,13){\circle*{0.15}}
\put(4.4,13){$\mathbb{DAPCSH}$}
\put(8,13){\circle*{0.15}}
\put(6.4,13){$\mathbb{BDQDSH}$}
\put(10,13){\circle*{0.15}}
\put(8.6,13){$\mathbb{D}{\rm ms}\mathbb{SH}$}

\put(4,15){\circle*{0.15}}
\put(2.2,15){$\mathbb{BDQS}\rm t \mathbb{SH}$}
\put(6,15){\circle*{0.15}}
\put(4.4,15){$\mathbb{DDPCSH}$}
\put(8,15){\circle*{0.15}}
\put(6.6,15){$\mathbb{DQDSH}$}

\put(4,17){\circle*{0.15}}
\put(2.4,17){$\mathbb{DQS}\rm t \mathbb{SH}$}
                     \put(8,17){\circle*{0.15}}
                     \put(6.6,17){$\mathbb{DSDSH}$}
\put(10,17){\circle*{0.15}}
\put(8.4,17){$\mathbb{OCKSH}$}
                                                             \put(6,18){\circle*{0.15}}

\put(8,19){\circle*{0.15}}
\put(6.4,19){$\mathbb{DHMSH}$}

\put(4,17){\line(0,-1){10}}
\put(6,15){\line(0,-1){6}}
\put(8,19){\line(0,-1){8}}
        \put(10,17){\line(0,-1){8}}
                          \put(6,15){\line(0,1){3}}
                          \put(4.5,18){$\mathbb{DSCSH}$}

\put(8,1){\line(1,1){2}}
\put(6,3){\line(1,1){2}}
\put(8,5){\line(2,1){4}}
\put(4,7){\line(1,1){2}}
\put(6,9){\line(1,2){2}}
\put(6,15){\line(1,1){2}}
           \put(4,17){\line(2,1){4}}

\put(6,3){\line(1,-1){2}}
\put(8,5){\line(1,-1){2}}
\put(4,7){\line(2,-1){4}}
\put(8,11){\line(1,-1){4}}
        \put(8,19){\line(1,-1){2}}
        \put(7.3,0.2){{\large{Figure 1}}}
\end{picture}
}
\begin{center}
{\bf PARTIAL POSET OF SUBVARIETIES OF $\mathbb{DHMSH}$}
\end{center}

Note that it is not yet known if $\mathbb{DPCSH} \subset \mathbb{SBDQDSH}$, although we conjecture it to be true.  So we mention the following open problem.

\noindent PROBLEM: Is $\mathbb{DPCSH} \subset \mathbb{SBDQDSH}$?  In particular, is $\mathbb{DPCH} \subset \mathbb{SBDQDH}$?

\begin{Corollary}  \label{Moisil_Logic_Theorem}

	The {\rm (Moisil's)} logic $\mathcal{LM}$ is equivalent to the logic $\mathcal{DMH}$.
\end{Corollary}
  
\begin{Proof}
We know from Moisil's result (see Monteiro \cite{monteiro1980algebres}) that $\mathcal{LM}$ corresponds to $\mathbb{DMH}$.  Also, observe from Theorem  \ref{Extensions}  
that the logic $\mathcal{DMH}$ correspond to $\mathbb{DMH}$ as well.
\end{Proof}

\medskip

\begin{definition}
Let $\mathcal{L}$ be an algebraizable logic. 
We say that $\mathcal{L}$ is a discriminator logic if its corresponding (quasi-) variety is a discriminator variety.  Furthermore,  $\mathcal{L}$ is a (quasi)primal logic if its corresponding (quasi-) variety is a variety generated by a (quasi)primal algebra.   
\end{definition}

The classical logic is the first well-known example of a primal logic (as the Boolean algebra $\mathbf{2}$ is a primal algebra).

\begin{Remark}
It was shown in \cite{sankappanavar2011expansions} that the variety $\mathbb{RDQDS}\rm t\mathbb{SH}_1$ is a discriminator variety.  Thus  $\mathcal{RDQDS}\rm t\mathcal{SH}_1$ is a discriminator logic.  Most of the logics considered in the rest of this paper are discriminator logics.  We will point them out as they appear. 
\end{Remark}
\medskip

We conclude this section by noting that 
the lattice of extensions of the logic $\mathcal{DMH}$ is an interval of the lattice of extensions of $\mathcal{DMSH}$, which, in turn, is an interval in the lattice of extensions of  $\mathcal{DHMSH}$.\\

\section{ The Deduction Theorem in extensions of $\mathcal{DHMSH}$}    
\label{deducTheoExtensions}
\
\indent In this section we first show that the ``usual'' form of the Deduction Theorem'' fails in the logic $\mathcal{DHMSH}$, and then characterize those extensions of $\mathcal{DHMSH}$ where it does hold.

	A logic $\mathcal L$ is said to have the {\it Deduction Property} for the connective $\to$ if the following statement holds:\\
	$$\Gamma, \alpha \vdash_{\mathcal{L}} \beta \text{ if and only if } \Gamma \vdash_{\mathcal{L}} \alpha \to \beta,$$ for all $\Gamma\cup\{\alpha, \beta\}\subseteq \mbox{\it Fm}$.

	In the logic $\mathcal{SI}$ the Deduction Property for the conective $\to_H$ is known to hold  \cite[Theorem 3.18]{cornejo2011semi}.   But, this property fails in the logic $\mathcal{DHMSH}$, as shown in the following remark.

\begin{Remark}
	First, we note, by Lemma {\rm\ref{lemma_properties_DMSH_logic}} {\rm (\ref{IL3_partial})},  that 
\begin{equation} \label{eqn_A}
	\phi \to_H \psi \ \vdash_\mathcal{DHMSH} \  \sim\psi \to_H \ \sim \phi.  
\end{equation}

\noindent Consider the algebra $\mathbf{L}_1^{dm} \in \mathbb{DHMSH}$ defined in Section \rm\ref{S8}.              
	Observe that $\mathbf{L}_1^{dm} \  \not\models_{\mathbb{DHMSH}} \ (x \to_H y) \to_H (y' \to_H x') \approx 1$ (by taking $x = 1$ and $y = a$).  Hence, $\mathbb{DHMSH} \not\models (x \to_H y) \to_H (y' \to_H x') \approx 1$ and therefore, by {\rm Theorem} {\rm \ref{completeness_DMSH}}, 
	$$\not\vdash_{\mathcal{DHMSH}} (\phi \to_H \psi) \to_H \  (\sim\psi \to_H \ \sim \phi).$$  
Thus, the Deduction Property fails in $\mathcal{DHMSH}$, in view of \eqref{eqn_A}.
\end{Remark}

We now wish to characterize the extensions of 
$\mathcal{DHMSH}$ in which the Deduction Property holds.  For this, we need a preliminary lemma.

\begin{Lemma} \label{lem_040417_04}
Let  $\mathcal E$ be an extension of the logic $\mathcal{DHMSH}$ such that 
$$\vdash_{\mathcal E} (\alpha \to_H \ \beta) \to_H (\sim \beta \to_H \ \sim\alpha).$$  
Then  $\mathcal E$ satisfies the Deduction Property for the connective $\to_H$.
\end{Lemma}

\begin{Proof}
Assume that  $\Gamma \cup \{\phi\} \vdash_\mathcal{E} \psi$. We shall prove $\Gamma \vdash_\mathcal{E} \phi \to_H \psi$ by induction on the proof for  $\psi$. By hypothesis,
\begin{equation} \label{050417_01}
 \vdash_{\mathcal E} (\alpha \to_H \ \beta) \to_H (\sim \beta \to_H \ \sim\alpha).
\end{equation}
If $\psi$ is an axiom of $\mathcal{E}$ or a formula in $\Gamma$, then $\Gamma \vdash_\mathcal{E} \psi$. By  Lemma \ref{lema_propiedades_SI}, part (\ref{300513_15}) we have  $\Gamma \vdash_\mathcal{E} \phi \to_H \psi.$

Let us assume that $\Gamma \cup \{\phi\} \vdash_\mathcal{E} \psi$ is the result of applying the rule (SMP). Then  we may assume that there is some formula $\alpha$ such that  $\Gamma \cup \{\phi\} \vdash_\mathcal{E} \alpha$ and $\Gamma \cup \{\phi\} \vdash_\mathcal{E} \alpha \to_H \psi$. So, by inductive hypothesis, we have,
\begin{enumerate}[1.]
	\item $\Gamma \vdash_\mathcal{E} \phi \to_H \alpha$,  \label{230513_02}
	\item $\Gamma \vdash_\mathcal{E} \phi \to_H (\alpha \to_H \psi),$ \label{230513_04}
	\item $\Gamma \vdash_\mathcal{E} \phi \to_H \phi$ by Lemma \ref{lema_propiedades_SI}, part (\ref{110613_12}), \label{230513_03}
	\item $\Gamma \vdash_\mathcal{E} \phi \to_H (\phi \wedge \alpha)$ by (A\ref{axioma_infimo_cota_superior}) and SMP applied to \ref{230513_02} and \ref{230513_03}, \label{230513_05}
	\item $\Gamma \vdash_\mathcal{E} (\phi \wedge \alpha) \to_H \psi$ by (A\ref{axioma_condicRes_ImplicAInf}) and  SMP applied to \ref{230513_04}, \label{230513_06}
	\item $\Gamma \vdash_\mathcal{E} \phi \to_H \psi$ by Lemma \ref{lema_propiedades_SI} (\ref{020713_42}) and SMP applied to \ref{230513_05} and \ref{230513_06}.
\end{enumerate}
Assume that 
$\Gamma \cup \{\phi\} \vdash_\mathcal{E} \psi$ is the result of applying the rule (SCP).  Hence  $\psi = \ \sim\beta \to_H \sim\alpha$ and $\Gamma \cup \{\phi\} \vdash_\mathcal{E} \alpha \to_H \beta$. By induction we have that 
\begin{enumerate}[1.]
	\item $\Gamma \vdash_\mathcal{E} \phi \to_H (\alpha \to_H \beta)$, \label{050417_02}
	\item $\Gamma \vdash_{\mathcal E} (\alpha \to_H \beta) \to_H (\sim\beta \to_H \ \sim\alpha)$ by (\ref{050417_01}), \label{050417_03}
	\item $\Gamma \vdash_{\mathcal E} \phi \to_H (\sim\beta \to_H \ \sim\alpha)$ by \ref{lema_propiedades_SI} (\ref{020713_42}) and SMP applied to \ref{050417_02} and \ref{050417_03}.
\end{enumerate}

For the other implication, we assume that  $\Gamma \vdash_\mathcal{E} \phi \to_H \psi$. Then $\Gamma  \cup \{\phi\} \vdash_\mathcal{E} \phi \to_H \psi$. Since $\Gamma \cup \{\phi\} \vdash_\mathcal{E} \phi$, we have $\Gamma  \cup \{\phi\} \vdash_\mathcal{E} \psi$ by (SMP).
\end{Proof}

\begin{Theorem}\label{Deduction}
	The Deduction Property holds in an extension $\mathcal E$ of the logic $\mathcal{DHMSH}$ for the connective $\to_H$ if and only if $\mathcal E \vdash (\alpha \to_H \ \beta) \to_H (\ \sim\beta \to_H \ \sim\alpha)$.
\end{Theorem}

\begin{Proof}
	Let us assume that the Deduction Property holds in $\mathcal E$ for the conective $\to_H$. Note that  $\alpha \to_H  \beta \vdash_{\mathcal E} \  \alpha \to_H   \beta $ and 
	$\alpha \to_H  \beta \vdash_{\mathcal E} \  \sim\beta  \to_H \ \sim\alpha$ by (SCP).  Hence
	$\vdash_{\mathcal E} (\alpha \to_H \ \beta) \to_H (\sim\beta \to_H \ \sim\alpha)$ by Deduction Property, or equivalently, $\mathcal E \vdash_{\mathcal{DHMSH}} (\alpha \to_H \ \beta) \to_H (\sim \beta \to_H \sim\alpha)$. 
		
	For the converse, let us assume that  $\mathcal E \vdash (\alpha \to_H \ \beta) \to_H (\sim\beta \to_H \sim\alpha)$. 	By Lemma \ref{lem_040417_04}, the Deduction Property holds in $\mathcal E$ for the conective $\to_H$.
\end{Proof}

Recall that semi-Heyting algebras are pseudocomplemented with $x^*:= x \to 0$ as the pseudocomplement of $x$.  A semi-Heyting algebra $\mathbf{L}$ is a \emph {Stone semi-Heyting algebra}  
 if $\mathbf{L}$ satisfies the Stone identity: $ x^* \lor x^{**} \approx 1$.  Let $\mathbb{S}\rm t\mathbb{SH}$
 denote the variety of Stone semi-Heyting algebras. 
Recall also that if $\mathbf{A}$ is a semi-Heyting algebra, then $\langle A, \lor, \land, \to_H 0, 1 \rangle$ is a Heyting algebra.

	\begin{Lemma} \label{lemma_270219_02}
		Let $\mathbf A \in \mathbb{DHMSH}$ such that $\mathbf A \models  (x \to_H y) \to_H (y' \to_H x') \approx 1$.  Then 
\begin{enumerate}[{\rm(a)}]
	\item $\mathbf A  \models x \wedge x' \approx 0,$ \label{270219_01}
	\item $\mathbf A  \models x^* \approx x'$, \label{260219_03}
	\item $\mathbf A  \models x^\ast \lor x^{\ast\ast} \approx 1.$
\end{enumerate}
	\end{Lemma}

	\begin{Proof}
		
		Let $a \in A$. 
\begin{enumerate}[(a)]
	\item Since $a \to_H (a' \to_H 0) = (1 \to_H a) \to_H (a' \to_H 0) = (1 \to_H a) \to_H (a' \to_H 1') = 1$ in view of hypothesis, we have that $a \wedge (a' \to_H 0) = a \wedge (a \to_H (a' \to_H 0)) = a \wedge 1 = a$. Hence
	\begin{equation}\label{equation_260219_02}
	a \wedge (a' \to_H 0) = a.
	\end{equation}
	Then, 
\noindent $	a \wedge a'$
$\overset{  (\ref{equation_260219_02}) 
}{=}  a \wedge (a' \to_H 0) \wedge a' $
$\overset{  
}{=}  a \wedge (a' \to 0) \wedge a' $
$\overset{  
}{=}  a \wedge a' \wedge 0 $
$\overset{  
}{=}  0,$  proving (a).\\ 

\item Observe that $a^* \to_H a' = a^* \to_H (1 \to_H a') = (a \to_H 0) \to_H (1 \to_H a') = (a \to_H 0) \to_H (0' \to_H a') = 1$ by hypothesis. Hence 
\begin{equation}\label{equation_260219_01}
\mathbf A \models x^* \leq x'.
\end{equation}

\noindent Next,\\
 $a' \wedge a^* $
$\overset{  
}{=}  a' \wedge (a \to 0) $
$\overset{  
}{=}  a' \wedge ((a' \wedge a) \to (a' \wedge 0)) $
$\overset{  (\ref{270219_01}) 
}{=}  a' \wedge (0 \to 0) $
$\overset{  
}{=}  a'. $
Hence $\mathbf A \models x' \leq x^*.$
\noindent Now, using (\ref{equation_260219_01}) we conclude that $a' = a^*$. \\

\item 		
\noindent $a^\ast \lor a^{\ast\ast} $
$\overset{  (\ref{260219_03}) 
}{=}  a' \lor a'' $
$\overset{  (E\ref{infDMlaw}) 
}{=}  (a \wedge a')' $
$\overset{  (\ref{270219_01}) 
}{=}  0' $
$\overset{  
}{=}  1, $
\end{enumerate}
proving the lemma. 
\end{Proof}

\begin{Lemma} \label{lemma_270219_01}
Let $\mathbf A \in \mathbb{DHMSH}$. Then the following conditions are equivalent in the algebra $\mathbf A$.
\begin{enumerate}[{\rm(1)}]
	\item $(x \to_H y) \to_H (y' \to_H x') \approx 1$, \label{260219_01}
	\item $x^* \approx x'$. \label{260219_02}
\end{enumerate}
\end{Lemma}

\begin{Proof} Let $a,b \in A$. 
Observe that (\ref{260219_01}) implies (\ref{260219_02}) from Lemma  \ref{lemma_270219_02}.   
For the converse, suppose $\mathbf A$ satisfies the identity (\ref{260219_02}).  Then, using (SH3) and the fact that $\to_H$ is a Heyting implication,  we have that $(a \to_H b) \to_H (b' \to_H a') = (a \to_H b) \to_H (b^\ast \to_H a^\ast) = (b^\ast \wedge (a \to_H b)) \to_H a^\ast = (b^\ast \wedge (((b^\ast \wedge a) \to_H (b^\ast \wedge b)) \to_H a^\ast = (b^\ast \wedge ((b^\ast \wedge a) \to_H 0)) \to_H a^\ast = (b^\ast \wedge (a \to_H 0)) \to_H a^\ast = (b^\ast \wedge  a^\ast) \to_H a^\ast = 1$, proving (1).
\end{Proof}

\begin{Lemma} \label{L6}
Let $\mathbf{L}$ be a Stone semi-Heyting algebra.  Let $\mathbf{L}^e$ be the expansion of $\mathbf{L}$ to the language $\langle \lor, \land,\to,  ', 0,1 \rangle$, where we define $'$ by: $x' := x^*$. 
Then 
\begin{enumerate}[{\rm(1)}]
	\item $\mathbf{L}^e \in \mathbb{DHMSH}$ and satisfies the identity: $ x' \approx x^*$,
	\item $\mathbf{L}^e \models (x \lor y)'' \approx x'' \lor y''$.
\end{enumerate}
\end{Lemma}

\begin{Proof}
The lemma clearly follows from the well-known facts that $\mathbf{L} \models (x \lor y)^* \approx x^* \land y^*$ and $\mathbf{L} \models (x \land y)^* \approx x^* \lor y^*$.
\end{Proof}\\

We will refer to the algebra $\mathbf{L}^e$ as an ``essentially a Stone semi-Heyting algebra''.\\

For $\mathbb{V}$ a subvariety of $\mathbb{S}t\mathbb{SH}$, we let 

$$\mathbb{V}^e:= \{\mathbf{L}^e: \ \mathbf{L} \in \mathbb V\}.$$

\medskip
It is clear that $\mathbb{V}^e$ is a subvariety of $\mathbb{DHMSH}$.\\

We are now ready to present our main result of this section that describes precisely those extensions of the logic $\mathcal{DHMSH}$ that have the Deduction Property.  The following theorem is immediate from Theorem \ref{Deduction}, Lemma \ref{lemma_270219_01} and Lemma \ref{L6}. 

\begin{Theorem}\label{Deduction5}
	The Deduction Property holds in an extension $\mathcal E$ of the logic $\mathcal{DHMSH}$ for the connective $\to_H$ if and only if the corresponding variety $\mathbb{E}$ is of the form 
$\mathbb{V}^e$, where $\mathbb{V} \subseteq \mathbb{S}\rm t\mathbb{SH}$. 
\end{Theorem}


\subsection{Deduction Theorem in the Extensions of the logic $\mathcal{DQDSH}$} \label{sec_six}

\  \\  

Recall 
that the variety $\mathbb{DQDSH}$ of dually quasi-De Morgan semi-Heyting algebras
and the corresponding extension $\mathcal{DQDSH}$ of $\mathcal{DHMSH}$ were defined in Section \ref{SecFive}.

\indent In this section we show that Theorem \ref{Deduction5} can be significantly improved for the extensions of the logic $\mathcal{DQDSH}$.   
In fact, we shall give an explicit description of the extensions of the logic $\mathcal{DQDSH}$ in which the Deduction Property holds.

For this purpose we need the following 2-element semi-Heyting algebras, $\mathbf{2}$ and  $\mathbf{\bar{2}}$
which are, up to isomorphism, the only two 2-element algebras in $\mathbb{SH}$. \\ 

$\mathbf{2}$:
\begin{tabular}{r|rr}
	$\to$: & 0 & 1\\
	\hline
	0 & 1 & 1 \\
	1 & 0 & 1
\end{tabular} \hspace{7cm}
$\mathbf{\bar{2}}$: 
\begin{tabular}{r|rr}
	$\to$: & 0 & 1\\
	\hline
	0 & 1 & 0 \\
	1 & 0 & 1
\end{tabular} \hspace{.5cm}
\smallskip
\begin{center} Figure 2
\end{center}

\medskip
The algebras $\mathbf{2^e}$ and  $\mathbf{\bar{2}^e}$ denote the expansions of the semi-Heyting algebras 
$\mathbf{2}$ and  $\mathbf{\bar{2}}$ by the unary operation $'$ defined as follows: $0'=1$ and $1'=0$.  It is clear that $\mathbf{2^e}$ and  $\mathbf{\bar{2}^e}$ are, up to isomorphism, the only two 2-element algebras in $\mathbb{DQDSH}$ (in fact, in $\mathbb{DMSH}).$  



\smallskip

\begin{Lemma} \label{teo_050417_05} 
	Let $\mathbb V$ be a subvariety of $\mathbb {DQDSH}$ such that $\mathbb V \models  (x \to_H y) \to (y' \to_H x') \approx 1$.  Then,
	$\mathbb V \subseteq \mathbb V(\mathbf{2^e}, \mathbf{\bar{2}^e})$, where $\mathbb V(\mathbf{2^e}, \mathbf{\bar{2}^e})$ denotes the variety generated by 
	$\{\mathbf{2^e}, \mathbf{\bar{2}^e}\}$.
\end{Lemma}

\begin{Proof}
	The hypothesis and Lemma \ref{lemma_270219_01} (\ref{260219_03}) imply that $\mathbb V \models x' \approx x^*$.  Hence
	$\mathbb V \  \subseteq \ \mathbb V(\mathbf{2^e}, \mathbf{\bar{2}^e})$ by \cite[Theorem 5.11]{sankappanavar2011expansions}.   
\end{Proof}

\ \\
The following theorem describes precisely those extensions of $\mathcal{DQDSH}$ in which the Deduction Property holds.  Let  $\mathbb{T}$ denote the trivial variety. 

\begin{Theorem} \label{Deduction4} 
	The Deduction Property holds in a logic $\mathcal E \in Ext(\mathcal{DQDSH})$ for $\to_H$ if and only if the corresponding variety is either either $\mathbb{T}$ or $\mathbb V(\mathbf{2^e})$ or $\mathbb V(\mathbf{\bar{2}^e})$ or $\mathbb V(\mathbf{2^e}, \mathbf{\bar{2}^e})$.   
\end{Theorem}

\begin{Proof}
	The theorem is immediate in view of Theorem \ref{completeness_DHMSH_extension},  Theorem \ref{Deduction},   
and Lemma \ref{teo_050417_05}.
\end{Proof}
\smallskip

\ 

Since $\mathbb{DMSH} \subseteq \mathbb{DQDSH}$ and $\mathbb{DPCSH} \subseteq \mathbb{DQDSH}$ (see \cite{sankappanavar2011expansions}), the following corollaries are immediate.

\begin{Corollary} \label{Deduction6} 
	The Deduction Property holds in a logic $\mathcal E \in Ext(\mathcal{DMSH})$ for $\to_H$ if and only if the corresponding variety is either            $\mathbb{T}$ or $\mathbb V(\mathbf{2^e})$ or $\mathbb V(\mathbf{\bar{2}^e})$ or $\mathbb V(\mathbf{2^e}, \mathbf{\bar{2}^e})$.  
\end{Corollary}

\begin{Corollary} \label{Deduction7} 
	The Deduction Property holds in a logic $\mathcal E \in Ext(\mathcal{DPCSH})$ for $\to_H$ if and only if the corresponding variety is either          $\mathbb{T}$ or $\mathbb V(\mathbf{2^e})$ or $\mathbb V(\mathbf{\bar{2}^e})$ or $\mathbb V(\mathbf{2^e}, \mathbf{\bar{2}^e})$.  
\end{Corollary}



\bigskip

\section{ Logics in $Ext(\mathcal{DQDSH})$ corresponding to subvarieties of $\mathbb{DQDSH}$ generated by finitely many finite algebras} \label{S7}  

\ 


 
 In this section, 
 as applications of Theorem \ref{teo_040417_01} and the algebraic results from \cite{sankappanavar2011expansions}, we will present several extensions of the logic 
$\mathcal{DQDSH}$ corresponding to subvarieties of $\mathbb{DQDSH}$ generated by finitely many finite algebras, 
thus providing a solution to PROBLEM 2.  

\subsection{2-valued Extensions of $\mathcal{DQDSH}$}  

\  \\

\indent It was shown in Theorem \ref{Deduction4} that the Deduction Property holds in an extension of the logic 
 $\mathcal{DQDSH}$ if and only if  the corresponding variety is a subvariety of $\mathbb V(\mathbf{2^e}, \mathbf{\bar{2}^e})$.  So, it is only natural to ask for the axiomatizations of the extensions of the logic $\mathcal{DQDSH}$ corresponding to the subvarieties of $\mathbb V(\mathbf{2^e}, \mathbf{\bar{2}^e})$.

The variety $\mathbb V(\mathbf{2^e}, \mathbf{\bar{2}^e})$ and its only non-trivial proper subvarieties $\mathbb V(\mathbf{2^e})$ and $\mathbb V(\mathbf{\bar{2}^e})$ were axiomatized in \cite[Theorem 5.11]{sankappanavar2011expansions}. 
$\mathbb V(\mathbf{2^e}, \mathbf{\bar{2}^e})$ is defined by the identity:  $x \leq x'^*$ (equivalently, $x \approx x'^*$), relative to the 
variety $\mathbb{DQDSH}$. 
 The varieties $\mathbb V(\mathbf{2^e})$ and $\mathbb V(\mathbf{\bar{2}^e})$ are defined, respectively, by the identities: $0 \to 1 \approx 1$ and $0 \to 1 \approx 0$, 
relative to $\mathbb V(\mathbf{2}^e, \mathbf{\bar{2}^e})$. 
\noindent In view of these observations,  
we obtain from Theorem \ref{teo_040417_01},  the following corollaries defining their corresponding logics.

 Let $\mathcal{L}(2^e,\bar{2}^e)$ (or $\mathcal{L}(\mathbb{V}(2^e,\bar{2}^e))$) be the extension of the logic $\mathcal{DQDSH}$ corresponding to the variety $\mathbb V(\mathbf{2}^e, \mathbf{\bar{2}^e})$.  Let $\alpha \Leftrightarrow_H \beta$ denote the formula: $(\alpha \rightarrow_H \beta) \land (\beta \rightarrow_H \alpha)$. 
 
\begin{Corollary}
The logic $\mathcal{L}(2^e,\bar{2}^e)$   
is defined, as an extension of the logic $\mathcal{DQDSH}$,
by the axiom:	
	\begin{enumerate}
	\item[] \qquad $(\sim\phi \to \bot) \Leftrightarrow_H \phi$.  
	\end{enumerate}
\end{Corollary}

 Let $\mathcal{L}({2^e})$ (or $\mathcal{L}(\mathbb{V}({2^e}$))) and $\mathcal{L}({\bar{2}^e})$ (or $\mathcal{L}(\mathbb{V}({\bar{2}^e})$)) denote, respectively, the extensions of the logic $\mathcal{L}({2^e,\bar{2}^e})$ corresponding to the varieties $\mathbb{V}(\mathbf{2}^e)$ and $\mathbb{V}(\mathbf{\bar{2}^e})$.
 
\begin{Corollary}
 The logic $\mathcal{L}({2^e})$ 
 is defined, as an extension of the logic $\mathcal{L}({2^e,\bar{2}^e})$, by the axiom:	
	\begin{enumerate}
	\item[]  \qquad $\bot \to \top.$  
          \end{enumerate}
\end{Corollary}

(We note that $\mathcal{L}({2^e})$ is yet another axiomatization of the classical logic.).
	
\begin{Corollary}
The logic $\mathcal{L}({\bar{2}^e})$ 
is defined, as an extension of the logic $\mathcal{L}({2^e,\bar{2}^e})$, by the axiom:	
	\begin{enumerate}
	\item[]  \qquad  $(\bot \to \top) \to_H \bot$.  
        \end{enumerate}
\end{Corollary}

\begin{Remark} Some features of the logics $\mathcal{L}(\mathbf{\bar{2}^e})$ and $\mathcal{L}(\mathbf{2^e})$:
\begin{itemize}
\item The logic $\mathcal{L}(\mathbf{\bar{2}^e})$ is ``anti-classical'' or ``contra-classical''
in the sense that the classically provable formula $\bot \to \top$ fails in it.  \\
\rm(It is somewhat perplexing to us that the intuitionists accept the principle that says, ``$False \to True =True''$.\rm )

\item The logics $\mathcal{L}(\mathbf{\bar{2}^e})$ and $\mathcal{L}(\mathbf{2^e})$ are coatoms in the lattice of extensions of the logic $\mathcal{DMSH}$ and hence that of $\mathcal{DQDSH}$ \rm(and of $\mathcal{DHMSH}$\rm).  

\item The implication $\to$ in $\mathcal{L}(\mathbf{\bar{2}^e})$ is commutative.  

\item The logics $\mathcal{L}(\mathbf{\bar{2}^e})$ and $\mathcal{L}(\mathbf{2^e})$ are not only disriminator logics, but, in fact, are primal logics, since $\mathbf{2^e}$ and $(\mathbf{\bar{2}^e})$ are  primal algebras. 
\item The logics $\mathcal{L}(\mathbf{\bar{2}^e})$ and $\mathcal{L}(\mathbf{2^e})$ do not have the Disjunction Property. (i.e., if $\alpha \lor \beta$ is provable, then $\alpha$ is provable or $\beta$ is provable.) 
\end{itemize}
\end{Remark}

More features of the logic $\mathcal{L}(\mathbf{\bar{2}^e})$ will be given in {\rm Remark} {\rm \ref{R_Connexive}}  of {\rm Section \ref{S8}}.

 \begin{Remark}
The Deduction Theorem holds only in the preceding three non-trivial logics in the lattice of extensions of the logic $\mathcal{DQDSH}$, in view of Theorem \ref{Deduction4}. 
\end{Remark}

\subsection{{\bf  3-valued Extensions of the logic $\mathcal{DQDSH}$}} \label{section_3valuedlogics}

\  \\

\indent It was shown in \cite{{sankappanavar2008semi}} that there are, up to isomorphism, ten $3$-element 
semi-Heyting algebras whose $\to$ operations are defined in figure $2$ below, where $0 < a < 1$.

\vspace{1cm}

\setlength{\unitlength}{1cm} \vspace*{.5cm}
\begin{picture}(15,2)
\put(1.5,-.4){\circle*{.15}}          

\put(1.5,1){\circle*{.15}}

\put(1.5,2.3){\circle*{.15}}

\put(1.1,-.4){$0$}                      

\put(1.1,.9){$a$}

\put(1.1,2.1){$1$}

\put(.1,.7){${\bf L}_1$ \ :}

\put(1.5,-.4){\line(0,1){2.7}}

\put(2.2,1){\begin{tabular}{c|ccc}
 $\to$ \ & \ 0 \ & \ $a$ \ & \ 1 \\ \hline
  & & & \\
  0 & 1 & 1 & 1 \\
  & & & \\
  $a$ & 0 & 1 & 1 \\
  & & & \\
  1 & 0 & $a$ & 1
\end{tabular}
}

\put(8.2,-.4){\circle*{.15}}                           

\put(8.2,1){\circle*{.15}}

\put(8.2,2.3){\circle*{.15}}

\put(7.8,-.4){$0$}                                     

\put(7.8,.9){$a$}

\put(7.8,2.1){$1$}

\put(6.7,.7){${\bf L}_2$ \ :}                        

\put(8.2,-.4){\line(0,1){2.7}}                        

\put(8.9,1){\begin{tabular}{c|ccc}                     
 $\to$ \ & \ 0 \ & \ $a$ \ & \ 1 \\ \hline
  & & & \\
  0 & 1 & $a$ & 1 \\
  & & & \\
  $a$ & 0 & 1 & 1 \\
  & & & \\
  1 & 0 & $a$ & 1
\end{tabular}
}
\end{picture}
\medskip

\setlength{\unitlength}{1cm} \vspace*{1.5cm}
\begin{picture}(15,2)
\put(1.5,-.4){\circle*{.15}}

\put(1.5,1){\circle*{.15}}

\put(1.5,2.3){\circle*{.15}}

\put(1.1,-.4){$0$}

\put(1.1,.9){$a$}

\put(1.1,2.1){$1$}

\put(.1,.7){${\bf L}_3$ \ :}

\put(1.5,-.4){\line(0,1){2.7}}

\put(2.2,1){\begin{tabular}{c|ccc}
 $\to$ \ & \ 0 \ & \ $a$ \ & \ 1 \\ \hline
  & & & \\
  0 & 1 & 1 & 1 \\
  & & & \\
  $a$ & 0 & 1 & $a$ \\
  & & & \\
  1 & 0 & $a$ & 1
\end{tabular}
}

\put(8.2,-.4){\circle*{.15}}

\put(8.2,1){\circle*{.15}}

\put(8.2,2.3){\circle*{.15}}

\put(7.8,-.4){$0$}

\put(7.8,.9){$a$}

\put(7.8,2.1){$1$}

\put(6.7,.7){${\bf L}_4$ \ :}

\put(8.2,-.4){\line(0,1){2.7}}

\put(8.9,1){\begin{tabular}{c|ccc}
 $\to$ \ & \ 0 \ & \ $a$ \ & \ 1 \\ \hline
  & & & \\
  0 & 1 & $a$ & 1 \\
  & & & \\
  $a$ & 0 & 1 & $a$ \\
  & & & \\
  1 & 0 & $a$ & 1
\end{tabular}
}
\end{picture}
\medskip

\setlength{\unitlength}{1cm} \vspace*{1.5cm}
\begin{picture}(15,2)
\put(1.5,-.4){\circle*{.15}}

\put(1.5,1){\circle*{.15}}

\put(1.5,2.3){\circle*{.15}}

\put(1.1,-.4){$0$}

\put(1.1,.9){$a$}

\put(1.1,2.1){$1$}

\put(.1,.7){${\bf L}_5$ \ :}

\put(1.5,-.4){\line(0,1){2.7}}

\put(2.2,1){\begin{tabular}{c|ccc}
 $\to$ \ & \ 0 \ & \ $a$ \ & \ 1 \\ \hline
  & & & \\
  0 & 1 & $a$ & $a$ \\
  & & & \\
  $a$ & 0 & 1 & 1 \\
  & & & \\
  1 & 0 & $a$ & 1
\end{tabular}
}

\put(8.2,-.4){\circle*{.15}}

\put(8.2,1){\circle*{.15}}

\put(8.2,2.3){\circle*{.15}}

\put(7.8,-.4){$0$}

\put(7.8,.9){$a$}

\put(7.8,2.1){$1$}

\put(6.7,.7){${\bf L}_6$ \ :}

\put(8.2,-.4){\line(0,1){2.7}}

\put(8.9,1){\begin{tabular}{c|ccc}
 $\to$ \ & \ 0 \ & \ $a$ \ & \ 1 \\ \hline
  & & & \\
  0 & 1 & 1 & $a$ \\
  & & & \\
  $a$ & 0 & 1 & 1 \\
  & & & \\
  1 & 0 & $a$ & 1
\end{tabular}
}
\end{picture}
\medskip

\setlength{\unitlength}{1cm} \vspace*{1.5cm}
\begin{picture}(15,2)
\put(1.5,-.4){\circle*{.15}}

\put(1.5,1){\circle*{.15}}

\put(1.5,2.3){\circle*{.15}}

\put(1.1,-.4){$0$}

\put(1.1,.9){$a$}

\put(1.1,2.1){$1$}

\put(.1,.7){${\bf L}_7$ \ :}

\put(1.5,-.4){\line(0,1){2.7}}

\put(2.2,1){\begin{tabular}{c|ccc}
 $\to$ \ & \ 0 \ & \ $a$ \ & \ 1 \\ \hline
  & & & \\
  0 & 1 & $a$ & $a$ \\
  & & & \\
  $a$ & 0 & 1 & $a$ \\
  & & & \\
  1 & 0 & $a$ & 1
\end{tabular}
}

\put(8.2,-.4){\circle*{.15}}

\put(8.2,1){\circle*{.15}}

\put(8.2,2.3){\circle*{.15}}

\put(7.8,-.4){$0$}

\put(7.8,.9){$a$}

\put(7.8,2.1){$1$}

\put(6.7,.7){${\bf L}_8$ \ :}

\put(8.2,-.4){\line(0,1){2.7}}

\put(8.9,1){\begin{tabular}{c|ccc}
 $\to$ \ & \ 0 \ & \ $a$ \ & \ 1 \\ \hline
  & & & \\
  0 & 1 & 1 & $a$ \\
  & & & \\
  $a$ & 0 & 1 & $a$ \\
  & & & \\
  1 & 0 & $a$ & 1
\end{tabular}
}
\end{picture}
\medskip

\setlength{\unitlength}{1cm} \vspace*{1.5cm}
\begin{picture}(15,2)
\put(1.5,-.4){\circle*{.15}}

\put(1.5,1){\circle*{.15}}

\put(1.5,2.3){\circle*{.15}}

\put(1.1,-.4){$0$}

\put(1.1,.9){$a$}

\put(1.1,2.1){$1$}

\put(.1,.7){${\bf L}_9$ \ :}

\put(1.5,-.4){\line(0,1){2.7}}

\put(2.2,1){\begin{tabular}{c|ccc}
 $\to$ \ & \ 0 \ & \ $a$ \ & \ 1 \\ \hline
  & & & \\
  0 & 1 & 0 & 0 \\
  & & & \\
  $a$ & 0 & 1 & 1 \\
  & & & \\
  1 & 0 & $a$ & 1
\end{tabular}
}

\put(8.2,-.4){\circle*{.15}}

\put(8.2,1){\circle*{.15}}

\put(8.2,2.3){\circle*{.15}}

\put(7.8,-.4){$0$}

\put(7.8,.9){$a$}

\put(7.8,2.1){$1$}

\put(6.7,.7){${\bf L}_{10}$ \ :}

\put(8.2,-.4){\line(0,1){2.7}}

\put(8.9,1){\begin{tabular}{c|ccc}
 $\to$ \ & \ 0 \ & \ $a$ \ & \ 1 \\ \hline
  & & & \\
  0 & 1 & 0 & 0 \\
  & & & \\
  $a$ & 0 & 1 & $a$ \\
  & & & \\
  1 & 0 & $a$ & 1
\end{tabular}
}

\put(6,-1.4){Figure 3} \label{page_algebra2}  
\end{picture}

\vspace{2cm}

Since $0'=1$ and $1'=0$, it is easy to see that there are exactly two expansions on each of the above 10 semi-Heyting algebras by a unary operation $'$ so that the expansions are  $\mathbb{DQDSH}$-algebras.  Ten of these that correspond to $a'=a$ 
are clearly in 
$\mathbb{DMSH}$.    
The other ten, that correspond to $a'=1$ 
are in $\mathbb{DPCSH}$.

To put it more precisely,
let ${\bf L}_i^{dm}$, $i=1,2, \dots, 10$, denote the expansion of ${\bf L}_i$ by adding the unary operation $'$ such that $0'=1$, $1'=0$, and $a'=a$, 
Similarly, let ${\bf L}_i^{dp}$, $i=1,2, \dots, 10$, denote the expansion of ${\bf L}_i$ by adding the unary operation $'$ such that $0'=1$, $1'=0$, and $a'=1$. 
Then, clearly,  ${\bf L}_i^{dm} \in \mathbf {DMSH}$ and  ${\bf L}_i^{dp} \in \mathbf {DPCSH}$.  

Let $\mathbf{C}^{dm} := \{{\bf L}_i^{dm} : i=1,2, \dots, 10\}$,  $\mathbf{C}^{dp} := \{{\bf L}_i^{dm} : i=1,2, \dots, 10\}$ and let $\mathbf{C}_{20} := \mathbf{C}^{dm} \,\cup \, \mathbf{C}^{dp} $.   Thus there are exactly 20 three-element  
$\mathbf{DQDSH}$-algebras, whose lattice reducts are chains.
Let      
$\mathbb {DQDSHC}^{3} :=\mathbb V(\mathbf C_{20})$, the subvariety of $\mathbb {DQDSH}$ generated by all the 20 3-element DQDSH-chains.  Also, let $\mathbb {DMSHC}^{3} :=\mathbb V(\mathbf C^{dm})$ and let $\mathbb {DPCSHC}^{3} :=\mathbb V(\mathbf C^{dp})$.   

We shall now present axiomatizations for the logics corresponding to $\mathbb {DQDSHC}^{3}$, $\mathbb {DMSHC}^{3}$, and $\mathbb {DPCSHC}^{3}$.  


The following theorem is immediate from 
 {\rm \cite[Lemma 10.2, Theorem 10.3, Corollary 10.4 and Theorem 11.1]{sankappanavar2011expansions}}.  Let $x^+ :=x'{^*}'$.  In the rest of the paper, ``equational base'' is abbreviated to ``base''.

\begin{Theorem}  \cite{sankappanavar2011expansions} \label{ThA}  
A base for $\mathbb {DQDSHC}^{3}$, relative to $\mathbb {DQDSH}$, is given by:
\begin{enumerate}
\item[\rm(i)] $x^{**} \approx x{^*}'$,
\item[\rm(ii)] $x \land x^+ \leq y \lor y^*$ \quad \rm(Regularity).
\end{enumerate}
\end{Theorem}

The following theorem follows from
Theorem \ref{teo_040417_01} and Theorem \ref{ThA}.  

\begin{Theorem}
The logic $\mathcal{DQDSHC}^3 $ corresponding to the variety 
$\mathbb {DQDSHC}^{3}$ is defined, as  an extension of $\mathcal {DQDSH}$,
by the following axioms:  	
	\begin{itemize}
	\item $[(\phi \to \bot) \to \bot] \Leftrightarrow_H \ \sim(\phi \to \bot)$,
       \item $[\phi \ \land \ \sim(\sim\phi  \to \bot)] \to_H [\psi \lor (\psi \to \bot)]$.
	\end{itemize}

\end{Theorem}

Since the logic $\mathcal{DQDSHC}^3 $ is finitely axiomatized and the corresponding variety $\mathbb{DQDSHC}^3 =\mathbb{V}(\mathbf C_{20})$ is finitely generated, the following corollary is  immediate.

\begin{Corollary}
The logic $\mathcal{DQDSHC}^3 $ is decidable. 
\end{Corollary}


\begin{Proof}
Observe that $\mathbf{C}_{20} \models x^* \lor x^{**} \approx 1$. 
\end{Proof}

\bigskip

\subsection{Logics $\mathcal {DMSHC}^{3}$ and $\mathcal{DPCSHC}^3$} 
 
 \ \\

We know from Section \ref{section_3valuedlogics} that $\mathbb {DMSHC}^{3}=\mathbb{V}(\mathbf{C^{dm}})$, $\mathbb {DPCSHC}^{3}=\mathbb{V}(\mathbf{C^{dp}})$, ${\bf L}_i^{dm} \in \mathbf {DMSH}$ and  
${\bf L}_i^{dp} \in \mathbf {DPCSH}$,  
$i=1,2, \dots, 10$. 
 
\indent The following theorem is immediate from Theorem \ref{ThA}. 

\begin{Theorem} \label{Theorem_280219_01}  
\begin{thlist}
\item[a] A base for $\mathbb {DMSHC}^{3}$, relative to $\mathbb {DQDSHC}^{3}$, is given by:
\begin{enumerate}
\item[] 
$x'' \approx x$.
\end{enumerate}
\smallskip
\item[b] A base for $\mathbb {DPCSHC}^{3}$, relative to $\mathbb {DQDSHC}^{3}$, is given by:
\begin{enumerate}
\item[] 
$x \lor x'\approx 1$.
\end{enumerate}
\end{thlist}
\end{Theorem}

Let $\mathcal{DMSHC}^3$ and $\mathcal{DPCSHC}^3$  
denote the extensions of 
the varieties 
$\mathbb{DMSHC}^3$ and $\mathbb{DPCSHC}^3$, respectively.
The following theorem is immediate from 
Theorem \ref{teo_040417_01} and Theorem \ref{Theorem_280219_01}.  \\

\medskip
\begin{Theorem}\label{Theorem_280219_01B}   
\begin{thlist}
\item[a]  $\mathcal{DMSHC}^3 $ 
 is defined, as an extension of $\mathcal{DQDSHC}^3 $  
by the following axioms:  
\smallskip	
	\begin{itemize}
	\item[] $\phi \to_H \ \sim\sim\phi$.
	\end{itemize}
\smallskip
\item[b] $\mathcal{DPCSHC}^3$ is defined, relative to the logic $\mathcal{DQDSHC}^3$  
	by the following axiom:  
	\smallskip	 		
	\begin{itemize}
			\item[] $\phi \ \lor \sim\phi $. 			
	\end{itemize}
\end{thlist}
\end{Theorem}

It is clear that the logics $\mathcal{DMSHC}^3$ and $\mathcal{DPCSHC}^3 $ are decidable.

\subsection{3-valued Extensions of $\mathcal{DMSHC}^3$ and of $\mathcal{DPCSHC}^3$}

\

\medskip

\indent We are ready to look at the problem of axiomatization for the logics associated with the 20 3-element chains 
in $\mathbf{C}_{20}$.
We need to recall another (algebraic) result from \cite{sankappanavar2011expansions} that gives a base for each of 3-chains in $\mathbf{C}^{dm}$ and $\mathbf{C}^{dp}$.  To this end, 
we need the
following identities from \cite{sankappanavar2011expansions}:

\medskip

\begin{enumerate}

  \item[(C1)] $x \vee (x \to y) \approx (x \to y)^* \to x$,  
  \item[(C2)] $x \vee [y \to (x \vee y)] \approx (0 \to x) \vee (x \to y)$, 
  \item[(C3)] $x \vee (y \to x) \approx [(x \to y) \to y] \to x$,   
  \item[(C4)] $x \vee (x \to y) \approx x \to [x \vee (y \to 1)]$, 
  \item[(C5)] $(x \to y) \to (0 \to y) \approx x \vee [(x \wedge y) \to 1]$, 
  \item[(C6)] $x^* \vee (x \to y) \approx (x \vee y) \to y$,   
  \item[(C7)] $x \vee (0 \to x) \vee (y \to 1) \approx x \vee [(x \to 1) \to (x \to y)]$,  
  \item[(C8)] $x \vee y \vee (x \to y) \approx x \vee [(x \to y) \to 1]$,  
  \item[(C9)] $x \vee [(0 \to y) \to y] \approx x \vee [(x \to 1) \to y]$, 
  \item[(C10)] $x \vee [x \to (y \wedge (0 \to y))] \approx x \to [(x \to y) \to y]$, 
  \item[(C11)] $(0 \to 1)^* = 0$, 
  \item[(C12)] $x \vee y \vee [y \to (y \to x)] \approx x \to [x \vee (0 \to y)]$, 
  \item[(C13)] $x \vee (x \to y) \approx x \vee [(x \to y) \to 1]$, 
  \item[(C14)] $0 \to 1 \approx 0$ (FTF identity), 
  \item[(C15)] $x \to y \approx y \to x$ (commutative identity). 
\end{enumerate}

In the sequel, we abbreviate
``is a base, relative to $\mathbb {DMSHC}^{3}$ [$\mathbb {DPCSHC}^{3}$]'' to just ``is a base''. \\ %

The reader should keep in mind that the following theorem is really a simultaneous presentation of two separate theorems (in order to keep the size of the paper within limits). One of the two theorems is regarding $\mathbb {DMSHC}^{3}$-algebras and the other is about $\mathbb {DPCSHC}^{3}$-algebras. 
As an illustration, item (i), when decoded, yields the following two (independent) statements:\\

(i$_{dm}$):  \qquad \{(C1)\} is a base,  relative to $\mathbb {DMSHC}^{3}$, for the variety $\mathbf{V(L_1^{dm})}$,  

\begin{center}
and    
\end{center}

(i$_{dp}$):  \qquad \{(C1)\} is a base,  relative to $\mathbb {DPCHC}^{3}$, for the variety $\mathbf{V(L_1^{dp})}$. \\ 

\noindent Similar remark applies to each of the other items of Theorem \ref{3bases} as well. 

Theorem \ref{3bases} is immediate from  {\rm \cite[Theorem 11.2]{sankappanavar2011expansions}}.

\begin{Theorem} \label{3bases} 
\begin{thlist}   
  \item[i]  {\rm \{(C1)\}} is a base for the variety $\mathbf{V(L_1^{dm})}$   {\rm [$\mathbf{V(L_1^{dp})}$]}, 
  \item[ii]  {\rm\{(C2), (C3)\}} is a base for the variety $\mathbf{V(L_2^{dm})}$  {\rm [$\mathbf{V(L_2^{dp})}$]},  
  \item[iii]   {\rm\{(C2), (C4)\}} is a base for the variety $\mathbf{V(L_3^{dm})}$  {\rm [$\mathbf{V(L_3^{dp})}$]},  
  \item[iv]   {\rm\{(C4), (C5)\}} is a base for the variety $\mathbf{V(L_4^{dm})}$  {\rm [$\mathbf{V(L_4^{dp})}$]},  
  \item[v]   {\rm\{(C7)\}} is a base for the variety $\mathbf{V(L_5^{dm})}$  {\rm [$\mathbf{V(L_5^{dp})}$]},
  \item[vi]   {\rm\{(C8)\}} is a base for the variety $\mathbf{V(L_6^{dm})}$   {\rm [$\mathbf{V(L_6^{dp})}$]}, 
  \item[vii]   {\rm\{(C9), (C10)\}} is a base for the variety $\mathbf{V(L_7^{dm})}$   {\rm [$\mathbf{V(L_7^{dp})}$]}, 
  \item[viii]   {\rm\{(C11), (C12)\}} 
 is a base for the variety $\mathbf{V(L_8^{dm})}$   {\rm [$\mathbf{V(L_8^{dp})}$]},
  \item[ix]  {\rm\{(C6), (C13), (C14)\}} is a base for the variety $\mathbf{V(L_9^{dm})}$  {\rm [$\mathbf{V(L_9^{dp})}$]},  
  \item[x]  {\rm\{(C15)\}} is a base for the variety $\mathbf{V(L_{10}^{dm})}$   {\rm [$\mathbf{V(L_{10}^{dp})}$]}.
\end{thlist}
\end{Theorem}

We are now ready to present the axiomatizations for the logics associated with the 20 3-element chains 
in $\mathbf{C}_{20}$.


Let $\mathcal{L}(\mathbf{ L}_i^{dm})$ (or $\mathcal{L}(\mathbb{V}(\mathbf{ L}_i^{dm}))$)  
denote the extension of the logic $\mathcal{DMSHC}^{3}$           
corresponding to the variety $\mathbb{V}(\mathbf{L}_i^{dm})$, 
for $i=1,2, \cdots, 10$.
Also, let $\mathcal{L}(\mathbf{ L_i^{dp}})$ (or $\mathcal{L}(\mathbb{V}(\mathbf{ L}_i^{dp}))$)  
denote the extension of the logic $\mathcal{DPCSHC}^{3}$ corresponding to the variety $\mathbb{V}(\mathbf{L}_i^{dp})$.

In what follows, ``defined, as an extension of the logic  
$\mathcal{DMSHC}_{3}$  {\rm [$\mathcal{DPCSHC}_{3}$]},
by'' is abbreviated to ``defined by''.  The following theorem will follow from
Theorem \ref{teo_040417_01},  Theorem \ref{Theorem_280219_01},  Theorem \ref{Theorem_280219_01B}, and Theorem \ref{3bases}.

\begin{Theorem} 
\begin{thlist}
\item[a] $\mathcal{L}(\mathbf{L}_1^{dm})$  {\rm [$\mathcal{L}(\mathbf{L}_1^{dp})$]} is defined   
by the following axiom: 
\smallskip
\begin{itemize}
	\item[] $[\phi \lor (\phi \to \psi)]  \Leftrightarrow_H [((\phi \to \psi) \to \bot) \to \phi]$. 
   	\end{itemize}
\smallskip
\item[b]$\mathcal{L}(\mathbf{L}_2^{dm})$  {\rm [$\mathcal{L}(\mathbf{L}_2^{dp})$]} is          
defined by the following  axioms:
\smallskip
\begin{itemize} 
     \item[{\rm (i)}] $[\phi \lor \{\psi \to (\phi \vee \psi)\}]  \Leftrightarrow_H [(\bot \to \phi) \vee (\phi \to \psi)]$,	
 \smallskip    
	\item[{\rm (ii)}] $[\phi \vee (\psi \to \phi)] \Leftrightarrow_H [\{(\phi \to \psi) \to \psi\} \to \phi] $.
   	\end{itemize} 
	\smallskip
\item[c] $\mathcal{L}(\mathbf{L}_3^{dm})$  {\rm [$\mathcal{L}(\mathbf{L}_3^{dp})$]} is defined 
by the following axioms:
\smallskip
\begin{itemize} 
      \item[{\rm (i)}]  $[\phi \lor \{\psi \to (\phi \vee \psi)\}]  \Leftrightarrow_H [(\bot \to \phi) \vee (\phi \to \psi)]$, 
	\smallskip
	\item[{\rm (ii)}] $[\phi \vee (\phi \to \psi)] \Leftrightarrow_H [\phi \to \{\phi \vee (\psi\ \to \top)\}] $.
   	\end{itemize}
\smallskip
\item[d] $\mathcal{L}(\mathbf{L}_4^{dm})$  {\rm [$\mathcal{L}(\mathbf{L}_4^{dp})$]} is defined 
by the following axioms:
\smallskip
\begin{itemize} 
	\item[{\rm (i)}] $[\phi \vee (\phi \to \psi)] \Leftrightarrow_H [\phi \to \{\phi \vee (\psi\ \to \top)\}] $,
	\smallskip
      \item[{\rm (ii)}] $[(\phi \to \psi) \to (\bot \to \psi)]  \Leftrightarrow_H [\phi \vee \{(\phi \land \psi) \to \top\}]$. 
   	\end{itemize}
\smallskip
\item[e] $\mathcal{L}(\mathbf{L}_5 ^{dm})$  {\rm [$\mathcal{L}(\mathbf{L}_5^{dp})$]} is defined
by the following axiom:
\smallskip
\begin{itemize} 

	\item[] $[\phi \vee (\bot \to \phi) \vee (\psi \to \top)] \Leftrightarrow_H [\phi \vee \{(\phi \to \top) \to (\phi\ \to \psi)\}] $.
   	\end{itemize}
\smallskip
\item[f] $\mathcal{L}(\mathbf{L}_6^{dm})$  {\rm [$\mathcal{L}(\mathbf{L}_6^{dp})$]} is defined 
by the following axiom:
\smallskip
\begin{itemize} 
	\item[] $[\phi \vee \psi \vee (\phi \to \psi)] \Leftrightarrow_H [\phi \vee \{(\phi \to \psi) \to \top\}] $.
   	\end{itemize}
\smallskip
\item[g] $\mathcal{L}(\mathbf{L}_7^{dm})$  {\rm [$\mathcal{L}(\mathbf{L}_7^{dp})$]} is
defined 
by the following axioms:
\smallskip
\begin{itemize} 
	\item[{\rm (i)}] $[\phi \vee \{(\bot \to \psi) \to \psi\}] \Leftrightarrow_H [\phi \vee \{(\phi \to \top) \to \psi\}] $,
	\smallskip
	\item[{\rm (ii)}] $[(\phi \vee \{\phi \to (\psi \land (\bot \to \psi))\}]  \Leftrightarrow_H [\phi \to \{(\phi \to \psi) \to \psi\}]$. 
   	\end{itemize}
\smallskip
\item[h] $\mathcal{L}(\mathbf{L}_8^{dm})$  {\rm [$\mathcal{L}(\mathbf{L}_8^{dp})$]} is 
defined  
by the following axioms:
\smallskip
\begin{itemize} 
	\item[{\rm (i)}]  $((\bot \to \top) \to \bot)  \to_H \bot  $,
	\smallskip
	\item[{\rm (ii)}] $[(\phi \vee \psi \vee \{\psi \to (\psi \to \phi)\}]  \Leftrightarrow_H [\phi \to \{\phi \vee (\bot \to \psi)\}]$. 
   	\end{itemize}
\smallskip
\item[i] $\mathcal{L}(\mathbf{L}_9^{dm})$  {\rm [$\mathcal{L}(\mathbf{L}_9^{dp})$]} is
defined 
by the following axioms:
\smallskip
\begin{itemize} 
	\item[{\rm (i)}] $[\phi^* \vee (\phi \to \psi)] \Leftrightarrow_H [(\phi \vee \psi) \to \psi] $,
	\smallskip
	\item[{\rm (ii)}] $[\phi \vee (\phi \to \psi)]  \Leftrightarrow_H [\phi \vee \{(\phi \to \psi) \to \top\}]$, 
	\smallskip
	\item[{\rm (iii)}] $(\bot \to \top) \to_H \bot$.
   	\end{itemize}
\smallskip
\item[j] $\mathcal{L}(\mathbf{L}_{10}^{dm})$  {\rm [$\mathcal{L}(\mathbf{L}_{10}^{dp})$]} is
defined 
by the following axioms:
\smallskip
\begin{itemize} 
	\item[] $(\phi\to \psi) \to_H (\psi \to \phi).  $\\
   	\end{itemize}
\end{thlist}
\end{Theorem}



\begin{Remark} Some features of these logics:
\begin{itemize}
\item The logics $\mathcal{L}(\mathbf{ L_i^{dp}}), i \in \{2,3, \dots,10\}$ and the logics $\mathcal{L}(\mathbf{ L_i^{dm}}), i \in \{2,3, \dots,10\}$ are, just 
like the logic $\mathcal{L}(\mathbf{\bar{2}^e})$, ``anti-classical''
in the sense that the classically provable formula $\bot \to \top$ fails in these logics. 

\item Each of the logics $\mathcal{L}(\mathbf{ L_i^{dp}}), i \in \{5,6,7,8\}$ and $\mathcal{L}(\mathbf{ L_i^{dp}}), i \in \{5,6,7,8\}$, just like $\mathcal{L}(\mathbf{2^e})$ and $\mathcal{L}(\mathbf{\bar{2}^e})$, is a coatom in the lattice of extensions of the logic 
$\mathcal{DQDSH}$.  

\item Each of the logics $\mathcal{L}(\mathbf{ L_i^{dp}}), i \in \{1,2,3,4\}$ and $\mathcal{L}(\mathbf{ L_i^{dm}}), i \in \{1,2,3,4\}$  is covered by the coatom  $\mathcal{L}(\mathbf{2^e})$, the classical propositional logic, while each of the logics $\mathcal{L}(\mathbf{ L_i^{dp}}), i \in \{9,10\}$ and $\mathcal{L}(\mathbf{ L_i^{dm}}), i \in \{9,10\}$   is covered by the coatom  $\mathcal{L}(\mathbf{\bar{2}^e})$.     

\item In the logics $\mathcal{L}(\mathbf{ L_i^{dp}}), i \in \{9, 10\}$ and $\mathcal{L}(\mathbf{ L_i^{dm}}), i \in \{9, 10\}$, 
Moreover,  in the logics $\mathcal{L}(\mathbf{ L_{10}^{dp}})$ and $\mathcal{L}(\mathbf{ L_{10}^{dm}})$,
 the connective $\to$ is commutative.  
 
\item The logics  $\mathcal{L}(\mathbf{ L}_i^{dm})$, 
 and
$\mathcal{L}(\mathbf{ L_i^{dp}})$, $i=1,2, \cdots, 10$,   
do not have the (DP) as the formula $\alpha^* \lor \alpha^{**}$ is provable in these logics.

\item The logics $\mathcal{L}(\mathbf{ L_i^{dp}}), i \in \{1,2,3, \cdots,10\}$ and the logics $\mathcal{L}(\mathbf{ L_i^{dm}}), i \in \{1,2,3, \cdots,10\}$ are quasiprimal in the sense that their corresponding varieties are quasiprimal {\rm (\cite{{sankappanavar2011expansions}}}.  

\item Each of the logics $\mathcal{L}(\mathbf{ L_i^{dp}}), i \in \{5,6,7,8\}$ and $\mathcal{L}(\mathbf{ L_i^{dp}}), i \in \{5,6,7,8\}$, just like $\mathcal{L}(\mathbf{2^e})$ and $\mathcal{L}(\mathbf{\bar{2}^e})$, is primal. 

\end{itemize}
\end{Remark}

Further features of some of these logics will be given in {\rm Remark} {\rm \ref{R_Connexive}}.

We note that all the logics mentioned in this subsection are decidable as their corresponding varieties are easily seen to have the finite model property.  

\subsection{$3$-valued {\L}ukasiewicz Logic Revisited} \label{section_Lukasieticz_revisited}
\ 
\medskip

\indent It is worthwhile to point out that the logic $\mathcal{L}\mathbf{(L_1^{dm})}$, defined earlier, has an interesting relationship with the well-known 3-valued {\L}ukasiewicz logic. 
Let us recall the definition of 3-valued {\L}ukasiewicz algebras.

An algebra $\mathbf{A}=\langle A, \lor, \land, ', d_1, d_2, 0,1 \rangle$ is a 3-valued  {\L}ukasiewicz algebras if
\begin{enumerate}[(1)]
\item  $\langle A, \lor, \land, ', 0,1 \rangle$ is a De Morgan algebra,
\item  $d_i (x \lor y)=d_i (x) \lor d_i (y)$, for $i=1,2$,
\item  $d_i (x) \lor  (d_i (x))'=1$, for $i=1,2$,
\item  $ d_i (d_j (x))=d_j (x)$, for $i=1,2$, 
\item  $d_i (x')=(d_{3-i} (x))'$, for $i=1,2$, 
\item  $d_1(x)  \leq d_2(x)$,
\item  If $d_1(x)=d_1(y)$ and $d_2 (x)=d_2 (y)$ then $x=y$.   
\end{enumerate}

Let $\mathbf{{L}} =\langle \{0,a,1\}, \lor, \land, '. d_1, d_2, 0,1 \rangle$ be the algebra such that
$\langle \{0,a,1\}, \lor, \land, ', 0,1 \rangle$ is a 3-element Kleene algebra with $0<a<1$,  and $d_1$ and $d_2$ are unary operations defined as follows: $d_1(0)=d_1 (a) =0, d_1 (1)=1$, and $d_2 (0)=0$, and $d_2 (1)= d_2 (a)=1$.  Then it is routine to verify that $\mathbf{{L}}$ is a 3-valued {\L}ukasiewicz algebra.
It is well-known that $\mathbb{V}\mathbf{({L})}$ is precisely the variety of all 3-valued {\L}ukasiewicz algebras.

\begin{Theorem}
The logic $\mathcal{L}\mathbf{(L_1^{dm})}$ is equivalent to the 3-valued {\rm{\L}ukasiewicz} logic. 
\end{Theorem}

\begin{Proof}
It suffices to prove that the variety $\mathbb{V}\mathbf{(L_1^{dm})}$ is term-equivalent to the variety $\mathbb{V}\mathbf{({\L})}$. 
Without loss of generality, we can assume that $\mathbf{L_1^{dm}}$ and $\mathbf{{\L}}$ have the same universe, say $L= \{0,a,1\}$ with $0 < a <1$.  Given $\mathbf{L_1^{dm}}$,
define the unary operations $d_1$ and $d_2$ on L  
by: $d_1(x)=x'^*$ and $d_2= x{^*}'$.  Then it is straightforward to verify that $\langle L; \lor,  \land, ', d_1, d_2, 0, 1 \rangle =  \mathbf{{\L}}$.  To prove the converse, let us first define the unary function $^*$ on $\mathbf{L}$ by:  $x^* := d_1((d_2(x))')$.  It is routine to verify that $^*$ is the pseudocomplement operation on $L$.  Using $^*$ we can now define the Katri\v{n}\'{a}k's implication $\to$ by:
$$x \to y := (x^* \lor y^{**}) \land [(x \lor x^*)'{^*}' \lor x^* \lor y \lor y^*].$$ 
\smallskip
Then, $\to$ is the Heyting implication (see \cite{Ka73}).  Hence, it follows that $\langle L;  \lor, \land, \to, ',  0, 1 \rangle =  \mathbf{L_1^{dm}}$.
The theorem is now proved.
\end{Proof}

\subsection{4-valued Extensions of  $\mathcal{DQDSH}$ with Boolean semi-Heyting Reducts} \label{SecEleven}
\ 
\medskip

Recall that the variety $\mathbb {DQDSH}$  was defined in Section \ref{SecFive}.
 An algebra $\mathbf{L}$  
 is a {\it dually quasi-De Morgan Boolean semi-Heyting algebra} 
 ($\mathbb {DQDBSH}$-algebra, for short) 
if its term-reduct \\ 
$\langle L, \lor, \land, ^*,0,1\rangle $ is a Boolean semi-Heyting algebra, that is, $\mathbf{L} \models x \lor x^* \approx 1$.  The variety of such algebras is denoted by $\mathbb {DQDBSH}$.  \\

The following theorem is now immediate, in view of Theorem \ref{teo_040417_01}. 

\begin{Theorem}
The logic $\mathcal{DQDBSH}$ is  
 defined, relative to $\mathcal{DQDSH}$ by the following
 axiom:  
\begin{itemize}
\item[\rm{(B)}] $\phi \lor (\phi \to \bot).$ 
\end{itemize}
\end{Theorem}

The concrete description of the lattice of subvarieties of $\mathbb {DQDBSH}$ was given in  \cite{sankappanavar2011expansions}.  We now wish to present the axiomatizations for corresponding extensions of the logic $\mathcal{DQDBSH}$. Toward this end, the following three algebras will be needed.

Figure 3 defines the $\to$ operation on the three $\rm{4}$-element algebras $\mathbf{D_1}$, $\mathbf{D_2}$ and $\mathbf{D_3}$, each of whose lattice reduct is the 4-element Boolean lattice
having the universe $\{0,a,b,1\}$, with $b$ as the complement of $a$, and $'$ is defined as follows: $a' = a$, $b' = b$, $0'=1$ and $1'=0$. \\
\vspace*{1.5cm}

\setlength{\unitlength}{.9cm}         
\begin{picture}(15,4)

\put(1,2.5){$\bf D_1$ \ :}

\put(8.5,2.5){$\bf D_2$ \ :}

\put(2.2,2.5){\begin{tabular}{c|cccc}
 $\to$ \ & \ 0 \ & \ $1$ \ & \ $a$\ &\ $b$ \\ \hline
  & & & &\\
  0 & 1 & 0 & $b$ & $a$ \\
  & & & &\\
  $1$ & 0 & 1 & $a$ & $b$ \\
  & & & &\\
  $a$ & $b$ & $a$ & 1 &0\\
  & & & &\\
  $b$ & $a$ & $b$ & 0 & 1\\
\end{tabular}
}
\put(9.5,2.5){\begin{tabular}{c|cccc}
 $\to$ \ & \ 0 \ & \ $1$ \ & \ $a$\ &\ $b$ \\ \hline
  & & & &\\
  0 & 1 & 1 & 1 & 1 \\
  & & & &\\
  $1$ & 0 & 1 & $a$ & $b$ \\
  & & & &\\
  a & $b$ & $1$ & 1 & $b$ \\
  & & & &\\
  $b$ & $a$ & 1 & $a$ & 1\\
\end{tabular}
}
 \end{picture}  

\hspace*{2cm}                          \setlength{\unitlength}{.5cm} 
                           \begin{picture}(15,8)       
                            \put(.4,2.5){$\bf D_3$ \ :}
\put(2.2,2.5){\begin{tabular}{c|cccc}                  
 $\to$ \ & \ 0 \ & \ $1$ \ & \ $a$\ &\ $b$ \\ \hline
  & & & &\\
  0 & 1 & $a$ & 1 & $a$ \\
  & & & &\\
  $1$ & 0 & 1 & $a$ & $b$ \\
  & & & &\\
  $a$ & $b$ & $a$ & 1 &0\\
  & & & &\\
  $b$ & $a$ & 1 & $a$ & 1\\
\end{tabular}
}
\put(3.6,-3.8){Figure 4}    
\end{picture}

\vspace{3cm}    

The algebras $\mathbf{D_1}$, $\mathbf{D_2}$, and $\mathbf{D_3}$ 
are the only simple (=subdirectly irreducible) algebras in $\mathbb {DQDBSH}$.

The following theorem, which follows immediately from \cite[Corollary 9.4]{sankappanavar2011expansions}, reveals the structure of $\mathbb {DQDBSH}$. 

\begin{Theorem}\label{DMBSH} 
$\mathbb {DQDBSH} = \mathbb V (\mathbf{D_1, D_2, D_3}) = \mathbb {DMBSH}$.
\end{Theorem}
  
The above theorem leads us to the following decidability result, in view of Theorem \ref{teo_040417_01}.

\begin{Corollary}
The logic $\mathcal{DQDBSH}$ is decidable.
\end{Corollary} 

We will now turn our attention to the axiomatization of logics corresponding to the varieties generated by these algebras.   
The following theorem is taken from {\rm \cite[Theorem 9.5]{sankappanavar2011expansions}.}
\begin{Theorem} \label{bases}  
\begin{thlist} 
\item[1]  A base for the variety $\mathbb V (\mathbf{D_1})$,  modulo $\mathbb{DQDBSH}$, is given by\\
$0 \to 1 \approx 0$,

\item[2]  A base for $\mathbb V (\mathbf{D_2})$, modulo $\mathbb{DQDBSH}$, is given by\\
$0 \to 1 \approx 1$,
\item[3]  A base for the variety $\mathbb V (\mathbf{D_3})$, modulo $\mathbb{DQDBSH}$, is given by\\
$(0 \to 1)' \approx 0 \to 1$.\\
\end{thlist}
\end{Theorem}

The following corollary will now follow as an application of
Theorem \ref{teo_040417_01}, Theorem \ref{DMBSH} and Theorem \ref{bases}.

Let $\mathcal{L}(\mathbf{D_i})$ (or $\mathcal{L}((\mathbb{V}(\mathbf{D_i})))$)
denote the extension of the logic $\mathcal{DMBSH}$ 
corresponding to the variety $\mathbb{V}(\mathbf{D_i})$ 
for $\rm i =1,2,3$.  

In the rest of this section, ``defined, relative to the logic  $\mathcal{DMBSH}$, by'' is abbreviated to ``defined by''.

\begin{Corollary}    
\begin{thlist}       
    \item[1] The logic $\mathcal{L}(\mathbf{D_1})$ is 
    defined by the axiom: \begin{itemize}
    		\item[] $(\bot \to \top) \to_H \bot$,
                 \end{itemize}
                 
     \item[2] The logic $\mathcal{L}(\mathbf{D_2})$ is 
     defined by the axiom: 
     \begin{itemize}
               \item[]   $ \bot \to \top$.
                \end{itemize}
       \item[3] The logic $\mathcal{L}(\mathbf{D_3})$ is 
       defined by the axiom: 
        \begin{itemize} 
               \item[] $\sim(\bot \to \top) \Leftrightarrow_H (\bot \to \top)$. 
         \end{itemize}                  
\end{thlist}  
\end{Corollary}

It is clear that the logics $\mathcal{L}(\mathbf{D_i}), \rm i \in \{1,2,3\}$, are decidable. 

\begin{Remark}  Some features of these logics:
\begin{thlist}      
\item  The logics $\mathcal{L}(\mathbf{D_i}), \rm i \in \{1,3\}$, 
 are anti-classical since the formula $\bot \to \top$ is not provable in each of them.  

\item  The logics $\mathcal{L}(\mathbf{D_1})$ and $\mathcal{L}(\mathbf{D_2})$, 
are covered, respectively, by the coatoms $\mathcal{L}(\mathbf{\bar{2}^e})$ and $\mathcal{L}(\mathbf{2^e})$ in the lattice of extensions of the logic $\mathcal{DQDSH}$.      

\item In the logic 
$\mathcal{L}(\mathbf{D_1})$,
the connective, $\to$, is commutative.  

\item The logics $\mathcal{L}(\mathbf{D_1})$ and $\mathcal{L}(\mathbf{D_2})$ are quasiprimal, in the sense that their corresponding varieties are generated by quasiprimal algebras $\mathbf{D_1}$ and $\mathbf{D_2}$ respectively.

\item  The logic $\mathcal{L}(\mathbf{D_3})$, just like $\mathcal{L}(\mathbf{2^e})$ and $\mathcal{L}(\mathbf{\bar{2}^e})$, is a coatom in the lattice of extensions of the logic $\mathcal{DMSH}$ and hence, of $\mathcal{DHMSH}$.

\item  The logic $\mathcal{L}(\mathbf{D_3})$ is primal.  


\end{thlist}   
\end{Remark}

 Further features of some of these logics 
 will be given in {\rm Remark} {\rm \ref{R_Connexive}}.  

\ \\
\bigskip

\section{Connection to Connexive Logics} \label{S8}
\ 

\medskip

Let {\bf L} be a language containing the connective symbols: $\to$ for implication and $\neg$ for negation.
 A logic $\mathcal L$ in {\bf L}  is a {\it connexive logic} (see \cite{Wa20}, for example) if the following Aristotle's Theses and Boethius' Theses are theorems in $\mathcal L$:\\

\noindent Aristotle's Theses:\\ 
\indent  (AT)  \quad $\neg(\neg \alpha \to \alpha)$, \\
\indent  (AT')  \quad $\neg(\alpha \to \neg \alpha)$,\\

\noindent Boethius' Theses:\\
\indent (BT)  \quad $(\alpha \to \beta) \to \neg(\alpha \to \neg \beta)$, \\
\indent  (BT') \quad $(\alpha \to \neg \beta) \to \neg(\alpha \to \beta).$\\

For more details on the motivation, the origin and the history of connexive logics, see \cite{Wa20} and \cite{JaMa19}.

Many of the extensions of the logics $\mathcal{SH}$ and $\mathcal{DHMSH}$, to our surprise, turn out to be connexive logics with $\neg \alpha := \alpha \to \bot$.  We present a few of these below.  (More will be said in the paper \cite{CoSa22} which is in preparation.)

Recall that $\mathbf{\bar{2}}$, $\mathbf{ L_9}$, and $\mathbf{ L_{10}}$, are in $\mathbb{SH}$.
 
\begin{Remark} \label{R_Connexive} 
\begin{thlist}  
\item The logics $\mathcal{L}(\mathbf{\bar{2}})$ and $\mathcal{L}(\mathbf{ L_9})$, 
which are extensions of the semi-intuitionistic logic $\mathcal{SI}$, are connexive logics since the corresponding varieties $\mathbb{V}(\mathbf{\bar{2}})$ and $\mathbb{V}(\mathbf{ L_9})$  
satisfy the following identities: 
      \begin{itemize}
      \item[\rm {i)}] $(x^* \to x)^* \approx 1$, 
      \item[{\rm (ii)}] $(x \to  x^*)^* \approx 1$,
      \item[{\rm (iii)}] $(x \to y) \to (x \to y^*)^* \approx 1$,
      \item[{\rm (iv)}] $(x \to y^*) \to (x \to y)^* \approx 1$.
      \end{itemize}
\item {\rm(AT)} and {\rm(AT')} are  theorems in the logic $\mathcal{L}(\mathbf{ L_{10}})$, since the corresponding variety $\mathbb{V}(\mathbf{ L_{10}})$ satisfies the identities {\rm(i)} and {\rm(ii)}, while {\rm(BT)} and {\rm(BT')}  are not theorems in the logics $\mathcal{L}(\mathbf{ L_{10}})$. 
\end{thlist} 
\end{Remark} 

\begin{Remark} \label{R_Connexive1} 
\begin{thlist}
\item The logics $\mathcal{L}(\mathbf{\bar{2}^e})$, $\mathcal{L}(\mathbf{ L_9^{dp}})$, $\mathcal{L}(\mathbf{ L_9^{dm}})$ and $\mathcal{D}_1$, which are extensions of the logic $\mathcal{DHMSH}$,
are connexive logics 
since it is easily verified that their corresponding varieties $\mathbb{V}(\mathbf{\bar{2}^e})$, $\mathcal{L}(\mathbf{ L_9^{dp}})$, $\mathcal{L}(\mathbf{ L_9^{dm}})$,  and $\mathbb V(\mathbf{D}_1)$ satisfy the identities {\rm(i)--(iv)}.

\item {\rm(AT)} and {\rm(AT')} are  theorems in the logics $\mathcal{L}(\mathbf{ L_{10}^{dp}})$ and $\mathcal{L}(\mathbf{ L_{10}^{dm}})$, since the corresponding varieties $\mathbb{V}(\mathbf{ L_{10}^{dp}})$ and $\mathbb{V}(\mathbf{ L_{10}^{dp}})$ satisfy the identities {\rm(i)} and {\rm(ii)}.
\item  {\rm(BT)} and {\rm(BT')}  are not theorems in the logics $\mathcal{L}(\mathbf{ L_{10}^{dp}})$ and $\mathcal{L}(\mathbf{ L_{10}^{dm}})$.
\end{thlist}         
\end{Remark}

Jarmuzek and Malinowski \cite{JaMa19} have recently introduced the notion of a ``quasi-connexive'' logic.  
A logic is quasi-connexive iff it is not connexive, but at least one of (AT), (AT'), (BT) and (BT') is a theorem in the logic. 
Thus, in view of the above remark, the logics 
$\mathcal{L}(\mathbf{ L_{10}^{dp}})$ and $\mathcal{L}(\mathbf{ L_{10}^{dm}})$, as well as the extension $\mathcal{L}(\mathbf{ L_{10}})$ of $\mathcal{SH}$, can be viewed as quasi-connexive logics.\\

We now mention a few facts about the relationships among the Aristotle's Theses and Boethius' Theses in the logics $\mathcal{SH}$ and $\mathcal{DHMSH}$ whose proofs will appear in a forthcoming paper.

\begin{Theorem} \label{Theorem_AT_BT_logics} In the logic $\mathcal{SH}$, and hence in $\mathcal{DHMSH}$,\\
{\rm(a)} {\rm(AT)} and {\rm(AT')} are equivalent.\\
{\rm(b)}  {\rm(AT)}, {\rm(AT')} and {\rm(BT')} are provable from {\rm(BT)}.\\
{\rm(c)}  {\rm(AT)}, {\rm(AT')} are provable from {\rm(BT')}, but {\rm(BT)} is not.
\end{Theorem}

\begin{Remark}
We propose that any logic in which the (classically provable) formula  $\bot \to \top$ is not provable be included in the family of connexive logics since such a logic would be not only anti-classical but also anti-intuitionistic logic.  Accordingly, the logics $\mathcal{L}(\mathbf{\bar{2}^e})$, $\mathcal{L}(\mathbf{ L_i}^{dp}), i = 5, \dots, 10$,  $\mathcal{L}(\mathbf{ L_i}^{dm}), i = 5, \dots, 10$, $\mathbb V(\mathbf{D}_1)$ and $\mathbb V(\mathbf{D}_3)$ can be considered as connexive logics. 
\end{Remark}

\vspace{.5cm}
\section{Two Infinite Chains of 
Extensions of the logic $\mathcal{DQDH}$}  \label{section_infinite_chains}
\ 

\indent Recall that the logic $\mathcal{DQDH}$ corresponds to the variety of dually quasi-De Morgan Heyting algebras. 
In this section, we present two infinite chains of logics that are extensions of the logic $\mathcal{DQDH}$. \\

\medskip
\subsection{De Morgan-G\"{o}del logic and its extensions: Logics corresponding to the subvarieties of the variety generated by all De Morgan G\"{o}del algebras}      

\ \\

Recall that the variety $\mathbb{DMH}$ of De Morgan Heyting algebras is the subvariety of $\mathbb{DQDH}$ defined by the axiom: $x'' \approx x$.  A De Morgan Heyting algebra whose lattice reduct is a chain is called a \emph{De Morgan Heyting  chain}.  
Let $\mathbb{DMG}$ (or $\mathbb{DMHC}$) denote the subvariety of $\mathbb{DMSH}$ generated by the De Morgan Heyting  chains.  It is proved in
 \cite[Theorem 12.5]{sankappanavar2011expansions} 
 that the lattice of subvarieties of the variety 
$\mathbb{DMG}$ is an $\omega +1$-chain.  Let $\mathcal{DMG}$ (or $\mathcal{DMHC}$) denote the extension of the logic $\mathcal{DMH}$, corresponding to $\mathbb{DMG}$.  We will refer to the logic $\mathcal{DMG}$ as ``De Morgan-G\"{o}del logic.''  
Then it follows 
that the lattice of extensions of $\mathcal{DMG}$ is a chain dual to $\omega +1$.

    In this subsection, we present axiomatizations for the logics corresponding to the subvarieties of 
$\mathbb{DMG}$.  For this purpose, we need  
the following algebraic result which was proved in 
  \cite[Theorem 12.3]{sankappanavar2011expansions}.\\

   \begin{Theorem}  \cite{sankappanavar2011expansions} 
A base for $\mathbb{DMG}$, relative to $\mathbb{DMSH}$, is given by:
\begin{itemize}
\item[{\rm (1)}] $x{^*}' \approx x^{**}$,
\item[{\rm (2)}]  $(x \to y) \lor (y \to x) \approx 1$.
\end{itemize}
\end{Theorem}

Hence we have the following axiomatization for the logic $\mathcal{DMG}$, relative to the logic $\mathcal{DMSH}$.

\begin{Corollary}
The logic $\mathcal{DMG}$, relative to the logic $\mathcal{DMSH}$ is defined by :  
\begin{itemize}
	\item[\rm(i)] $\sim(\alpha  \to \bot) \Leftrightarrow_H ((\alpha  \to \bot) \to \bot)$,
	\item[(ii)] $(\alpha \to \beta) \lor (\beta \to \alpha)$.  
\end{itemize}
\end{Corollary}

In view of the axiom (ii)   
it is clear that the logic $\mathcal{DMG}$ does not have the Disjunction Property.  

\medskip
Next we will present an axiomatization for the logic $\mathcal{DMG}_n$ corresponding to 
each subvariety $\mathbb{DMG}_n$ of $\mathbb{DMG}$ generated by the $n$-element $\mathbb{DMH}$-chain, where $n \in \mathbb{N}$ with $n \geq 2$.  

\begin{Theorem} {\rm \cite{sankappanavar2011expansions} \label{chain}}
Let $n \in \omega$ such that $n \geq 2$.  Then
 $\mathbb{DMG}_n$ is defined, mod $\mathbb{DMG}$, by 
the following axiom:\\
 $$  (\bigvee_{i=1}^{i=n} x_i ) \lor [\bigvee_{i=1}^{i=n-1} (x_i \to x_{i+1})] \approx 1.$$
\end{Theorem}

Hence we have the following axiomatization of the logic $\mathcal{DMG}_n$.  

\begin{Corollary}
Let $n \in \omega$ such that $n \geq 2$.  Then the logic $\mathcal{DMG}_n$, relative to the logic $\mathcal{DMG}$, is defined by\\
$$(\bigvee_{i=1}^{i=n} \alpha_i ) \lor [\bigvee_{i=1}^{i=n-1} (\alpha_i \to \alpha_{i+1})].$$
\end{Corollary}

\medskip



\bigskip

\subsection{{\bf Dually pseudocomplemented G\"{o}del Logic and its extensions: Logic corresponding to the variety generated by dually pseudocomplemented Heyting chains}}\label{11.2}

\ \\

A $\mathbb{DPCSH}$-algebra $\mathbf{L} = \langle L, \land, \vee, \to, ', 0, 1 \rangle$, 
whose lattice reduct is a chain,  
is called a $\mathbb{DPCSH}${\it-chain}.   Let $\mathbb{DPCG}$ (or $\mathbb{DPCHC}$) denote the subvariety of $\mathbb{DPCH}$ generated by the $\mathbb{DPCH}$-chains.  Observe that $\mathbb{DPCHC} = \mathbb{DS}\rm t\mathbb{HC}$.  
It was implicit in 
 \cite[Section 13]{sankappanavar2011expansions} that the lattice of subvarieties of the variety 
$\mathbb{DPCG}$  
is an $\omega +1$-chain 
and was explicitly proved in  \cite[Theorem 4.7]{sankappanavar2017}.   We let $\mathcal{DPCG}$ (or $\mathcal{DPCHC}$) denote the extension of the logic $\mathcal{DPCH}$ corresponding to $\mathbb{DPCG}$.  The logic $\mathcal{DPCG}$ will be referred to as ``dually pseudocomplemented G\"{o}del logic''. 
 It follows from the just mentioned algebraic result 
that the extensions of $\mathcal{DPCG}$ form a chain dual to $\omega +1$.



In this subsection, we present axiomatizations for the logics corresponding to the subvarieties of 
$\mathbb{DPCG}$.  For this purpose, we need the following algebraic result which
was proved in \cite[Theorem 13.2]{sankappanavar2011expansions}.  Let $x^+ := x'{^*}'$.

\begin{Theorem}
The following identities form an equational base, mod $\mathbb{DQDSH}$, for $\mathbb{DPCG}$.
\begin{enumerate}
\item[{\rm(i)}] $x^+\approx x'$, 
\item[{\rm(ii)}] $(x \to y) \lor (y \to x) \approx 1$.
\end{enumerate}
\end{Theorem}

\begin{Corollary}
The logic $\mathcal {DPCG}$ corresponding to the variety $\mathbb{DPCG}$ is defined, as an extension of the logic $\mathcal {DQDSH}$ by

\begin{enumerate}
\item[{\rm(i)}] $\alpha^+ \Leftrightarrow_H \ \sim\alpha$, where $\alpha^+:= \ \sim(\sim\alpha \to \bot)$.
\item[{\rm(ii)}] $(\alpha \to \beta) \lor (\beta \to \alpha) $. 
\end{enumerate}
\end{Corollary}

In view of the axiom (ii)  
it is clear that the logic $\mathcal{DPCG}$ does not have the Disjunction Property. 

Let $n \in \omega$ such that $n \geq 2$ and let $\mathbb{DPCG}_n$ denote the variety generated by the $n$-element $\mathbb{DPCH}$-chain.
 The following theorem, which follows from \cite[Theorem 13.3]{sankappanavar2011expansions}, gives an equational base for each subvariety $\mathbb{DPCG}_n$ of $\mathbb{DPCG}$.   

\begin{Theorem}
Let $n \in \omega$ such that $n \geq 3$.  Then,
$\{{\rm  (An)}\}$ is an equational base, mod $\mathbb{DPCG}$, for $\mathbb{DPCG}_n$, 
where  
{\rm (An)} is the following axiom:

  $$\bigvee_{j=1}^{j=n} x_j  \lor \bigvee_{j=1}^{j=n-1} (x_j \to x_{j+1}) \approx 1.$$
\end{Theorem}
\begin{Corollary}
Let $n \in \omega$ such that $n \geq 2$.  Then the logic $\mathcal {DPCG}_n$ corresponding to the variety $\mathbb{DPCG}_n$ is defined, relative to the logic $\mathcal {DPCG}$ by
\begin{enumerate}
\item[\rm{($\Lambda_n$)}]  $\bigvee_{j=1}^{j=n}  \alpha_j  \lor \bigvee_{j=1}^{j=n-1} (\alpha_j \to \alpha_{j+1})$.     
\end{enumerate}
\end{Corollary}



\medskip
\section{Logics corresponding to subvarieties of regular dually quasi-De Morgan Stone semi-Heyting algebras} 
\label{Sec_Rdqdstsh}
In the rest of the paper we will give axiomatizations for more new logics that are extensions of $\mathcal{DQDSH}$, as applications of Theorem \ref{teo_040417_01}, 
and the algebraic results from \cite{sankappanavar2014a, sankappanavar2014b,  sankappanavar2016, sankappanavar2017, sankappanavar2019}). 
Recall from Section \ref{S7}.2 that 
 $\mathbf{C}_{20}= \mathbf{C}^{dm} \cup  \mathbf{C}^{dp}$, where  
$\mathbf{C}^{dm} := \{{\bf L}_i^{dm} : i=1,2, \dots, 10\}$,  $\mathbf{C}^{dp} := \{{\bf L}_i^{dm} : i=1,2, \dots, 10\}$ and 
 that the algebras ${\mathbf L}_i^{dm}$, ${\mathbf L}_i^{dp}$ were defined in Section \ref{S8}.2  and the three $4$-element algebras $\mathbf{D}_1$, $\mathbf{D}_2$ and $\mathbf{D}_3$  were defined in 
Section \ref{SecEleven}.  
Recall also that      
 $\mathbb V(\mathbf C_{20})$ is the subvariety of $\mathbb {DQDSH}$ generated by all the 20 3-element DQDSH-chains.  Let $\mathbb {DQDSHC}^{3} :=\mathbb V(\mathbf C_{20}).$ \\
 
The notion of regularity has played an important role in \cite{Va71, adams2019, sankappanavar2011expansions, sankappanavar2014a, sankappanavar2014b, sankappanavar2016}.\\

An algebra $\mathbf{A} \in \mathbb{DQDSH}$ is called {\it regular} (\cite{sankappanavar2014b}) if $\mathbf{A}$ satisfies:
\begin{equation} \tag{R}
  \qquad x \land x^+ \leq y \vee y^*,$$
\end{equation}
where $x^+ := x'{^*}'$.\\

The subvariety of $\mathbb{DQDSH}$ of regular algebras is denoted by $\mathbb{RDQDSH}$.   
(We caution the reader that  the term ``regular'' was used in \cite{sankappanavar2011expansions} to mean something else.) 


Observe from Theorem \ref{ThA} that $\mathbb {DQDSHC}^{3} \subset \mathbb {RDQDSH}$.\\

The concept of \emph {level} has played an important role in finding discriminator subvarieties of $\mathbb{DQDSH}$ (see \cite[Corollary 8.2]{sankappanavar2011expansions}).  Here we only define 
$\mathbb{DQDSH}$-algebras of level 1.\\

 An algebra $\mathbf{A} \in \mathbb{DQDSH}$ is of {\it level }$1$ if $\mathbf{A}$ satisfies:
 $$x \land x'^* \approx x \land x'^* \land x'{^*}'^*.$$
 For the varieties of level 1 considered in the rest of the paper, the above definition of ``level 1'' is equivalent to the following:
$$x \land x'^* \approx (x \land x'^*)'^*.$$
Let  $\mathbb{DQDSH}_1$ denote the variety of $\mathbb{DQDSH}$-algebras of level 1.   
Let  $\mathbb{DQDS}\rm t\mathbb{SH}$ denote the subvariety of $\mathbb{DQDSH}$
  that satisfies the Stone identity:\\ 
\begin{equation} \tag{St}   
  \qquad x^\ast \  \vee \ x^{\ast \ast} \approx 1.
\end{equation}  
 $\mathbb{DQDS}\rm t\mathbb{SH}_1$ denotes the subvariety of $\mathbb{DQDS}\rm t\mathbb{SH}$ of level 1, while  $\mathbb{RDQDS}\rm t\mathbb{SH}_1$ denotes the subvariety 
 of $\mathbb{DQDS}\rm t\mathbb{SH}_1$ defined by (R).
 
In this section we present axiomatizations for new logics corresponding to several subvarieties of the variety $\mathbb{RDQDS}\rm t\mathbb{SH}_1$ of regular dually quasi-De Morgan Stone semi-Heyting algebras of level 1.

{\bf In what follows,  $\mathcal{V}$ (or $\mathcal{L}(\mathbb{{}V}))$ denotes the logic corresponding to the subvariety $\mathbb{V}$
of $\mathbb{DQDSH}$-algebras.} 

 (Thus, for example, the logic $\mathcal{DQDS}\rm{t}\mathcal{SH}_1$ corresponds to the variety $\mathbb{DQDS}\rm{t}\mathbb{SH}_1$.) 
 
The following corollary is immediate from the above definitions and Theorem \ref{teo_040417_01}.

\begin{Corollary} 
\begin{thlist}
\item[a] The logic $\mathcal{DQDS}\mathbf{t}\mathcal{SH}_1$ 
is defined, as an extension of 
the logic $\mathcal{DQDSH}$, by the following axioms:
\begin{itemize}
\item[{\rm(1)}]  $[\sim\{(\alpha \wedge \ (\sim\alpha)^*\}]^* \Leftrightarrow_H \ \alpha \wedge \ (\sim\alpha)^*$, 
\item[{\rm(2)}] $\alpha^* \lor \alpha^{**}.$  
\end{itemize}
\smallskip
\item[b] The logic $\mathcal{RDQDS}\rm{t}\mathcal{SH}_1$ 
is defined, as an extension of 
the logic $\mathbb{DQDS}\rm{t}\mathbb{SH}_1$ by the following axiom:
\begin{itemize} 
\item[]  $(\alpha \land \alpha^+) \ \to_H \ (\beta \lor \beta^*) $.
\end{itemize}
\end{thlist}
\end{Corollary}



The following result is taken from \cite[Theorem 3.1]{sankappanavar2014b}. 

\begin{Theorem} \label{TH_Main} 

 $\mathbb{RDQDS}\rm t\mathbb{SH}_1
                                         =  \mathbb{V}(\mathbf{C}_{20} \cup \{\mathbf{D}_1, \mathbf{D}_2, \mathbf{D}_3\}). $
In particular, \\
$\mathbb{RDQDS}{\rm t}\mathbb{H}_1 = \mathbb{V}(\{\mathbf{L}^{dm}_1, \mathbf{L}^{dp}_1, \mathbf D_2\}).$ 
\end{Theorem}  

The following corollary is immediate from Theorem \ref{TH_Main} and Theorem \ref{teo_040417_01}, as the variety  $\mathbb{RDQDS}\rm t\mathbb{SH}_1$ is finitely axiomatized and is generated by a finite set of finite algebras.  

\begin{Corollary}
 The logics $\mathcal{RDQDS}\rm{t}\mathcal{SH}_1$ and $\mathcal{RDQDS}\rm{t}\mathcal{H}_1$ are decidable. 
\end{Corollary}

 In view of the above corollary, 
 it would be of interest to know if 
the logic $\mathcal{DQDS}\rm{t}\mathcal{SH}_1$ is decidable; 
in particular, if the logic $\mathcal{DQDS}\rm{t}\mathcal{H}_1$ is decidable. 
This naturally leads us to the following open problem.\\

\noindent {\bf PROBLEM:} Is the variety $\mathbb{DQDS}\rm{t}\mathcal{SH}_1$ (or even $\mathbb{DQDS}\rm{t}\mathbb{H}_1$) generated by its finite members?





\begin{Remark}
It was shown in \cite{sankappanavar2011expansions} that the variety $\mathbb{RDQDS}\rm t\mathbb{SH}_1$ is a discriminator variety.  Thus  $\mathcal{RDQDS}\rm t\mathcal{SH}_1$ is a discriminator logic.  
\end{Remark}

Let  $\mathbb{RDMS}\rm t\mathbb{SH}_1$ and  $\mathbb{RDMS}\rm t\mathbb{H}_1$ denote, respectively, the varieties of regular De Morgan Stone semi-Heyting algebras and regular De Morgan Stone semi-Heyting algebras.  Similarly, the varieties  
$\mathbb{RDPCS}\rm t\mathbb{SH}_1$ and $\mathbb{RDPCS}\rm t\mathbb{H}_1$ denote, respectively,
the varieties of regular dually pseudocomplemented Stone semi-Heyting algebras and regular dually pseudocomplemented Stone Heyting algebras.   Note that all these varieties are subvarieties of  $\mathbb{RDQDS}\rm t\mathbb{SH}_1$.  \\

Recall $\mathbb {DMSHC}^{3} =\mathbb V(\mathbf C^{dm})$ and let $\mathbb{DPCSHC}^3 =\mathbb V(\mathbf C^{dp})$.   
Observe that $\mathbb {DMSHC}^{3} \subset \mathbb{RDMSH}$ and $\mathbb {DPCSHC}^{3} \subset \mathbb{RDPCSH}$.  \\   

The following corollary is immediate from Theorem \ref{TH_Main}, where ``is defined by'' means 
``is defined, as an extension of $\mathcal{RDQDS}{\rm t}\mathcal{SH}_1$, by ''.

\begin{Corollary} 
\begin{thlist}
\item[a] The logic $\mathcal{RDMS}\rm t\mathcal{SH}_1$ is defined 
by 
\begin{itemize}
\item[]  $\alpha \ \to_H \alpha''$.
\end{itemize}
\item[b]  The logic $\mathcal{RDPCS}\rm t\mathcal{SH}_1$ is defined 
by 
\begin{itemize}
\item[]  $\alpha \ \lor \  \alpha'$.
\end{itemize}
\item[c] The logic $\mathcal{RDMS}\rm t\mathcal{H}_1$ is defined 
by  
\begin{itemize}
\item[]   $(\alpha \land \beta) \to  \alpha$.
\end{itemize} 
\item[d]
The logic $\mathcal{RDPCS}\rm t\mathcal{H}_1$ is defined 
by 
\begin{itemize}
\item[]   $(\alpha \land \beta) \to  \alpha$.
\end{itemize}
\end{thlist}
\end{Corollary}

The following theorem was recently proved in  \cite{sankappanavar2019}.

\begin{Theorem}{\rm \cite[Corollary 3.4]{sankappanavar2019}} \label{Th11.A}
\quad $\mathbb{DMSH}_1= \mathbb{DMS}\rm t\mathbb{SH}_1$.\\
In particular, 
$\mathbb{RDMSH}_1= \mathbb{RDMS}{\rm t}\mathbb{SH}_1.$
\end{Theorem}


The following theorem is immediate from Theorem \ref{Th11.A} and \cite[Corollary 3.4]{sankappanavar2014b}. 
\medskip
\begin{Theorem}  \label{Theo_050719_01}   
\begin{thlist}
\item[a]
 $ \mathbb{RDMSH}_1 = \mathbb{RDMS}{\rm t}\mathbb{SH}_1 =  \mathbb{V}(\mathbf{C}^{dm}) \lor \mathbb{V}
 (\{\mathbf{D_1}, \mathbf{D_2}, \mathbf{D_3}\})$,
\item[b]  $ \mathbb{RDMH}_1 = \mathbb{RDMS}{\rm t}\mathbb{H}_1 =  \mathbb{V}(\{\mathbf{L}^{dm}_1, \mathbf{D}_2\})= \mathbb{V}(\mathbf{L}^{dm}_1) \lor \mathbb{V}(\mathbf{D}_2)$,
\item[c]
 $\mathbb{RDPCS}\rm t\mathbb{SH}_1 = \mathbb{V}(\mathbf{C}^{dp})$,  
\item[d]   $\mathbb{RDPCS}\rm t\mathbb{H}_1 =  \mathbb{V}(\mathbf{L}^{dp}_1)$.       
\end{thlist}
\end{Theorem}

It is clear from Theorem \ref{Theo_050719_01} that the logics $\mathcal{RDM}\mathcal{SH}_1$ and $\mathcal{RDMS}{\rm t}\mathcal{SH}_1$ are equivalent and so are  $\mathcal{RDMH}_1$ and $\mathcal{RDMS}{\rm t}\mathcal{H}_1$.

The following corollary is immediate from Theorem \ref{Theo_050719_01}.

\begin{Corollary}
The logics $\mathcal{RDM}\mathcal{SH}_1$ and $\mathcal{RDPCS}\rm t\mathcal{SH}_1$
are decidable.
\end{Corollary}  

Let $\mathbb{RDQD}{\rm cm}\mathbb{S}{\rm t}\mathbb{SH}_1$ be the subvariety of 
$\mathbb{RDQD}\mathbb{S}{\rm t}\mathbb{SH}_1$ defined by the commutative law: 

\qquad $x \to y \approx y \to x$.

\begin{Corollary} \label{coro_150719_01}
The logic $\mathcal{RDQD} cm\mathcal{S}{\rm t}\mathcal{H}_1$ is defined, as an extension of $\mathcal{RDQDS}{\rm t}\mathcal{SH}_1$, by 
\begin{itemize}
\item[]  $(\alpha \to \beta) \to_H \ (\beta \to \alpha)$.
\end{itemize} 
\end{Corollary}

 The following theorem is an immediate consequence of Theorem \ref{TH_Main} and Theorem \ref{Theo_050719_01}.
\begin{Theorem} {\rm \cite[Corollary 3.5]{sankappanavar2014b}} \label{Theorem_180719_02} 
\begin{thlist}
\item[a]
$\mathbb{RDQD}{\rm cm}\mathbb{S}{\rm t}\mathbb{SH}_1$ =$ \mathbb{V}(\mathbf{L}^{dm}_{10}) \lor \mathbb{V}(\mathbf{L}^{dp}_{10}) \lor \mathbb{V}(\mathbf{D}_1)$,  
\item[b] $\mathbb{RDM}{\rm cm}\mathbb{SH}_1 = \mathbb{RDM}{\rm cm}\mathbb{S}{\rm t}\mathbb{SH}_1 =  \mathbb{V}(\{\mathbf{L}^{dm}_{10},  \mathbf{D}_1\})$,    
\item[c]  $\mathbb{RDPC}{\rm cm}\mathbb{S}{\rm t}\mathbb{SH}_1 =  \mathbb{V}((\mathbf{L}^{dp}_{10})$,
\item[d]  $\mathbb{RDM}{\rm cm}\mathbb{SH}_1 \cap \mathbb{RDPC}{\rm cm}\mathbb{S}{\rm t}\mathbb{SH}_1 =  \mathbb{V}(\bar{2}^e)$.
\end{thlist} 
\end{Theorem}

It follows from the preceding theorem that the logics $\mathcal{RDQD}\rm cm\mathcal{S}{\rm t}\mathcal{SH}_1$, $\mathcal{RDM}{\rm cm}\mathcal{SH}_1$  and $\mathcal{RDPC}{\rm cm}\mathcal{S}{\rm t}\mathcal{SH}_1$ are decidable.  \\

{\bf In the rest of this section, unless otherwise stated, the phrase  ``defined, modulo 
$\mathbb{RDQDS}{\rm t}\mathbb{SH}_1$, by''  is abbreviated to the phrase ``defined by''
  in the context of varieties.  
Similarly,  the phrase 
``defined, as an extension of the logic 
$\mathcal{RDQDS}{\rm t}\mathbb{SH}_1$, by'' is also abbreviated to the phrase ``defined by'' in the case of logics.}\\

The theorems that appear in the rest of this section were proved in \cite{sankappanavar2014b}.
Each of the corollaries appearing below follows from the theorem immediately preceding it and Theorem \ref{teo_040417_01}.

\begin{Theorem}
The variety $\mathbf{V(\{ L^{dm}_1, L^{dp}_1, L^{dm}_3, L^{dp}_3, D_2\})}$ is defined 
by the identity:
\begin{itemize}
\item[ ] $(x \to y) \to (0 \to y) \approx (x \to y) \to 1$.
\end{itemize}
\end{Theorem}

\begin{Corollary}
The logic $\mathcal{L}(\mathbb{V}\mathbf{(L_1^{dm}, L_1^{dp}, L_3^{dm}, L_3^{dp}, D_2}))$ is defined 
by
\begin{itemize}
\item[] $[(\alpha \to \beta) \to (\bot \to \beta)] \Leftrightarrow_H [(\alpha \to \beta) \to \top]$.
\end{itemize} 
\end{Corollary}

The variety generated by $\mathbf{D}_1$ was axiomatized earlier in Section \ref{S7}.  
Here are two more bases for it.

\begin{Theorem}
 $\mathbb{V}(\mathbf{D}_1)$ is defined 
by 
\begin{itemize}
\item[ ] $x \to (y \to z) \approx z \to (x \to y)$.
\end{itemize}      
It is also defined 
by
\begin{itemize}
\item[ ] $(x \to y) \to (u \to w) \approx (x \to u) \to (y \to w)$   \quad    \rm{(Medial Law)}.
\end{itemize}     
\end{Theorem}

\begin{Corollary}
The logic $\mathcal{L}(\mathbb{V}\mathbf{(D_1}))$ is defined 
by
\begin{itemize}
\item[] $\alpha \to (\beta \to \gamma) \Leftrightarrow_H \gamma \to (\alpha \to \beta)$.
\end{itemize}      
It is also defined by
\begin{itemize}
\item[] $(\alpha \to \beta) \to (\gamma \to \delta) \Leftrightarrow_H ((\alpha \to \gamma) \to (\beta \to \delta))$,   \quad    \rm{(Medial Law)}.
\end{itemize}     
\end{Corollary}

\begin{Theorem}
The variety $\mathbb{V}(\mathbf{\{ L^{dm}_1, L^{dp}_1, L^{dm}_2, L^{dp}_2, D_2\})}$ is defined
by: 
\begin{itemize}
\item[ ] $y \leq x \to y$. 
\end{itemize}  
It is also defined  
by: 
\begin{itemize}
\item[ ] $[(x \to y) \to y] \to (x \to y) \approx x \to y$.  
\end{itemize}  
It is also defined 
by
\begin{itemize}
\item[ ] $x \to (y \to z) \approx (x \to y) \to (x \to z)$   \quad     \rm{(Left distributive law)}.  
\end{itemize}  
\end{Theorem}

\begin{Corollary}
The logic $\mathcal{L}(\mathbb{V}(\mathbf{\{ L^{dm}_1, L^{dp}_1, L^{dm}_2, L^{dp}_2, D_2\})}$ is defined 
by
\begin{itemize}
\item[ ] $\beta \land (\alpha  \to \beta)$. 
\end{itemize}  
It is also defined by: 
\begin{itemize}
\item[] $[(\alpha \to \beta) \to \beta] \to (\alpha \to \beta) \Leftrightarrow_H (\alpha  \to \beta)$.  
\end{itemize}  
It is also defined by
\begin{itemize}
\item[] $\alpha \to (\beta \to \gamma) \Leftrightarrow_H [(\alpha \to \beta) \to (\alpha \to \gamma)$.    
\end{itemize}
\end{Corollary}

\begin{Theorem}
The variety  
$\mathbb{V}(\mathbf{\{ L^{dm}_1, L^{dp}_1, L^{dm}_2, L^{dp}_2, L^{dm}_5, L^{dp}_5, L^{dm}_6, L^{dp}_6, D_2\}})$ is defined 
by: 
\begin{itemize}
\item[ ] $[x \to (y \to x)] \to x \approx x$.
\end{itemize}  
\end{Theorem}

\begin{Corollary}
The logic $\mathcal{L}(\mathbb{V}(\mathbf{\{ L^{dm}_1, L^{dp}_1, L^{dm}_2, L^{dp}_2, L^{dm}_5, L^{dp}_5, L^{dm}_6, L^{dp}_6, D_2\})})$ is defined 
by
\begin{itemize}
\item[] $ [\{\alpha \to (\beta \to \alpha)\} \to \alpha] \Leftrightarrow_H \alpha$.
\end{itemize}  
\end{Corollary}

\begin{Theorem}
$\mathbb{V}(\mathbf{\{L^{dp}_1, L^{dp}_2, L^{dp}_5, L^{dp}_6\}) }$ is defined
by: 
\begin{itemize}
\item[{\rm(1)}] $[x \to (y \to x)] \to x \approx x$,  
\item[{\rm(2)}] $x \lor x' \approx 1$.
\end{itemize}  
\end{Theorem}

\begin{Corollary}
The logic $\mathcal{L}(\mathbb{V}(\mathbf{\{ L^{dp}_1, L^{dp}_2, L^{dp}_5, L^{dp}_6\})})$ is defined 
by
\begin{itemize}
\item[{\rm(1)}] $[\{\alpha \to (\beta \to \alpha)\} \to \alpha] \Leftrightarrow_H \alpha$, 
\item[{\rm(2)}]  $\alpha \ \lor \ \sim\alpha$.
\end{itemize}  
\end{Corollary}

 Recall that $x^+ := x'{^*}'$.

\begin{Theorem}
The variety  
$\mathbb{V}(\mathbf{\{ L^{dm}_1,
  L^{dp}_1, L^{dm}_2, L^{dp}_2, L^{dm}_3, L^{dp}_3, 
L^{dm}_4, L^{dp}_4, L^{dm}_5, L^{dm}_6,} $\\
$\mathbf{L^{dm}_7, L^{dm}_8, D_2, D_3\})}$ is defined 
by the identity:
\begin{itemize}
\item[ ]  $(0 \to 1)^+ \to (0 \to 1)'  \approx 0 \to 1$.  
 \end{itemize}  
\end{Theorem}

\begin{Corollary}
The logic $\mathcal{L}((\mathbb{V}(\mathbf{\{ L^{dm}_1,
  L^{dp}_1, L^{dm}_2, L^{dp}_2, L^{dm}_3, L^{dp}_3,  L^{dm}_4, L^{dp}_4, L^{dm}_5, L^{dm}_6,}$\\
$\mathbf{  L^{dm}_7, L^{dm}_8, D_2, D_3\}}))$ is defined 
by
\begin{itemize}
\item[]  $[(\bot \to \top)^+ \to \ \sim(\bot \to \top)]  \Leftrightarrow_H (\bot \to \top)$.
 \end{itemize} 
\end{Corollary}

\begin{Theorem}
The variety  
$\mathbf{ 
 V(\{ L^{dp}_1,  L^{dp}_2,  L^{dp}_3, 
 L^{dp}_4\})}$ 
 is defined
by the identities:
\begin{itemize}
\item[{\rm(1)}]  $(0 \to 1)^+ \to (0 \to 1)'  \approx 0 \to 1$,  
\item[{\rm(2)}]  $x \lor x' \approx 1$.
 \end{itemize}  
\end{Theorem}

\begin{Corollary}
The logic $\mathcal{L}(\mathbb{V}(\mathbf{\{ L^{dp}_1,
  L^{dp}_2, L^{dp}_3, L^{dp}_4\}}))$ is defined 
  by
\begin{itemize}
\item[{\rm(1)}]  $[(\bot \to \top)^+ \to \ \sim(\bot \to \top)] \Leftrightarrow_H (\bot \to \top)$,
\item[{\rm(2)}]  $\alpha \ \lor \sim\alpha$.
 \end{itemize}
\end{Corollary}


\begin{Theorem}
The variety \\
$\mathbf{V(\{ L^{dm}_1,
   L^{dm}_2,  L^{dm}_3,  
L^{dm}_4,  L^{dm}_5, L^{dm}_6, L^{dm}_7, L^{dm}_8},$ 
$\mathbf{ D_2, D_3\})}$ is defined  
by the identities:
\begin{itemize}
\item[{\rm(1)}]  $(0 \to 1)^+ \to (0 \to 1)'  \approx 0 \to 1$, 
\item[{\rm(2)}] $x'' \approx x$.
 \end{itemize}  
\end{Theorem}

\begin{Corollary}
The logic $\mathcal{L}(\mathbb{V}(\mathbf{\{ L^{dm}_1,
  L^{dm}_2, L^{dm}_3, L^{dm}_4, L^{dm}_5, L^{dm}_6, L^{dm}_7, L^{dm}_8, 
 D_2, D_3\}}))$ is defined 
 by
\begin{itemize}
\item[{\rm(1)}]  $[(\bot \to \top)^+ \to \ \sim(\bot \to \top)]  \Leftrightarrow_H (\bot \to \top)$, 
\item[{\rm(2)}] $\alpha \to_H \ \sim\sim\alpha$.
 \end{itemize}  
\end{Corollary}

\begin{Theorem}
The variety 
$\mathbb{V}({\mathbf{\{L^{dm}_5, L^{dm}_6, L^{dm}_7, L^{dm}_8, D_3\}})}$ is defined     
by the identity:
\begin{itemize}
\item[]  $(0 \to 1)^+ \to (0 \to 1)  \approx (0 \to 1)'$. 
 \end{itemize}  
\end{Theorem}

\begin{Corollary}
The logic $\mathcal{L}(\mathbb{V}(\mathbf{\{ L^{dm}_5,
  L^{dm}_6, L^{dm}_7, L^{dm}_8, D_3\}}))$ is defined 
 by
\begin{itemize}
\item[]  $[(\bot \to \top)^+ \to (\bot \to \top)]  \Leftrightarrow_H \ \sim(\bot \to \top)$.
 \end{itemize}
\end{Corollary}

$\mathbb{V}\mathbf{(D_3)}$ was axiomatized in Section \ref{S7}. 
Here is another base for it. 

\begin{Theorem}
$\mathbb{V}(\mathbf{D_3)}$ is defined
by the identities:
\begin{itemize}
\item[{\rm(1)}]  $(0 \to 1)^+ \to (0 \to 1) \approx (0 \to 1)'$,
\item[{\rm(2)}]   $x \lor x^* \approx 1$.  
 \end{itemize}  
\end{Theorem}

\begin{Corollary}
The logic $\mathcal{L}(\mathbb{V}(\mathbf{\{D_3\}}))$ is defined 
by
\begin{itemize}
\item[{\rm(1)}]  $(\bot \to \top)^+ \to (\bot \to \top) \Leftrightarrow_H \ \sim(\bot \to \top)$,
\item[{\rm(2)}]   $\alpha \lor \alpha^*$.
\end{itemize}
\end{Corollary}

\begin{Theorem}
The variety generated by the algebras
$\mathbf{ 
  L^{dm}_1, L^{dm}_2,}$
  $\mathbf{ L^{dm}_3,  L^{dm}_4,}$ 
$\mathbf{D_2, D_3}$ is defined
by the identities:
\begin{itemize}
\item[{\rm(1)}]  $(0 \to 1)^+ \to (0 \to 1)' \approx 0 \to 1$,
\item[{\rm(2)}]  $(0 \to 1)^+ \to (0 \to 1){^*}'{^*} \approx 0 \to 1$,
\item[{\rm(3)}]   $x'' \approx x$.  
 \end{itemize}  
\end{Theorem}

\begin{Corollary}
The logic $\mathcal{L}(\mathbb{V}(\mathbf{\{ L^{dm}_1,
  L^{dm}_2, L^{dm}_3, L^{dm}_4, D_2, D_3\})})$ is defined  
  by
\begin{itemize}
\item[{\rm(1)}]  $(\bot \to \top)^+ \to \ \sim(\bot \to \top) \Leftrightarrow_H (\bot \to \top)$,
\item[{\rm(2)}]  $[(\bot \to \top)^+ \to \ (\sim(\bot \to \top){^*}){^*}] \Leftrightarrow_H (\bot \to \top)$,
\item[{\rm(3)}]   $\alpha \to_H \ \sim\sim\alpha$.  
 \end{itemize}  
\end{Corollary}

\begin{Theorem}
The variety generated by the algebras\\ 
$\mathbf{ L^{dm}_5,
  L^{dp}_5, L^{dm}_6, L^{dp}_6, L^{dm}_7, L^{dp}_7, 
L^{dm}_8, L^{dp}_8, L^{dp}_9, L^{dm}_9, L^{dp}_{10}, L^{dm}_{10}},$ \\
$\mathbf{ D_1, D_3}$ is defined
by the identity:
\begin{itemize}
\item[ ]  $(0 \to 1)^+ \to (0 \to 1)' \approx (0 \to 1)'$.
 \end{itemize}  
\end{Theorem}

\begin{Corollary}
The logic \\
$\mathcal{L}(\mathbb{V}(\mathbf{\{ L^{dm}_5,
  L^{dp}_5, L^{dm}_6, L^{dp}_6, L^{dm}_7, L^{dp}_7, L^{dm}_8, L^{dp}_8,\\
   L^{dm}_9, L^{dp}_9, L^{dm}_{10}, L^{dp}_{10}, D_1, D_3\}}))$ is defined  
   by
\begin{itemize}
\item[]  $[(\bot \to \top)^+ \to \ \sim(\bot \to \top)]  \Leftrightarrow_H \ \sim(\bot \to \top)$.
 \end{itemize} 
\end{Corollary}

\begin{Theorem}
The variety generated by the algebras\\ 
$\mathbf{ 
  L^{dp}_5,  L^{dp}_6,  L^{dp}_7, 
 L^{dp}_8, L^{dp}_9, L^{dp}_{10}} $ 
is defined 
by the identities:
\begin{itemize}
\item[{\rm(1)}]  $(0 \to 1)^+ \to (0 \to 1)' \approx (0 \to 1)'$,
\item[{\rm(2)}]   $x \lor x' \approx 1$.
 \end{itemize}  
\end{Theorem}

\begin{Corollary}
The logic $\mathcal{L}(\mathbb{V}(\mathbf{\{L^{dp}_5, L^{dp}_6, L^{dp}_7, L^{dp}_8,
  L^{dp}_9, L^{dp}_{10}\}}))$ is defined \\
by
\begin{itemize}
\item[{\rm(1)}]  $[(\bot \to \top)^+ \to \ \sim(\bot \to \top)]  \Leftrightarrow_H \ \sim(\bot \to \top)$,
\item[{\rm(2)}]    $\alpha \ \lor \ \sim\alpha$.
 \end{itemize}  
\end{Corollary}


\begin{Theorem}
The variety generated by the algebras\\ 
$\mathbf{ L^{dm}_5,
   L^{dm}_6,  L^{dm}_7,  
L^{dm}_8, L^{dm}_9, L^{dm}_{10}}$,  
$\mathbf{ D_1, D_3}$ is defined 
by the identities:
\begin{itemize}
\item[{\rm(1)}]  $(0 \to 1)^+ \to (0 \to 1)' \approx (0 \to 1)'$,
\item[{\rm(2)}]   $x'' \approx x$.
 \end{itemize}  
\end{Theorem}

\begin{Corollary}
The logic $\mathcal{L}(\mathbb{V}(\mathbf{\{L^{dm}_5, L^{dm}_6, L^{dm}_7, L^{dm}_8,
  L^{dm}_9, L^{dm}_{10}, D_1, D_3\}}))$ is defined \\
 by
  \begin{itemize}
  \item[{\rm(1)}]  $[(\bot \to \top)^+ \to \ \sim(\bot \to \top)]  \Leftrightarrow_H \ \sim(\bot \to \top)$,
\item[{\rm(2)}]   $\alpha \to_H \ \sim\sim\alpha$.  
 \end{itemize}
\end{Corollary}

\begin{Theorem}
The variety generated by the algebras\\ 
$\mathbf{ D_1, D_3}$ is defined 
by the identities:
\begin{itemize}
\item[{\rm(1)}]  $(0 \to 1)^+ \to (0 \to 1)' \approx (0 \to 1)'$,
\item[{\rm(2)}]   $x \lor x^* \approx 1$.
 \end{itemize}  
\end{Theorem}

\begin{Corollary}
The logic $\mathcal{L}(\mathbb{V}(\mathbf{\{D_1, D_3\}}))$ is defined \\
by
\begin{itemize}
 \item[\rm{(1)}]  $[(\bot \to \top)^+ \to \ \sim(\bot \to \top)]  \Leftrightarrow_H \ \sim(\bot \to \top)$,
\item[{\rm(2)}]    $\alpha \lor \alpha^*$.
 \end{itemize} 
\end{Corollary}


\begin{Theorem}
The variety generated by the algebras\\ 
$\mathbf{ L^{dm}_5,
   L^{dm}_6,  L^{dm}_7,  
L^{dm}_8}$, 
$\mathbf{ D_3}$ is defined
by the identities:
\begin{itemize}
\item[{\rm(1)}]  $(0 \to 1)^+ \to (0 \to 1)' \approx (0 \to 1)'$,
\item[{\rm(2)}]    $(0 \to 1)^+ \to (0 \to 1)' \approx (0 \to 1)$.
 \end{itemize} 
It is also defined
by 
\begin{itemize}
\item[ ]     $(0 \to 1)' \approx 0 \to 1$. 
 \end{itemize} 
\end{Theorem}

\begin{Corollary}
The logic $\mathcal{L}(\mathbb{V}(\mathbf{\{L^{dm}_5, L^{dm}_6, L^{dm}_7, L^{dm}_8,
  D_3\}}))$ is defined \\
 by
\begin{itemize}
 \item[{\rm(1)}]  $[(\bot \to \top)^+ \to \ \sim(\bot \to \top)]  \Leftrightarrow_H \ \sim(\bot \to \top)$,
 \item[{\rm(2)}]  $[(\bot \to \top)^+ \to \ \sim(\bot \to \top)]  \Leftrightarrow_H (\bot \to \top)$.
 \end{itemize} 
 It is also defined by 
\begin{itemize}
\item[ ]     $\sim(\bot \to \top) \Leftrightarrow_H (\bot \to \top)$. 
 \end{itemize} 
\end{Corollary}

\begin{Theorem}
The variety generated by the algebras\\ 
$\mathbf{ L^{dm}_1,
  L^{dp}_1, L^{dm}_2, L^{dp}_2, L^{dm}_3, L^{dp}_3, 
L^{dm}_4, L^{dp}_4, L^{dp}_5,\\
 L^{dp}_6, L^{dp}_7, L^{dp}_8}$,  
$\mathbf{L^{dm}_9, L^{dp}_9, L^{dm}_{10}, L^{dp}_{10},  D_1, D_2}$ is defined
by the identity:
\begin{itemize}
\item[ ]  $(0 \to 1)' \to (0 \to 1)  \approx 0 \to 1$.  
 \end{itemize}  
\end{Theorem}

\begin{Corollary}
The logic corresponding to the variety generated by the algebras\\
$\mathbf{L^{dm}_1,
  L^{dp}_1, L^{dm}_2, L^{dp}_2, L^{dm}_3, L^{dp}_3, 
L^{dm}_4, L^{dp}_4, L^{dp}_5, L^{dp}_6, L^{dp}_7, L^{dp}_8,  
L^{dm}_9, L^{dp}_9, L^{dm}_{10}, L^{dp}_{10}, D_1,  D_2}$ is defined
by
\begin{itemize}
\item[]   $[\sim(\bot \to \top) \to (\bot \to \top)]  \Leftrightarrow_H (\bot \to \top)$.
 \end{itemize} 
\end{Corollary}

\begin{Theorem}
The variety generated by the algebras\\ 
$\mathbf{ L^{dm}_1,
   L^{dm}_2,  L^{dm}_3,  
L^{dm}_4, L^{dm}_9,  L^{dm}_{10}},$ 
$\mathbf{ D_1, D_2}$ 
is defined
by the identities:
\begin{itemize}
\item[{\rm(1)}]  $(0 \to 1)' \to (0 \to 1)  \approx 0 \to 1$,  
\item[{\rm(2)}]  $x'' \approx x$.
 \end{itemize}  
\end{Theorem}

\begin{Corollary}
The logic $\mathcal{L}(\mathbb{V}(\mathbf{\{L^{dm}_1,
   L^{dm}_2, L^{dm}_3,  
L^{dm}_4,   
L^{dm}_9, L^{dm}_{10}, D_1, D_2\}}))$ is defined 
 by
\begin{itemize}
\item[{\rm(1)}]   $[\sim(\bot \to \top) \to (\bot \to \top)]  \Leftrightarrow_H (\bot \to \top)$,
\item[{\rm(2)}]   $\alpha \to_H \ \sim\sim\alpha$.  
 \end{itemize}  
\end{Corollary}

\begin{Theorem}
The variety generated by the algebras\\ 
$\mathbf{ D_1, D_2}$ is defined
by the identities:
\begin{itemize}
\item[{\rm(1)}]  $(0 \to 1)' \to (0 \to 1)  \approx 0 \to 1,$
\item[{\rm(2)}]   $x \lor x^* \approx 1$.
 \end{itemize}  
\end{Theorem}

\begin{Corollary}
The logic $\mathcal{L}(\mathbb{V}(\mathbf{\{D_1, D_2\}}))$ is defined 
by
\begin{itemize}
\item[{\rm(1)}]   $[\sim(\bot \to \top) \to (\bot \to \top)]  \Leftrightarrow_H (\bot \to \top)$,
\item[{\rm(2)}]   $\alpha \ \lor \ \alpha^*$.
 \end{itemize}
\end{Corollary}


\begin{Theorem}
The variety generated by the algebras\\
 $\mathbf{L_1^{dm}, L_1^{dp}, L_3^{dm}, L_3^{dp}, L_6^{dm}, L_6^{dp}, L_8^{dm}, L_8^{dp},
D_1, D_2, D_3}$ is defined by the identity:
 \begin{itemize}
\item[ ] $x \lor [y \to (x \lor y)] \approx (0 \to x) \lor x \lor y$.   
 \end{itemize}  
\end{Theorem}

\begin{Corollary}
The logic $\mathcal{L}(\mathbb{V}(\mathbf{\{L_1^{dm}, L_1^{dp}, L_3^{dm}, L_3^{dp},\\
 L_6^{dm}, L_6^{dp}, L_8^{dm}, L_8^{dp},
D_1, D_2, D_3\}}))$ is defined
by
 \begin{itemize}
\item[] $[\alpha \lor \{\beta \to (\alpha \lor \beta)\}] \Leftrightarrow_H [(\bot \to \alpha) \lor \alpha \lor \beta]$.   
 \end{itemize}  
\end{Corollary}

\begin{Theorem}
The variety generated by the algebras\\
 $\mathbf{ L^{dm}_2, L^{dp}_2, L^{dm}_5, L^{dp}_5, D_2}$ is defined by the identity:
 \begin{itemize}
\item[ ] $x \lor (y \to x) \approx [(x \to y) \to y] \to x$.  
\end{itemize}  
\end{Theorem}

\begin{Corollary}
The logic $\mathcal{L}(\mathbb{V}(\mathbf{\{L_2^{dm}, L_2^{dp}, L_5^{dm}, L_5^{dp},
  D_2\}}))$ is defined
by
 \begin{itemize}
\item[] $[\alpha \lor (\beta \to \alpha)] \Leftrightarrow_H [\{(\alpha \to \beta) \to \beta\} \to \alpha]$.
\end{itemize}
\end{Corollary}

\begin{Theorem}
The variety generated by the algebras\\
 $\mathbf{ L^{dm}_3, L^{dp}_3, L^{dm}_4, L^{dp}_4,  D_1, D_2, D_3 }$ is defined
by the identity:
 \begin{itemize}
\item[ ] $x \lor (x \to y) \approx x \to [x \lor (y \to 1)]$. 
\end{itemize}  
\end{Theorem}

\begin{Corollary}
The logic $\mathcal{L}(\mathbb{V}(\mathbf{\{L_3^{dm}, L_3^{dp}, L_4^{dm}, L_4^{dp},
D_1, D_2, D_3\}}))$ is defined
by
 \begin{itemize}
\item[] $[\alpha \lor (\alpha \to \beta)] \Leftrightarrow_H [\alpha \to \{\alpha \lor (\beta \to \top)\}]$.
\end{itemize}  
\end{Corollary}

\begin{Theorem}
The variety generated by the algebras\\
 $\mathbf{ L^{dm}_5, L^{dp}_6, L^{dm}_7, L^{dp}_8,  
 D_3}$ is defined by the identity:
 \begin{itemize}
\item[ ]  $(0 \to 1)^* \to (0 \to 1) \approx (0 \to 1)'$.   
\end{itemize}  
\end{Theorem}

\begin{Corollary}
The logic $\mathcal{L}(\mathbb{V}(\mathbf{\{L_5^{dm}, L_6^{dp}, L_7^{dm}, L_8^{dp},
 D_3\}}))$ is defined
by
 \begin{itemize}
\item[] $[(\bot \to \top)^* \to (\bot \to \top)] \Leftrightarrow_H \ \sim(\bot \to \top)$.
\end{itemize}  
\end{Corollary}

\begin{Theorem}
The variety generated by the algebras\\
 $\mathbf{ L^{dm}_1, L^{dp}_1, L^{dm}_2, L^{dp}_2, L^{dm}_3, L^{dp}_3}$, 
$\mathbf{L^{dm}_4, L^{dp}_4, D_2}$ is defined 
by the
 identity:
\begin{itemize}
\item[ ] $0 \to 1 \approx 1$   \rm{(FTT identity)}.
\end{itemize}  
\end{Theorem}

\begin{Corollary}
The logic $\mathcal{L}(\mathbb{V}(\mathbf{\{L_1^{dm}, L_1^{dp}, L_2^{dm}, L_2^{dp},
 L_3^{dm}, L_3^{dp}, L_4^{dm}, L_4^{dp},
 D_2\}}))$ is defined
by
\begin{itemize}
\item[]  $\bot \to \top $   \rm{(FTT identity)}.
\end{itemize}  
\end{Corollary}

\begin{Theorem}
The variety generated by the algebras\\
 $\mathbf{ L^{dm}_1, L^{dp}_1, L^{dm}_3, L^{dp}_3, L^{dm}_6, L^{dp}_6}$, 
$\mathbf{L^{dm}_8, L^{dp}_8, D_1, D_2, D_3}$ is \\
defined
by the identity:
\begin{itemize}
\item[ ] $x \lor (y \to x) \approx (x \lor y) \to x$.
\end{itemize}  
\end{Theorem}

\begin{Corollary}
The logic $\mathcal{L}(\mathbb{V}(\mathbf{\{L_1^{dm}, L_1^{dp}, L_3^{dm}, L_3^{dp},
 L_6^{dm}, L_6^{dp}, L_8^{dm}, L_8^{dp},
D_1, D_2, D_3\}})$ \\ 
is defined 
by
\begin{itemize}
\item[] $[\alpha \lor (\beta \to \alpha)] \Leftrightarrow_H  [(\alpha \lor \beta) \to \alpha]$.
\end{itemize}  
\end{Corollary}

\begin{Theorem}
The variety generated by the algebras\\
 $\mathbf{ L^{dp}_1,  L^{dp}_3,  L^{dp}_6}$, 
$\mathbf{L^{dp}_8}$ is 
defined 
by the identities:
\begin{itemize}
\item[{\rm(1)}] $x \lor (y \to x) \approx (x \lor y) \to x$,
\item[{\rm(2)}]  $x \lor x' \approx 1$.
\end{itemize}  
\end{Theorem}

\begin{Corollary}
The logic $\mathcal{L}( \mathbb{V}(\mathbf{\{L^{dp}_1,  L^{dp}_3,  L^{dp}_6, L^{dp}_8\}}))$ is defined
by
\begin{itemize}
\item[{\rm(1)}] $[\alpha \lor (\beta \to \alpha)] \Leftrightarrow_H  [(\alpha \lor \beta) \to \alpha]$,
\item[{\rm(2)}]   $\alpha \ \lor \ \sim\alpha$.
\end{itemize}  
\end{Corollary}


\begin{Theorem}
The variety generated by the algebras\\
 $\mathbf{ L^{dm}_1, L^{dm}_3,  L^{dm}_6}$, 
$\mathbf{L^{dm}_8,  D_1, D_2, D_3}$ is 
defined 
by the identities:
\begin{itemize}
\item[{\rm(1)}] $x \lor (y \to x) \approx (x \lor y) \to x$,
\item[{\rm(2)}] $x'' \approx x$.
\end{itemize}  
\end{Theorem}

\begin{Corollary}
The logic $\mathcal{L}(\mathbb{V}(\mathbf{\{L_1^{dm}, L_3^{dm}, L_6^{dm},
 L_8^{dm},  
D_1, D_2,
 D_2\}}))$ is defined
by
\begin{itemize}
\item[{\rm(1)}] $[\alpha \lor (\beta \to \alpha)] \Leftrightarrow_H  [(\alpha \lor \beta) \to \alpha]$,
\item[{\rm(2)}]   $\alpha \to_H \ \sim\sim\alpha$. 
\end{itemize}
\end{Corollary}

\begin{Theorem}
The variety generated by the algebras\\
 $\mathbf{ L^{dm}_1, L^{dp}_1, L^{dm}_2, L^{dp}_2, L^{dm}_5, L^{dp}_5, L^{dm}_6, L^{dp}_6, L^{dm}_9, L^{dp}_9}$, \\
 $\mathbf{D_1, D_2, D_3}$ 
is defined 
by the identity:
 \begin{itemize}
\item[ ] $x^* \lor (x \to y) \approx (x \lor y) \to y$.   
\end{itemize}  
\end{Theorem}

\begin{Corollary}
The logic \\
$\mathcal{L}(\mathbb{V}(\mathbf{\{L_1^{dm}, L^{dp}_1, L_2^{dm}, L^{dp}_2, L^{dm}_5, L^{dp}_5,\\ L_6^{dm}, L^{dp}_6,
 L_9^{dm}, L_9^{dp},  
D_1, D_2,
 D_2\}}))$ is defined
by
 \begin{itemize}
\item[] $[\alpha^* \lor (\alpha \to \beta)] \Leftrightarrow_H [(\alpha \lor \beta) \to \beta]$.  
\end{itemize}  
\end{Corollary}

\begin{Theorem}
 $\mathbf{V(\{ L^{dp}_1,  L^{dp}_2,  L^{dp}_5,  L^{dp}_6,  L^{dp}_9\})}$   
is defined 
by the identity:
 \begin{itemize}
\item[{\rm(1)}] $x^* \lor (x \to y) \approx (x \lor y) \to y$,   
\item[{\rm(2)}] $x \lor x' \approx 1$.
\end{itemize}  
\end{Theorem}

\begin{Corollary}
The logic $\mathcal{L}(\mathbb{V}(\mathbf{ \{L^{dp}_1, L^{dp}_2, L^{dp}_5, L^{dp}_6,
  L_9^{dp}\}}))$ is defined
 by
 \begin{itemize}
 \item[{\rm(1)}] $[\alpha^* \lor (\alpha \to \beta)] \Leftrightarrow_H [(\alpha \lor \beta) \to \beta]$,   
\item[{\rm(2)}]    $\alpha \ \lor \ \sim\alpha$.
 \end{itemize} 
\end{Corollary}


\begin{Theorem}
The variety generated by the algebras\\
 $\mathbf{ L^{dm}_1,  L^{dm}_2,  L^{dm}_5,  L^{dm}_6,  L^{dm}_9}$, 
 $\mathbf{D_1, D_2, D_3}$ 
is defined 
by the identity:
 \begin{itemize}
\item[{\rm(1)}] $x^* \lor (x \to y) \approx (x \lor y) \to y$,   
\item[{\rm(2)}]  $x'' \approx x$.
\end{itemize}  
\end{Theorem}

\begin{Corollary}
The logic $\mathcal{L}(\mathbb{V}(\mathbf{\{L^{dm}_1, L^{dm}_2, L^{dm}_5, L^{dm}_6,
  L_9^{dm}, D_1, D_2, D_3\}}))$ is defined
  by
 \begin{itemize}
 \item[{\rm(1)}] $[\alpha^* \lor (\alpha \to \beta)] \Leftrightarrow_H [(\alpha \lor \beta) \to \beta]$,   
\item[{\rm(2)}]   $\alpha \to_H \ \sim\sim\alpha$.
 \end{itemize} 
\end{Corollary}

\begin{Theorem}
The variety generated by the algebras\\
 $\mathbf{L^{dm}_5, L^{dp}_5, D_2}$  
 is defined 
by the identity:
 \begin{itemize}
\item[ ] $x \lor (0 \to x) \lor (y \to 1) \approx x \lor [(x \to 1) \to (x \to y)]$.
\end{itemize}  
\end{Theorem}

\begin{Corollary}
The logic $\mathcal{L}(\mathbb{V}(\mathbf{\{L^{dm}_5, L^{dp}_5. 
  D_2\}}))$ is defined
by
 \begin{itemize}
\item[] $[\alpha \lor (\bot \to \alpha) \lor (\beta \to \top)] \Leftrightarrow_H \alpha \lor [(\alpha \to \top) \to (\alpha \to \beta)]$.
\end{itemize}
\end{Corollary}

\begin{Theorem}
The variety generated by the algebras\\
 $\mathbf{ L^{dm}_6, L^{dp}_6, D_2}$  
 defined 
by the identity:
\begin{itemize}
\item[ ] $x \lor y \lor (x \to y) \approx x \lor [(x \to y) \to 1]$.
\end{itemize}  
\end{Theorem}

\begin{Corollary}
The logic $\mathcal{L}(\mathbb{V}(\mathbf{\{L^{dm}_6, L^{dp}_6, 
  D_2\}}))$ is defined
by
\begin{itemize}
\item[] $\alpha \lor \beta \lor (\alpha \to \beta) \Leftrightarrow_H \alpha \lor [(\alpha \to \beta) \to \top]$.
\end{itemize}
\end{Corollary}

\begin{Theorem}
The variety generated by the algebras\\
 $\mathbf{ L^{dm}_1, L^{dp}_1, L^{dm}_7, {L^{dp}_7, D_2}}$  
is defined
by the identity:
 \begin{itemize}
\item[ ] $x \lor [(0 \to y) \to y] \approx x \lor [(x \to 1) \to y]$.
\end{itemize}  
\end{Theorem}

\begin{Corollary}
The logic $\mathcal{L}(\mathbb{V}(\mathbf{\{L^{dm}_1, L^{dp}_1, L^{dm}_7, L^{dp}_7,
  D_2\}}))$ is defined
 by
 \begin{itemize}
\item[] $[\alpha \lor \{(\bot \to \beta) \to \beta\}] \Leftrightarrow_H  [\alpha \lor [(\alpha \to \top) \to \beta]$.
\end{itemize}  
\end{Corollary}

\begin{Theorem}
The variety generated by the algebras\\
 $\mathbf{ L^{dm}_7, L^{dp}_7, L^{dm}_8, L^{dp}_8, D_1, D_2, D_3}$  
is defined 
by the identity:
 \begin{itemize}
\item[ ] $x \ \lor [x \to (y \land (0 \to y))] \approx x \to [(x \to y) \to y]$.  
\end{itemize}  
\end{Theorem}

\begin{Corollary}
The logic $\mathcal{L}(\mathbb{V}(\mathbf{\{L^{dm}_7, L^{dp}_7, L^{dm}_8, L^{dp}_8, D_1, D_2,
  D_3\}}))$ is defined
by
 
 \begin{itemize}
 \item[] $[\alpha \lor [\alpha \to \{\beta \land (\bot \to \beta) \}]] \Leftrightarrow_H  [\alpha \to [(\alpha \to \beta) \to \beta]$.
\end{itemize} 
 
\end{Corollary}

\begin{Theorem}
The variety generated by the algebras\\
 $\mathbf{ L^{dm}_8, L^{dp}_8,  D_1, D_2, D_3}$  
is defined by the identity:
 \begin{itemize}
\item[ ] $x \lor y \lor [y \to (y \to x)] \approx x \to [x \lor (0 \to y)]$.    
\end{itemize} 
It is also defined
by the identity:
\begin{itemize}
\item[ ] $x \lor  [y \to \{0 \to (y \to x)\}] \approx x \lor y \lor (y \to x)$.   
\end{itemize}   
\end{Theorem}

\begin{Corollary}
The logic $\mathcal{L}(\mathbb{V}(\mathbf{\{L^{dm}_8, L^{dp}_8, 
 D_1, D_2, D_3\}}))$ is defined
by
 \begin{itemize}
\item[] $[\alpha \lor \beta \lor \{\beta \to (\beta \to \alpha)\}] \Leftrightarrow_H  [\alpha \to \{\alpha \lor (0 \to \beta)\}]$.    
\end{itemize} 
It is also defined by:
\begin{itemize}
\item[] $[\alpha \lor \{\beta \to (0 \to (\beta \to \alpha))\}] \Leftrightarrow_H [\alpha \lor \beta \lor (\beta \to \alpha)]$.
\end{itemize} 
\end{Corollary}

\begin{Theorem}
The variety generated by the algebras\\
 $\mathbf{ L^{dm}_7, L^{dp}_7, L^{dm}_8, L^{dp}_8, L^{dm}_9, L^{dp}_9}$, 
$\mathbf{L^{dm}_{10}, L^{dp}_{10}, D_1, D_2, D_3}$ is defined 
by the identity:
\begin{itemize}
\item[ ] $x \lor (x \to y)  \approx x \lor [(x \to y) \to 1]$.  
\end{itemize}  
\end{Theorem}

\begin{Corollary}
The logic $\mathcal{L}(\mathbb{V}(\mathbf{\{ L^{dm}_7, L^{dp}_7, L^{dm}_8, L^{dp}_8, L^{dm}_9, L^{dp}_9,
L^{dm}_{10}, L^{dp}_{10}, D_1, D_2, D_3\}}))$ is defined
  by
\begin{itemize}
\item[ ] $[\alpha \lor (\alpha \to \beta)]  \Leftrightarrow_H [\alpha \lor \{(\alpha \to \beta) \to \top\}]$. 
\end{itemize} 
\end{Corollary}

\begin{Theorem}
The variety generated by the algebras\\
 $\mathbf{2^e,  L^{dp}_7,  L^{dp}_8,  L^{dp}_9}$, 
$\mathbf{L^{dp}_{10}}$ is defined 
by the identities:
\begin{itemize}
\item[{\rm(1)}] $x \lor (x \to y)  \approx x \lor [(x \to y) \to 1]$,  
\item[{\rm(2)}] $x \lor x'  \approx 1$.
\end{itemize}  
\end{Theorem}

\begin{Corollary}
The logic $\mathcal{L}(\mathbb{V}(\mathbf{\{2^e, L^{dp}_7, L^{dp}_8, L^{dp}_9, 
 L^{dp}_{10} \}})$ is defined
  by
\begin{itemize}
\item[{\rm(1)}] $[\alpha \lor (\alpha \to \beta)]  \Leftrightarrow_H [\alpha \lor \{(\alpha \to \beta) \to \top\}]$,  
\item[{\rm(2)}]   $\alpha \ \lor \ \sim\alpha$.
\end{itemize}  
\end{Corollary}


\begin{Theorem}
The variety generated by the algebras\\
 $\mathbf{ L^{dm}_7,  L^{dm}_8,  L^{dm}_9,}$ 
$\mathbf{L^{dm}_{10}, D_1, D_2, D_3}$ is defined 
by the identities:
\begin{itemize}
\item[{\rm(1)}] $x \lor (x \to y)  \approx x \lor [(x \to y) \to 1]$,  
\item[{\rm(2)}] $x'' \approx x$.
\end{itemize}  
\end{Theorem}

\begin{Corollary}
The logic $\mathcal{L}(\mathbb{V}(\mathbf{\{ L^{dm}_7, L^{dm}_8, L^{dm}_9,
L^{dm}_{10}, D_1, D_2, D_3\}}))$ is defined
  by
\begin{itemize}
\item[{\rm(1)}] $[\alpha \lor (\alpha \to \beta)]  \Leftrightarrow_H [\alpha \lor \{(\alpha \to \beta) \to \top\}]$,  
\item[{\rm(2)}]   $\alpha \to_H \ \sim\sim\alpha$.
\end{itemize}  
\end{Corollary}

\begin{Theorem}
The variety generated by the algebras\\
 $\mathbf{L^{dm}_9, L^{dp}_9, L^{dm}_{10}, L^{dp}_{10},  D_1}$ is  
 defined 
by the identity:
 \begin{itemize}
\item[ ] $0 \to 1 \approx 0$    {\rm (FTF identity)}.  
\end{itemize}  
\end{Theorem}

\begin{Corollary}
The logic $\mathcal{L}(\mathbb{V}(\mathbf{\{ L^{dm}_9, L^{dp}_9, L^{dm}_{10}, L^{dp}_{10}, D_1\}}))$ is defined
  by
 \begin{itemize}
 
\item[] $(\bot \to \top) \Leftrightarrow_H \bot$  \qquad   {\rm (FTF identity)}.  
\end{itemize}  
\end{Corollary}

\begin{Theorem}
The variety generated by the algebras\\
 $\mathbf{{L}^{dm}_{10}, L^{dp}_{10},  D_1}$ is  
 defined by 
the identity:
 \begin{itemize}
\item[ ] $x \to y \approx y \to x$  \qquad   {\rm (commutative identity)}.  
\end{itemize}  
\end{Theorem}

\begin{Corollary}
The logic $\mathcal{L}(\mathbb{V}(\mathbf{\{L^{dm}_{10}, L^{dp}_{10}, D_1\}}))$ is defined
  by
 \begin{itemize}
\item[ ] $(\alpha \to \beta) \Leftrightarrow_H (\beta \to \alpha)$   {\rm (commutative identity)}.  
\end{itemize}  
\end{Corollary}





\begin{Theorem}
The variety $\mathbb{V}(\mathbf{C}_{20})$
 is defined 
by 
\begin{itemize}
\item[ ] $x^* \leq x'  $.
\end{itemize}
\end{Theorem}

\begin{Corollary}
The logic $\mathcal{L}(\mathbb{V}(\mathbf{C_{20}}))$
 is defined
  by
\begin{itemize}
\item[ ] $\alpha^* \to_H \ \sim\alpha  $.
\end{itemize}
\end{Corollary}


\begin{Theorem}
The variety $\mathbb{V}(\mathbf{D_2)}$
 is defined 
by 
\begin{itemize}
\item[ ] $(x \to y)^* \approx x \land y^* $.
\end{itemize}
\end{Theorem}

\begin{Corollary}
The logic $\mathcal{L}(\mathbb{V}(\mathbf{D_2})$ is defined
by
\begin{itemize}
\item[] $(\alpha \to \beta)^* \Leftrightarrow_H \alpha \land \beta^* $.
\end{itemize}
\end{Corollary}

\begin{Theorem}
The variety generated by the algebras in \\
$\{\mathbf{L^{dp}_i}: i=1, \dots, 8\} \cup \{\mathbf{L^{dm}_i}: i=1, \dots, 8\} \cup \mathbf{\{D_2\}}$ 
 is defined 
by the identity:
\begin{itemize}
\item[ ]   $(x \to y)^* \approx (x \land y^*)^{**}$.     
\end{itemize}
It is also defined by
\begin{itemize}
\item[ ]   $(0 \to 1)^* \approx 0$.     
\end{itemize}
\end{Theorem}

\begin{Corollary}
The logic $\mathcal{L}(\mathbb{V}(\{\mathbf{L^{dp}_i}: i=1, \dots, 8\} \cup \{\mathbf{L^{dm}_i}: i=1, \dots, 8\} \cup \mathbf{\{D_2\})}$ is defined
by
\begin{itemize}
\item[]   $(\alpha \to \beta)^* \to_H (\alpha \land \beta^*)^{**}$.
\end{itemize}
It is also defined by
 \begin{itemize}
\item[]   $(\bot \to \top)^* \Leftrightarrow_H \bot$.    
\end{itemize}
\end{Corollary}

\begin{Theorem}
The variety generated by the algebras 
 $\mathbf{L^{dp}_i}$, $ i=1, \dots, 8$, is    
 defined 
by the identities:
\begin{itemize}
\item[{\rm(1)}] $(x \to y)^* \approx (x \land y^*)^{**}$,
\item[{\rm(2)}]  $x \lor x'  \approx 1$.    
\end{itemize}
\end{Theorem}


\begin{Corollary}
The logic $\mathcal{L}(\mathbb{V}(\mathbf{\{ L^{dp}_i : i=1, \dots, 8\}}))$ is    
  defined
  by
\begin{itemize}
\item[{\rm(1)}]   $(\alpha \to \beta)^* \Leftrightarrow_H (\alpha \land \beta^*)^{**}$,
\item[{\rm(2)}]   $\alpha  \ \lor \ \sim\alpha$.  
\end{itemize}
\end{Corollary}


\begin{Theorem}
The variety generated by the algebras 
 $\mathbf{L^{dm}_i}$, $ i=1, \dots, 8$, and $\mathbf{D_2}$ is    
 defined 
by the identities:
\begin{itemize}
\item[{\rm(1)}] $(x \to y)^* \approx (x \land y^*)^{**}$,
\item[{\rm(2)}]  $x''  \approx x$.    
\end{itemize}
\end{Theorem}

\begin{Corollary}
The logic $\mathcal{L}(\mathbb{V}(\mathbf{\{ L^{dm}_i : i=1, \dots, 8\}\cup \{D_2\}}))$, is    
  defined
  by
\begin{itemize}
\item[{\rm(1)}]   $(\alpha \to \beta)^* \Leftrightarrow_H (\alpha \land \beta^*)^{**}$,
\item[{\rm(2)}]   $\alpha  \to_H \ \sim\sim \alpha$.  
\end{itemize}
\end{Corollary}

\begin{Theorem}
The variety generated by the algebras\\
 $\mathbf{L^{dm}_1,  L^{dp}_1, L^{dm}_2, L^{dp}_2, D_2}$ 
is  defined 
by the identity:
\begin{itemize}
\item[ ] $x \land z \leq y \lor (y \to z) $.           
\end{itemize}
\end{Theorem}

\begin{Corollary}
The logic $\mathcal{L}(\mathbb{V}(\mathbf{\{ L^{dm}_1, L^{dp}_1, L^{dm}_2, L^{dp}_2,  D_2\}}))$ is defined 
   by
\begin{itemize}
\item[] $(\alpha \land \gamma) \to_H \{\beta \lor (\beta \to \gamma)\}$.      
\end{itemize}
\end{Corollary}

\begin{Theorem}
The variety generated by \\
$\mathbf{L^{dm}_1, L^{dp}_1, L^{dm}_2, L^{dp}_2, L^{dm}_5, L^{dp}_5, L^{dm}_6, L^{dp}_6, D_2} $
 is defined 
by the identity:
\begin{itemize}
\item[ ]   $x \lor y \leq (x \to y) \to y$.    
\end{itemize}
\end{Theorem}

\begin{Corollary}
The logic $\mathcal{L}(\mathbb{V}(\mathbf{\{ L^{dm}_1, L^{dp}_1, L^{dm}_2, L^{dp}_2, L^{dm}_5, L^{dp}_5, 
L^{dm}_6, L^{dp}_6, D_2\}}))$ is defined
  by
\begin{itemize}
\item[ ]  $(\alpha \lor \beta) \to_H [(\alpha \to \beta) \to \beta]$.    
\end{itemize}
\end{Corollary}

\begin{Theorem}
The variety generated by \\
$\mathbf{ L^{dp}_1,  L^{dp}_2,  L^{dp}_5,  L^{dp}_6} $
 is defined 
by the identity:
\begin{itemize}
\item[{\rm(1)}]   $x \lor y \leq (x \to y) \to y$,    
\item[{\rm(2)}]  $x \lor x' \approx 1$.
\end{itemize}
\end{Theorem}

\begin{Corollary}
The logic $\mathcal{L}(\mathbb{V}(\mathbf{\{ L^{dp}_1, L^{dp}_2, L^{dp}_5, L^{dp}_6 \}}))$ is defined 
  by
\begin{itemize}
\item[{\rm(1)}]  $(\alpha \lor \beta) \to_H [(\alpha \to \beta) \to \beta]$,   
\item[{\rm(2)}]   $\alpha \ \lor \ \sim\alpha$.
\end{itemize}
\end{Corollary}


\begin{Theorem}
The variety generated by \\
$\mathbf{L^{dm}_1, L^{dm}_2,  L^{dm}_5,  L^{dm}_6,  D_2} $
 is defined 
by the identity:
\begin{itemize}
\item[{\rm(1)}]   $x \lor y \leq (x \to y) \to y$,     
\item[{\rm(2)}]  $x'' \approx x$.
\end{itemize}
\end{Theorem}

\begin{Corollary}
The logic $\mathcal{L}(\mathbb{V}(\mathbf{\{ L^{dm}_1, L^{dm}_2, L^{dm}_5, L^{dp}_6,  
D_2\}}))$ is defined
  by
\begin{itemize}
\item[{\rm(1)}]  $(\alpha \lor \beta) \to_H [(\alpha \to \beta) \to \beta]$,  
\item[{\rm(2)}]  $\alpha \to_H \ \sim\sim\alpha$.
\end{itemize}  
\end{Corollary}

The variety $\mathbb{V}(\mathbf{{\{D_1, D_2, D_3\}}})$ was axiomatized in Theorem \ref{DMBSH}.
Here are two more bases for it.

\begin{Theorem} \label{D1D2D3}
The variety $\mathbb{V}(\mathbf{{\{D_1, D_2, D_3\}}})$
 is defined 
by the identity:
\begin{itemize}
\item[ ]   $x \lor (y \to z) \approx (x \lor y) \to (x \lor z)$ \quad \rm(Strong JID).    
\end{itemize}
It is also defined by the identity:
\begin{itemize}
\item[ ] $x{'^*}{'^*} \approx x$.
\end{itemize}
\end{Theorem}
\begin{Corollary}
The logic $\mathcal{L}(\mathbb{V}(\mathbf{{\{D_1, D_2, D_3\}}}))$ is defined
  by
  
\begin{itemize}
\item[]   $(\alpha \lor (\beta \to \gamma)) \Leftrightarrow_H [(\alpha \lor \beta) \to (\alpha \lor \gamma)]$.     
\end{itemize}
It is also defined 
by the identity:
\begin{itemize}
\item[] $(\sim((\sim\alpha)^*))^* \Leftrightarrow_H \alpha$.
\end{itemize}
\end{Corollary}

\begin{Theorem}
The variety generated by 
$\mathbf{L^{dm}_2, L^{dp}_2, D_2} $ 
 is defined 
by the identity:
\begin{itemize}
\item[ ]   $(x \to y) \to x \approx x $.    
\end{itemize}
\end{Theorem}

\begin{Corollary}
The logic $\mathcal{L}(\mathbb{V}(\mathbf{\{ L^{dm}_2, L^{dp}_2,  D_2\}}))$ is defined
  by
 
	\begin{itemize}
\item[]   $((\alpha \to \beta) \to \alpha) \Leftrightarrow_H \alpha $.
\end{itemize}

\end{Corollary}

 $\mathbb{V}(\mathbf{D}_2)$ was axiomatized in Theorem \ref{bases}. 
 Here are some more bases for it.

\begin{Theorem}
$\mathbb{V}(\mathbf{D}_2)$  
 is defined by the identity:
\begin{itemize}
\item[ ]   $x \lor y \approx (x \to y) \to y$.     
\end{itemize}
It is also defined by the identities:
\begin{itemize}
\item[{\rm(1)}]  $x \lor (y \to z) \approx (x \lor y) \to (x \lor z)$,
\item[{\rm(2)}]   $(x \to y) \to x \approx x$.
\end{itemize}
It is also defined by the identity:
\begin{itemize}
\item[ ]   $x \lor (x \to y) \approx  x \lor ((x \lor y) \to 1)$.     
\end{itemize}
\end{Theorem}

\begin{Corollary}
The logic $\mathcal{L}(\mathbb{V}(\mathbf{D}_2))$ is axiomatized 
by
 
\begin{itemize}
\item[{}]   $(\alpha \lor \beta) \Leftrightarrow_H ((\alpha \to \beta) \to \beta)$.
\end{itemize}
This logic has an interesting property in that $\lor$ is definable in terms of $\to$. \\
It is also axiomatized  
by    
\begin{itemize}
\item[{\rm(1)}]  $(\alpha \lor (\beta \to \gamma)) \Leftrightarrow_H  [(\alpha \lor \beta) \to (\alpha \lor \gamma)]$,
\item[{\rm(2)}]   $((\alpha \to \beta) \to \alpha) \Leftrightarrow_H \alpha$.
\end{itemize}
It is also axiomatized  
by
\begin{itemize}
\item[]   $(\alpha \lor (\alpha \to \beta)) \Leftrightarrow_H  [\alpha \lor \{(\alpha \lor \beta) \to \top\}]$.
\end{itemize}
\end{Corollary}

\begin{Theorem}
The variety generated by \\
$\mathbf{L^{dm}_1, L^{dp}_1, L^{dm}_2, L^{dp}_2, L^{dm}_9, L^{dp}_9, D_1, D_2, D_3} $
 is defined by the identity:
\begin{itemize}
\item[]   $x \to (y \to z) \approx y \to (x \to z)$.     
\end{itemize}
\end{Theorem}

\begin{Corollary}
The logic $\mathcal{L}(\mathbb{V}(\mathbf{\{ L^{dm}_1, L^{dp}_1, L^{dm}_2, L^{dp}_2, L^{dm}_9, L^{dp}_9,
D_1, D_2, D_3\}}))$ is defined
by
\begin{itemize}
\item[]   $[\alpha \to (\beta \to \gamma)] \Leftrightarrow_H [\beta \to (\alpha \to \gamma)]$.     
\end{itemize}
\end{Corollary}

\begin{Theorem}
The variety generated by \\
$\mathbf{L^{dm}_1, L^{dp}_1, L^{dm}_2, L^{dp}_2, L^{dm}_5, L^{dp}_5,  D_2} $
 is defined by the identity:
\begin{itemize}
\item[ ]   $(x \to y) \to z \leq ((y \to x) \to z) \to z$  
\end{itemize}
\end{Theorem}

\begin{Corollary}
The logic $\mathcal{L}(\mathbb{V}(\mathbf{\{ L^{dm}_1, L^{dp}_1, L^{dm}_2, L^{dp}_2, L^{dm}_5, L^{dp}_5, 
D_2\}}))$ is defined
   by
  \begin{itemize}
\item[ ]   $[(\alpha \to \beta) \to \gamma] \to_H [((\beta \to \alpha) \to \gamma) \to \gamma]$ \quad   
\end{itemize}

\end{Corollary}

We note that 
a new extension of each of the logic defined in this section is obtained by adding the axiom $\alpha'' \Leftrightarrow_H \alpha$, as an extension of the logic $\mathcal{DMSH}$.  
Similarly, the addition of the axiom: $\alpha \lor \alpha'$ yields new extensions to the logics, over the logic $\mathcal{DPCSH}$, defined in the preceding corollaries.

We conclude this section by remarking that all the logics described in this section are discriminator logics and also are decidable.\\

\bigskip
\section{{\bf Logics corresponding to subvarieties of    
Regular De Morgan Semi-Heyting Algebras of level 1}} \label{Sec_Rdmsh}

\


In this section, we present axiomatizations for logics corresponding to several subvarieties of the variety 
$\mathbb{RDMSH}_1$ of regular De Morgan semi-Heyting algebras of level 1.
The algebraic results mentioned (or referred to) in this section were proved in \cite {sankappanavar2016}.
Recall that $\mathcal{DMSH}_1$ denotes the logic corresponding to the variety $\mathbb{DMSH}_1$.
The following corollary is immediate from Theorem \ref{completeness_DHMSH_extension} and definitions.

{\bf In what follows,  $\mathcal{V}$ (or $\mathcal{L}(\mathbb{{}V}))$ denotes the logic corresponding to the variety $\mathbb{V}$.}\\

Recall that the variety $\mathbb{DMSH}_1$ was defined in Section 10.

\begin{Corollary} \label{coro_160719_01} 
\begin{thlist}
\item The logic $\mathcal{DMSH}_1$ 
 is defined, 
relative to $\mathcal{DMSH}$, by

  
   \qquad $\alpha \ \land (\sim\alpha)^* \Leftrightarrow_H \ [\sim(\alpha \ \land (\sim\alpha)^*)]^*.$

		\item The logic $\mathcal{RDMSH}_1$ is defined, relative to $\mathcal{DMSH}_1$, by
		 
		$(\alpha \land \alpha^+)  \to_H (\beta \lor \beta^*)$.

		\item  The logic $\mathcal{RDMH}_1$ is defined, relative to $\mathcal{RDMSH}_1$, 
		by
		 
		$(\alpha \land \beta) \to \alpha $.

		\item  The logic $\mathcal{RDM}cm\mathcal{SH}_1$  is defined, relative to $\mathcal{RDMSH}_1$, by\\  
		 $(\alpha \to \beta) \to_H (\beta \to \alpha).$
	\end{thlist}
\end{Corollary}
It follows from Theorem \ref{Th11.A} that the logic $\mathcal{RDMH}_1$  is decidable.  However, teh following problem is still open.\\

{\bf PROBLEM}: Is the logic $\mathcal{RDMH}_1$  is decidable? \\

Let $\mathbf{L} \in \mathbb{DHMSH}$.  We say $\mathbf{L}$ is pseudocommutative if 
$\mathbf{L}$ satisfies the identity:\\
(PCM)   $x^* \to y^* \approx  y^* \to x^*$.

$\mathbb{RDM}{\rm pcm}\mathbb{SH}$ denotes the variety of regular De Morgan pseudocommutative semi-Heyting algebras.

The following corollary is immediate from Theorem \ref{completeness_DHMSH_extension} and the above definition.
\begin{Corollary}
The logic $\mathcal L(\mathbb{RDM}{\rm pcm}\mathbb{SH}_1)$ is defined 
by\\

 \qquad  $(\alpha^* \to \beta^*) \Leftrightarrow_H (\beta^* \to \alpha^*)$.
\end{Corollary}

\begin{Theorem} \cite {sankappanavar2016}
 $\mathbb{RDM}{\rm pcm}\mathbb{SH}_1= \mathbf{V(L_9^{dm}, L_{10}^{dm}, D_1)}$.
\end{Theorem}

\begin{Corollary}
The logic $\mathcal L(\mathbb{RDM}{\rm pcm}\mathbb{SH}_1)$ is decidable.

\end{Corollary}


 In the rest of this section, unless otherwise stated, the phrase  ``defined, modulo 
$\mathbb{RDMSH}_1$, by''  is abbreviated to the phrase ``defined by''
  in the context of varieties.  
Similarly,  the phrase 
``defined, as an extension of the logic 
$\mathcal{RDMSH}_1$, by'' is also abbreviated to the phrase ``defined by'' in the case of logics.

 The theorems that appear in the rest of this section were proved in \cite{sankappanavar2016}.
Each of the corollaries given below  
follows from the theorem immediately preceding it and Theorem \ref{teo_040417_01}.

Here is another axiomatization for $\mathbb{RDM}{\rm pcm}\mathbb{SH}$.

\begin{Theorem}
 The variety 
 $\mathbb{RDM}{\rm pcm}\mathbb{SH}$
 is defined  
 by 
\begin{itemize}
\item[ ] \begin{itemize}
\item[ ]  $(x \to y)^* \approx  (y \to x)^*$.
\end{itemize}
\end{itemize}
\end{Theorem}

\begin{Corollary}
The logic $\mathcal L(\mathbb{RDM}{\rm pcm}\mathbb{SH})$ is defined  
by\\

$(\alpha \to \beta)^* \Leftrightarrow_H (\beta \to \alpha)^*$.
\end{Corollary}

\begin{Theorem}
 The variety $\mathbb{V}(\mathbf{L_1^{dm}, L_2^{dm}, L_3^{dm}, L_4^{dm},  D_2, D_3})$ is defined 
 by 
\begin{itemize}
\item[ ]  

        $(0 \to 1)^+  \to \ [\sim\{(0 \to 1){^*}\}]^* \approx  0 \to 1$.
\end{itemize}
\end{Theorem}

\begin{Corollary}
The logic $\mathcal L(\mathbb{V}\mathbf{(L_1^{dm}, L_2^{dm}, L_3^{dm}, L_4^{dm},  D_2, D_3)})$ is defined 
by 
\begin{itemize}
\item[]    $((\bot \to \top)^+  \to \ [\sim\{(\bot \to \top){^*}\}]^*) \Leftrightarrow_H  (\alpha \to \top)$.
\end{itemize}
\end{Corollary}

The variety $\mathbf{V(D_1, D_2, D_3}) (= \mathbf{DQDBSH})$ was axiomatized earlier.  
Here are some more bases for $\mathbf{V(D_1, D_2, D_3})$.

\begin{Theorem}
Each of the following identities is a base for
 the variety $\mathbf{V(D_1, D_2, D_3})$: 
\begin{itemize}
\item[{\rm(1)}]  $x \to y  \approx  y^* \to x^*$  {\rm (Law of contraposition)},
\item[{\rm(2)}]   $ [\{x \lor (x \to y^*)\} \to (x \to y^*)] \lor (x \lor y^*)=1$. 
\end{itemize}
\end{Theorem} 

\begin{Corollary}  
\begin{thlist}
\item[1]  The logic $\mathcal L(\mathbb{V}\mathbf{(D_1,  D_2, D_3)})$ is defined 
by \\
 $(\alpha \to \beta)  \Leftrightarrow_H (\beta^* \to \alpha^*)$.   
\item[2]
The logic $\mathcal L(\mathbb{V}\mathbf{(D_1,  D_2, D_3)})$ is also defined 
by\\
  $ [\{\alpha \lor (\alpha \to \beta^*)\} \to (\alpha \to \beta^*)] \lor (\alpha \lor \beta^*)$.
\end{thlist}
\end{Corollary}

\begin{Theorem} 
 The variety 
 $\mathbb{V}\mathbf{(L_1^{dm}, L_2^{dm}, L_5^{dm}, L_6^{dm}, L_9^{dm}, D_1, D_2, D_3})$ is defined
 by 
\begin{itemize}
\item[ ]  $x \to y^*  \approx  y \to x^*$.
\end{itemize}
\end{Theorem}

\begin{Corollary}
The logic $\mathcal{L}(\mathbb{V}\mathbf{(L_1^{dm}, L_2^{dm}, L_5^{dm}, L_6^{dm}, L_9^{dm}, D_1, D_2, D_3}))$ is defined 
by\\

 $(\alpha \to \beta^*)  \Leftrightarrow_H (\beta \to \alpha^*)$.  

\end{Corollary}

\begin{Theorem}
 The variety $\mathbf{V(L_7^{dm}, L_8^{dm}, L_9^{dm}, L_{10}^{dm}, D_1, D_2, D_3})$ is defined 
 by 
\begin{itemize}
\item[ ]  $x \lor  (x \to y)  \approx  x \lor [(x \to y) \to 1]$.
\end{itemize}
\end{Theorem}

\begin{Corollary}
The logic $\mathcal{L}(\mathbb{V}\mathbf{(L_7^{dm}, L_8^{dm}, L_9^{dm}, L_{10}^{dm}, D_1, D_2, D_3}))$ is defined 
by
\begin{itemize}
\item[]  $[\alpha \lor  (\alpha \to \beta)]  \Leftrightarrow_H  [\alpha \lor \{(\alpha \to \beta) \to \top\}]$,
\end{itemize}
\end{Corollary}

\begin{Theorem}
 The variety $\mathbf{V(L_7^{dm}, L_8^{dm},  D_2})$ is defined 
 by 
\begin{itemize}
\item[{\rm(1)}]  $x \lor  (x \to y)  \approx  x \lor [(x \to y) \to 1]$,
\item[{\rm(2)}] $ (0 \to 1)^{**} \approx 1$.
\end{itemize}
\end{Theorem}

\begin{Corollary}
The logic $\mathcal{L}(\mathbb{V}\mathbf{(L_7^{dm}, L_8^{dm}, D_2}))$ is defined 
by
\begin{itemize}
\item[{\rm(1)}]  
 $[\alpha \lor  (\alpha \to \beta)]  \Leftrightarrow_H  [\alpha \lor \{(\alpha \to \beta) \to \top\}]$,

\item[{\rm(2)}] 
   $(\bot \to \top)^{**}$.
\end{itemize}
\end{Corollary}

\begin{Theorem} 
 The variety $\mathbf{V(2^e, L_7^{dm}, L_8^{dm}, L_9^{dm},  L_{10}^{dm}})$ is\\ defined 
 by 
\begin{itemize}
\item[{\rm(1)}]  $x \lor  (x \to y)  \approx  x \lor [(x \to y) \to 1]$,
\item[{\rm(2)}]  $x{^*}' \approx  x^{**}$ \quad {\rm ($\star$-regular)}.
\end{itemize}
\end{Theorem}

We caution the reader that in \cite{sankappanavar2011expansions}, (2) was referred to as ``regular''.

\begin{Corollary}
The logic $\mathcal{L}(\mathbb{V}\mathbf{(2^e, L_7^{dm}, L_8^{dm}, L_9^{dm},  L_{10}^{dm}}))$ is defined 
by
\begin{itemize}
\item[{\rm(1)}]  $[\alpha \lor  (\alpha \to \beta)] \Leftrightarrow_H  [\alpha \lor \{(\alpha \to \beta) \to 1\}]$,
\item[{\rm(2)}]   $\sim(\alpha{^*}) \Leftrightarrow_H  \alpha^{**}$.  
\end{itemize}
\end{Corollary}

\begin{Theorem}
 The variety 
 $\mathbf{V(2^e,  L_9^{dm},  L_{10}^{dm}})$ is defined 
 by 
\begin{itemize}
\item[{\rm(1)}]  $x \lor  (x \to y)  \approx  x \lor [(x \to y) \to 1]$,
\item[{\rm(2)}]  $x{^*}' \approx  x^{**}$,
\item[{\rm(3)}] $ (0 \to 1) \lor (0 \to 1)^* \approx 1$.
\end{itemize}
\end{Theorem}

\begin{Corollary}
The logic $\mathcal{L}(\mathbb{V}\mathbf{(2^e, L_9^{dm},  L_{10}^{dm}}))$ is defined 
by
\begin{itemize}
\item[{\rm(1)}]   $[\alpha \lor  (\alpha \to \beta)] \Leftrightarrow_H  [\alpha \lor \{(\alpha \to \beta) \to \top\}]$,
\item[{\rm(2)}]  $\sim(\alpha{^*}) \Leftrightarrow_H  \alpha^{**}$,  
\item[{\rm(3)}] $ (\bot \to \top) \lor (\bot \to \top)^* $.
\end{itemize}
\end{Corollary}


\begin{Theorem}
 The variety \\
 $\mathbf{V(L_9^{dm},  L_{10}^{dm}})$ is defined 
 by 
\begin{itemize}
\item[{\rm(1)}]  $x \lor  (x \to y)  \approx  x \lor [(x \to y) \to 1]$,
\item[{\rm(2)}]  $x{^*}' \approx  x^{**}$,
\item[{\rm(3)}] $ (0 \to 1)^{*} \approx 1$.
\end{itemize}
\end{Theorem}

\begin{Corollary}
The logic $\mathcal{L}(\mathbb{V}\mathbf{(L_9^{dm},  L_{10}^{dm}}))$ is defined 
by
\begin{itemize}
\item[{\rm(1)}]   $[\alpha \lor  (\alpha \to \beta)] \Leftrightarrow_H  [\alpha \lor \{(\alpha \to \beta) \to 1\}]$,
\item[{\rm(2)}]   $\sim(\alpha{^*}) \Leftrightarrow_H  \alpha^{**}$,  
\item[{\rm(3)}] $ (\bot \to \top)^* $.
\end{itemize}
\end{Corollary}

\begin{Theorem}
 The variety \\
 $\mathbf{V(L_1^{dm}, L_2^{dm}, L_3^{dm}, L_4^{dm}, L_5^{dm}, L_6^{dm}, L_7^{dm}, L_8^{dm}})$ is defined 
 by 
\begin{itemize}
\item[{\rm(1)}] $x{^*}' \approx  x^{**}$,
\item[{\rm(2)}] $ (0 \to 1)^{**} \approx 1$.
\end{itemize}
\end{Theorem}

\begin{Corollary}
The logic $\mathcal{L}(\mathbb{V}\mathbf{( (L_1^{dm}, L_2^{dm}, L_3^{dm}, L_4^{dm}, L_5^{dm}, L_6^{dm}, L_7^{dm}, L_8^{dm}})$ is defined 
by
\begin{itemize}
\item[{\rm(1)}]  $\sim(\alpha{^*}) \Leftrightarrow_H  \alpha^{**}$.  
\item[{\rm(2)}] $ (\bot \to \top)^{**} $.
\end{itemize}
\end{Corollary}

\begin{Theorem}
 The variety $\mathbf{V(L_1^{dm}, L_2^{dm}, L_3^{dm}, L_4^{dm}, D_2})$ is  defined 
 by 
\begin{itemize}
\item[{\rm(1)}]  $(0 \to 1) \lor (0 \to 1)^* \approx 1$,
\item[{\rm(2)}] $ (0 \to 1)^{**} \approx 1$.
\end{itemize}
\end{Theorem}

\begin{Corollary}
The logic $\mathcal{L}(\mathbb{V}\mathbf{( (L_1^{dm}, L_2^{dm}, L_3^{dm}, L_4^{dm}, D_2}))$ is defined 
by
\begin{itemize}
\item[{\rm(1)}]   $(\bot \to \top) \lor (\bot  \to \top)^* $,
\item[{\rm(2)}] $ (\bot \to \top)^{**} $.
\end{itemize}
\end{Corollary}


\begin{Theorem}
 The variety $\mathbf{V(L_1^{dm}, L_3^{dm}, D_1, D_2, D_3})$ is 
 defined  
 by 
\begin{itemize}
\item[{\rm(1)}]  $x \lor (y \to x)  \approx (x \lor y) \to x$,
\item[{\rm(2)}] $(0 \to 1) \lor (0 \to 1)^* \approx 1 $.
\end{itemize}
\end{Theorem}

\begin{Corollary}
The logic $\mathcal{L}(\mathbb{V}\mathbf{( (L_1^{dm}, L_3^{dm}, D_1, D_2, D_3}))$ is defined 
by
\begin{itemize}
\item[{\rm(1)}]  $[\alpha \lor  (\beta \to \alpha)] \Leftrightarrow_H  [(\alpha \lor \beta) \to \alpha]$,
\item[{\rm(2)}]   $(\bot \to \top) \lor (\bot  \to \top)^*. $
\end{itemize}
\end{Corollary}


\begin{Theorem}
 The variety $\mathbf{V(L_1^{dm},  L_3^{dm}, D_2})$ is defined 
 by 
\begin{itemize}
\item[{\rm(1)}]  $x \lor (y \to x)  \approx (x \lor y) \to x$,
\item[{\rm(2)}] $(0 \to 1) \lor (0 \to 1)^* \approx 1 $,
\item[{\rm(3)}]  $(0 \to 1)^{**} \approx 1$.
\end{itemize}
\end{Theorem}

\begin{Corollary}
The logic $\mathcal{L}(\mathbb{V}\mathbf{( (L_1^{dm}, L_3^{dm}, D_2}))$ is defined 
by
\begin{itemize}
\item[{\rm(1)}]  $[\alpha \lor  (\beta \to \alpha)] \Leftrightarrow_H  [(\alpha \lor \beta) \to \alpha]$,
\item[{\rm(2)}]   $(\bot \to \top) \lor (\bot  \to \top)^*, $
\item[{\rm(3)}] $(\bot \to \top)^{**}$.
\end{itemize}
\end{Corollary}


\begin{Theorem}
 The variety $\mathbf{V(L_1^{dm}, L_2^{dm}, L_8^{dm}, D_1, D_2, D_3})$ is defined 
 by 
\begin{itemize}
\item[ ]  $y \lor (y \to (x \lor  y)) \approx (0 \to x) \lor (x \to y)$.
\end{itemize}
\end{Theorem}

\begin{Corollary}
The logic  $\mathcal{L}(\mathbb{V}\mathbf{( (L_1^{dm}, L_2^{dm}, L_8^{dm}, D_1, D_2, D_3}))$ is defined 
by
\begin{itemize}
\item[ ]  $[\beta \lor (\beta \to (\alpha \lor  \beta))] \Leftrightarrow_H [(\bot \to \alpha) \lor (\alpha \to  \beta)]$.
\end{itemize}
\end{Corollary}

\begin{Theorem}
 The variety $\mathbf{V(L_8^{dm}, D_1, D_2, D_3)}$ is defined 
 by 
\begin{itemize}
\item[ ]  $x \lor [y \to (0 \to (y \to x))] \approx x \lor y \lor (y \to x)$.
\end{itemize}
\end{Theorem}

\begin{Corollary}
The logic $\mathcal{L}(\mathbb{V}\mathbf{(L_8^{dm}, D_1, D_2, D_3}))$ is defined 
by
\begin{itemize}
\item[ ]   $[\alpha \lor \{\beta \to (\bot \to (\beta \to  \alpha))\}] \Leftrightarrow_H [\alpha \lor \beta \lor (\beta \to \alpha)]$.
\end{itemize}
\end{Corollary}



\begin{Theorem}
 The variety $\mathbb{V}(\mathbf{C}^{dm})$ is defined 
 by 
\begin{itemize}
\item[ ]  $x \land x' \leq y \lor y'$  {\rm (Kleene identity)}.
\end{itemize}
\end{Theorem}

\begin{Corollary}
The logic $\mathcal{L}(\mathbb{V}\mathbf{(C^{dm}))}$ is defined 
by
\begin{itemize}
\item[ ]   $(\alpha \ \land \ \sim\alpha) \to_H (\beta \ \lor \ \sim\beta$ \quad  (Kleene identity)
\end{itemize}
\end{Corollary}

\begin{Theorem}
 The variety $\mathbb{V}(\mathbf{L}_{10}^{dm})$ is defined 
 by 
\begin{itemize}
\item[{\rm(1)}]  $x \land x' \leq y \lor y'$  {\rm (Kleene identity)},
\item[{\rm(2)}]  $x \to y \approx y \to x$.
\end{itemize}
\end{Theorem}

\begin{Corollary}
The logic $\mathcal{L}(\mathbb{V}\mathbf{(L}_{10}^{dm}))$ is defined 
by
\begin{itemize}
\item[{\rm(1)}]    $(\alpha \ \land \ \sim\alpha) \to_H (\beta \ \lor \ \sim\beta$) \quad  (Kleene identity),
\item[{\rm(2)}]     $\alpha \to \beta \Leftrightarrow_H  \beta \to \alpha$.
\end{itemize}
\end{Corollary}

 $\mathbf{V(D_2)}$ was axiomatized in Section \ref{S7}. 
 Here are some more bases for it, but relative to $\mathbf{RDMH_1}$.

\begin{Theorem}
Each of the following identities is a base for 
 $\mathbf{V(D_2)}$, mod  $\mathbf{RDMH_1}$: 
  
\begin{itemize}
\item[{\rm(1)}]  $[y \to \{0 \to (y \to x)\}] \approx  y \lor (y \to x)$.

\item[{\rm(2)}]  $x \lor (y \to z) \approx (x \lor  y) \to (x \lor z)$.


\item[{\rm(3)}]  $[\{x \lor (x \to y^*)\} \to (x \to y^*)] \lor x \lor y^* \approx 1.$
\end{itemize}
\end{Theorem}

\begin{Corollary} Each of the following axioms defines the logic $\mathcal{L}(\mathbb{V}\mathbf{(D}_2)$, relative to $\mathcal{RDMH}_1$:
\begin{itemize}
\item[{\rm(1)}]  $ [\beta \to \{\bot \to (\beta \to \alpha)\}] \Leftrightarrow_H  [\beta \lor (\beta \to \alpha)]$,

\item[{\rm(2)}]  $[\alpha \lor (\beta \to \gamma)] \Leftrightarrow_H [(\alpha \lor  \beta) \to (\alpha \lor \gamma)]$,


\item[{\rm(3)}]  $[\{\alpha \lor (\alpha \to \beta^*)\} \to (\alpha \to \beta^*)] \lor \alpha \lor \beta^*.$
\end{itemize}
\end{Corollary}


$\mathbb{V}(\mathbf{D}_1)$ was axiomatized in Section \ref{S8}. 
Here are more bases for it.  
Let $\mathbb{RDM}{\rm cm}\mathbb{SH}_1$ denote the subvariety of 
$\mathbb{RDM}\mathbb{SH}_1$ defined by: $x \to y \approx y \to x$.

\begin{Theorem}
Each of the following identities is an equational base for 
 $\mathbb{V}(\mathbf{D}_1)$, mod  $\mathbb{RDM}{\rm cm}\mathbb{SH}_1$:

\begin{itemize}

\item[{\rm(1)}] $y \lor (y \to (x \lor  y)) \approx (0 \to x) \lor (x \to y)$,  

\item[{\rm(2)}]   $x \lor [y \to (y \to x)^*] \approx x \lor y \lor (y \to x)$,

\item[{\rm(3)}]   $[\{x \lor (x \to y^*)\} \to (x \to y^*)] \lor x \lor y^* \approx 1$,

\item[{\rm(4)}]  $x \lor (y \to z) \approx (x \lor  y) \to (x \lor z)$.
 \end{itemize}
\end{Theorem}

\begin{Corollary}
Each of the following axioms defines the logic $\mathcal{L}(\mathbb{V}\mathbf{(D}_1))$, relative to $\mathcal{RDM}{\rm cm}\mathcal{SH}_1$:
\begin{itemize}

\item[{\rm(1)}] $[\beta \lor (\beta \to (\alpha \lor  \beta))] \Leftrightarrow_H [(\bot \to \alpha) \lor (\alpha \to \beta)]$,  

\item[{\rm(2)}]   $[\alpha \lor \{\beta \to (\beta \to \alpha)^*\}] \Leftrightarrow_H [\alpha \lor \beta \lor (\beta \to \alpha)]$,

\item[{\rm(3)}]   $[\{\alpha \lor (\alpha \to \beta^*)\} \to (\alpha \to \beta^*)] \lor \alpha \lor \beta^* $,

\item[{\rm(4)}]  $[\alpha \lor (\beta \to \gamma)] \Leftrightarrow_H [(\alpha \lor  \beta) \to (\alpha \lor \gamma)]$. \\
 \end{itemize}
\end{Corollary}

We conclude this section with the remark that all logics introduced in this section are discriminator logics as their corresponding varieties are discriminator varieties.\\



\medskip

\section{Extensions of the logic $\mathcal{JIDSH}_1$} \label{Sec_JID}

\ 

Algebras closely related to $\mathbb{DS}\rm t\mathbb{SH}$-algebras, called ``JI-distributive semi-Heyting algebras'', were introduced in \cite{sankappanavar2017}. 

An algebra $\mathbf{A}$ in $\mathbb{DQDSH}$ is  JI-distributive if $\mathbf{A}$ satisfies:

\begin{itemize}
\item[{\rm (JID)}]  $x' \lor (y \to z) \approx (x' \lor y) \to (x' \lor z)$ \quad ((restricted) {{\bf J}oin over {\bf I}mplication \bf D}istributivity).
\end{itemize}
\smallskip
We note that the identity (JID) is obtained by slightly weakening the identity (Strong JID) that has appeared earlier in Theorem \ref{D1D2D3}.
Let $\mathbb{JIDSH}$ denote the variety of JI-distributive $\mathbb{DQDSH}$-algebras and
let $\mathbb{JIDSH}_1$ (or  $\mathbb{JID}_1$, for short) denote the subvariety of $\mathbb{JIDSH}$ of level 1. 


{\bf In what follows,  $\mathcal{V}$ (or $\mathcal{L}(\mathbb{{}V}))$ denotes the logic corresponding to the subvariety $\mathbb{V}$.}

In this section we present axiomatizations of the logics 
corresponding to the subvarieties of $\mathbb{JIDSH}_1$ which we denote simply by $\mathbb{JID}_1$. 

\begin{Corollary}
The logic $\mathcal {JID}_1$ corresponding to $\mathbb{JID}_1$ is defined, as an extension of $\mathcal {DQDSH}$,  by 
\begin{enumerate}
\item[{\rm(a)}] $(\sim\alpha \lor (\beta \to \gamma)) \Leftrightarrow_H ((\sim\alpha \lor \beta) \to (\sim\alpha \lor \gamma))$,
\item[{\rm(b)}] 
                            
                            $\alpha \ \land (\sim\alpha)^* \Leftrightarrow_H \ [\sim(\alpha \ \land \ (\sim\alpha)^*)]^*.$
\end{enumerate}
\end{Corollary}

Let $\mathbb{DS}\rm t$ [$\mathbb{DS}\rm t\mathbb{H}$] denote the variety of dually Stone semi-Heyting [Heyting] algebras.  
The following theorem was proved in \cite[Corollary 5.10]{sankappanavar2017}. 

\begin{Theorem}
$\mathbb{JID}_1 = \mathbb{DS}\rm t
\lor \mathbb{V}(\mathbf{D_1}, \mathbf{D_2}, \mathbf{D_3}). $ 
\end{Theorem}



The preceding Theorem leads us naturally to raise the following open problems.\\


\noindent{\bf PROBLEM }: Is the logic $\mathcal {DS}t$ decidable? \\


 We make the following conjecture about an extension of $\mathcal {DS}t$:\\

\noindent{\bf CONJECTURE}:  The logic $\mathcal {DS}t\mathcal{H}$ is decidable. \\


We let $\mathbb{JIDL}_1$ denote the subvariety of $\mathbb{JID}$-algebras of level 1 defined by 
\begin{itemize}
\item[{\rm (L)}]  $(x \to y) \lor (y \to x) \approx 1$.  
\end{itemize}

\medskip
The results in the rest of this section depend on the corresponding algebraic results of \cite{sankappanavar2017}.
The relevant results, however, are stated here for the convenience of the reader.  The following corollary is immediate from the above definitions in view of Theorem \ref{teo_040417_01}.

\begin{Corollary}
The logic $\mathcal {JIDL}_1$ corresponding to $\mathbb{JIDL}_1$ is defined, modulo $\mathcal {JID}_1$,  by 
\begin{enumerate}
\item[] $(\alpha \to \beta) \vee (\beta \to \alpha)$.
\end{enumerate}
\end{Corollary}

Let $\mathbb{DS}\rm t\mathbb{L}$ denote the subvariety of $\mathbb{DS}\rm t$ defined by the identity (L) and  $\mathbb{DS}\rm t\mathbb{HC}$ denote the subvariety of $\mathbb{DS}\rm t\mathbb{H}$ generated by its chains.  

\begin{Theorem} {\rm \cite{sankappanavar2017}.} \label{C1} 
$\mathbb{JIDL}_1 = \mathbb{DS}\rm t\mathbb{HC}  \lor \mathbb{V}(\mathbf{D_2)}.$  
\end{Theorem}







For $n \in \mathbb{N}$, let $\mathbf{C_n^{dp}}$ denote the n-element $\mathbb {DS}\rm t\mathbb{H}$-chain (= $\mathbb {DPCH}$-chain)
denotes the variety generated by
$\mathbf{C_n^{dp}}$.  (Note that $\mathbf {C_3^{dp}}$= $\mathbf {L_1^{dp}}$.) 

Since the variety of Boolean algebras is the smallest non-trivial subvariety of $\mathbb{JIDL}_1$, we denote by 
$\mathbf {L_V}^+(\mathbb{JIDL}_1)$ the latttice of non-trivial subvarieties of $\mathbb{JIDL}_1$.

The following theorem was proved in \cite[Corollary 7.1] {sankappanavar2017}.
\begin{Theorem} \label{Th1}
\begin{thlist}
\item[1]
  $\mathbf{L_V}^+\mathbb{(JIDL}_1) \cong 
   [\mathbf{(\omega + 1)} \times \mathbf{2}]$, where $\times$ represents the direct product. 
\item[2]   $\mathbb{JIDL}_1$ and $\mathbb{DS}\rm t\mathbb{HC}$ are the only two elements of infinite height in the lattice $\mathbf{L_V}^+\mathbb{(JIDL}_1)$.
\item[3] 
    $\mathbb{V} \in  \mathbf{L_V}^+ \mathbb{(JIDL}_1)$ is of finite height if and only if $\mathbb{V}$ is either $\mathbb{V}(\mathbf{D_2)}$ or $\mathbb{V}(\mathbf{C}_n^{dp})$, for some $n \in \mathbb{N} \setminus \{1\}$, or $\mathbb{V}(\mathbf{C}_m^{dp})  \, \lor \, \mathbb{V}(\mathbf{D_2)}$, for some $m \in \mathbb{N} \setminus \{1\}$.     
\end{thlist}
\end{Theorem}

The following corollary is immediate from the preceding theorem and Theorem \ref{completeness_DHMSH_extension}. 

\begin{Corollary} 
The logic $\mathcal {JIDL}_1$ has the finite model property and hence it is decidable.
\end{Corollary}

Bases for all subvarieties of $\mathbb {JIDL}_1$ were given in \cite{sankappanavar2017}. 
The theorems presented below are taken from \cite{sankappanavar2017} and each of the corollaries given below follows from the theorem that precedes it and Theorem \ref{completeness_DHMSH_extension}. 

 In the rest of this section, the phrase ``defined, modulo $\mathbb{JIDL}_1$, by'' is abbreviated to ``defined by'', 
 in the context of varieties. 
Similarly,  the phrase 
``defined, as an extension of the logic 
$\mathcal{JIDL}_1$, by'' is also abbreviated to the phrase ``defined by'' in the case of logics.



 The theorems that appear below were proved in \cite{sankappanavar2016}.
Each of the corollaries given below  
follows from the theorem immediately preceding it and Theorem \ref{teo_040417_01}.

\begin{Theorem}
The variety  $\mathbb {DS}\rm t\mathbb{HC}$ is defined 
by 
\begin{enumerate}
\item[ ] $x \lor x' \approx 1$.
\end{enumerate}
\end{Theorem}

\begin{Corollary}  \label{cor_6_2}
The logic $\mathcal {DS}\rm t\mathcal{HC}$ 
is defined 
by 
\begin{enumerate}
\item[ ] $\alpha \ \lor  \ \sim\alpha$.  
\end{enumerate}
\end{Corollary}

The variety $\mathbb{V}(\mathbf{D}_2)$ was axiomatized earlier.  Here is another one. 

\begin{Theorem}\label{cor_6_3}
The variety $\mathbb{V}(\mathbf{D}_2)$ is defined 
 by 
\begin{enumerate}
\item[ ] $x'' \approx x$.
\end{enumerate}
\end{Theorem}

\begin{Corollary}
The logic $\mathcal{L}(\mathbb{V}(\mathbf{D}_2))$ 
 is defined  
 by 
\begin{enumerate}
\item[ ] $\alpha \Leftrightarrow_H  \ \sim\sim\alpha$.  
\end{enumerate}
\end{Corollary}

Let $n \in \mathbb{N}$ such that $n \geq 2$.
 
\begin{Theorem} \label{cor_6_4}
  The variety $\mathbb{V}(\mathbf{C}_n^{dp})\ \lor \ \mathbb{V}(\mathbf{D}_2)$ is defined 
by 
\begin{enumerate}
\item[{\rm (E$_n$)}] $x_1  \lor  x_2  \lor \cdots \lor x_n  \lor  (x_1 \to x_2)  \lor (x_2 \to x_3) \lor \cdots 
\lor (x_{n-1} \to x_n)=1.$
\end{enumerate}
\end{Theorem}

\begin{Corollary}
The logic $\mathcal{L}(\mathbb{V}(\mathbf{C_n}^{dp})\ \lor \ \mathbb{V}(\mathbf{D}_2))$  
 is defined 
 by 
\begin{enumerate}
\item[{\rm ($\mathcal E_n$)}] $\alpha_1  \lor  \alpha_2  \lor \cdots \lor \alpha_n  \lor  (\alpha_1 \to \alpha_2)  \lor (\alpha_2 \to \alpha_3) \lor \cdots 
\lor (\alpha_{n-1} \to \alpha_n).$
\end{enumerate}
\end{Corollary}


\begin{Theorem}
 The variety $\mathbb{V}(\mathbf{C}_n^{dp})$ is defined  
by 
\begin{enumerate}
\item[{\rm(1)}] $x \lor x' \approx 1,$
\item[{\rm(2)}] $x_1  \lor  x_2  \lor \cdots \lor x_n  \lor  (x_1 \to x_2)  \lor (x_2 \to x_3) \lor \cdots 
\lor (x_{n-1} \to x_n)=1.$
\end{enumerate}
\end{Theorem}

\begin{Corollary}
The logic $\mathcal{L}(\mathbb{V}(\mathbf{C}_n^{dp})) $    
 is defined 
 by 
\begin{enumerate}
\item[{\rm (a)}] $\alpha \ \lor \ \sim\alpha ,$
\item[{\rm (C$_n$)}] $\alpha_1  \lor  \alpha_2  \lor \cdots \lor \alpha_n  \lor  (\alpha_1 \to \alpha_2)  \lor (\alpha_2 \to \alpha_3) \lor \cdots 
\lor (\alpha_{n-1} \to \alpha_n).$
\end{enumerate}
\end{Corollary}


\begin{Theorem}
The variety $\mathbb {V}(\mathbf{C_3^{dp})}\,\lor \, \mathbb {V}(\mathbf{D_2)}$ is defined 
by 
\begin{enumerate}
\item[] $x \land x^+ \ \leq \ y \lor y^*$   {\rm (Regularity)}.
\end{enumerate}
It is also defined  
by 
\begin{enumerate}
\item[] $x \land x'  \leq y \lor y^*$.
\end{enumerate}
\end{Theorem}

\begin{Corollary}
The logic $\mathcal{L}(\mathbb{V}(\mathbf{Ch}_3^{dp})\ \lor \ \mathbb{V}(\mathbf{D}_2))$  
 is defined 
 by 
\begin{enumerate}
\item[]  $(\alpha \land \alpha^+) \to_H (\beta \lor \beta^*)$.   
\end{enumerate}
It is also defined  
by 
\begin{enumerate}
\item[]  $(\alpha \ \land \ \sim\alpha) \to_H (\beta \lor \beta^*)$. 
\end{enumerate}
\end{Corollary}

The variety $\mathbb{V}(\mathbf{L}_1^{dp})$ is axiomatized in   Corollary 
\ref{3bases}.  
Here is yet another axiomatization for it.

\begin{Theorem}
The variety $\mathbb{V}(\mathbf{L}_1^{dp})$ is defined 
by
\begin{itemize}
\item[{\rm(1)}] $x \land x^+ \leq y \lor y^*$  {\rm (Regularity)},
\item[{\rm(2)}] $x{^*}'=x^{**}$.
\end{itemize}
\end{Theorem}

\begin{Corollary} \label{cor_6_5}
The logic $\mathcal{L}(\mathbb{V}(\mathbf{L}_1^{dp})) $    
 is defined 
 by 
\begin{enumerate}
\item[{\rm(1)}] \quad  $(\alpha \land \alpha^+) \to_H (\beta \lor \beta^*)$,   
\item[{\rm(2)}] \quad $\sim\alpha{^*} \to_H \alpha^{**}$.
\end{enumerate}
\end{Corollary}

We note that the extensions of $\mathcal {JIDL}_1$ are all decidable. 

We conclude this aection with a partial poset of subvarieties of $\mathbb{DQDSH}$ discussed in the last sections.  Its dual will give the partial poset of the axiomatic extensions of the logic $\mathcal{DQDSH}$.  Note that the links in the poset do not, in general, represent the covers.

\newpage
\begin{center}
{\bf PARTIAL POSET OF SUBVARIETIES OF $\mathbb{DQDSH}$}
\end{center}

{\small
	\setlength{\unitlength}{1cm}
	\begin{picture}(20,16)(0,0)

\put(2,6){\circle*{0.15}}
\put(2.5,5.8){$\mathbb{RDMSH}_1=\mathbb{RDMS}\rm t\mathbb{SH}_1 = \mathbb V(C_{10}^{dm} \cup \{D_1, D_2, D_3\})$}	
\put(2,8){\circle*{0.15}}
 \put(0,8){$\mathbb{DM}\mathbb{SH}_1$}
 \put(0,7.5){$\mathbb{= DMS}\rm t\mathbb{SH}_1$}
\put(2,9.8){\circle*{0.15}}
      \put(0.5,9.8){$\mathbb{DMS}\rm t\mathbb{SH}$}
\put(2,12){\circle*{0.15}}
\put(0.7,12){$\mathbb{DMSH}$}

\put(4,10){\circle*{0.15}}
  \put(2.5,9.8){$\mathbb{RDMSH}$}

\put(8,2){\circle*{0.15}}
\put(8.2,1.6){$\mathbb V(\mathbf{2}, \mathbf{\bar{2}})$}
\put(8,8){\circle*{0.15}}
\put(5.9,8.1){$\mathbb{RDQDS}\rm t\mathbb{SH}_1$}
\put(8,10){\circle*{0.15}}
   \put(6.1,10){$\mathbb{RDQDSH}_1$}
\put(8,12){\circle*{0.15}}
\put(6.3,12){$\mathbb{RDQDSH}$}
\put(8,14){\circle*{0.15}}
\put(8,14.3){$\mathbb{DQDSH}$}

\put(12,4){\circle*{0.15}}
\put(12.2,3.8){$\mathbb{RDPCS}\rm t\mathbb{SH}_1 = \mathbb V(C_{10}^{dp})$}
     \put(12.2,3.3){$= \mathbb{RDS}\rm t\mathbb{S}\rm t\mathbb{SH}_1$}
\put(12,6){\circle*{0.15}}
\put(12.2,5.8){$\mathbb{RDS}\rm t\mathbb{SH}_1 = \mathbb{RDS}\rm t\mathbb{SH}$}
\put(12.2,5.3){$ = \mathbb{RDPCSH}$}
\put(12,8){\circle*{0.15}}
\put(12.2,7.8){$\mathbb{DS}\rm t\mathbb{SH}_1 = \mathbb{DS}\rm t\mathbb{SH}$}
\put(12,12){\circle*{0.15}}
\put(12.2,11.8){$\mathbb{DPCSH}$}

\put(16,10){\circle*{0.15}}
\put(14.7,10.1){$\mathbb{JIDSH}_1$}
   \put(16,12){\circle*{0.15}}
\put(14.6,11.7){$\mathbb{JIDSH}$}

\put(2,6){\line(0,1){6}}
\put(8,8){\line(0,1){6}}
\put(12,4){\line(0,1){8}}
\put(16,10){\line(0,1){2}}

\put(2,6){\line(1,2){2}}
\put(2,6){\line(3,1){6}}
\put(2,12){\line(3,1){6}}
\put(4,10){\line(2,1){4}}
\put(8,2){\line(2,1){4}}
\put(12,8){\line(2,1){4}}

\put(2,6){\line(3,-2){6}}
\put(2,12){\line(1,-1){2}}
\put(8,8){\line(2,-1){4}}
\put(8,14){\line(2,-1){4}}
\put(8,14){\line(4,-1){8}}
	\put(7.2,0.5){{\large{Figure 5}}}
\end{picture}
}

\medskip

\medskip
\section{Concluding Remarks} \label{section_concluding_remarks}

\

It is, perhaps, worthwhile to mention here that we know from \cite{sankappanavar2011expansions} that every simple algebra in $\mathbb{RDQDS}\rm t\mathbb{H}_1$ is quasiprimal.
Of all the 25 simple algebras in $\mathbb{RDQDS}\rm t\mathbb{SH}_1$ (Section \ref{S7}), $\mathbf{2^e}$, 
 $\mathbf{\bar{2}^e}$, and $\mathbf{L_i}, i = 5,6,7,8$, and $\mathbf{D}_3$ are primal algebras and the rest, except $\mathbf{D}_1$ and $\mathbf{D}_2$, are semiprimal algebras.  \\

We will now collect here all the open problems that were mentioned in the earlier sections.\\

\noindent {\bf PROBLEM 1}:  Is $RDQDSH \subseteq BDQDSH$? In particular, is $RDQDStSH imply (B)$?\\
\noindent {\bf PROBLEM 2}: Is $\mathbb{DPCSH} \subset \mathbb{SBDQDSH}$?  In particular, is $\mathbb{DPCH} \subset \mathbb{SBDQDH}$?\\
\noindent {\bf PROBLEM 3}: Is the variety $\mathbb{DQDS}\rm{t}\mathcal{SH}_1$ (or even $\mathbb{DQDS}\rm{t}\mathbb{H}_1$) generated by its finite members?\\
\noindent {\bf PROBLEM 4}: Is the logic $\mathcal{RDMH}$  decidable?  We conjecture that  the logic $\mathcal{RDMH}_1$ is decidable.\\
\noindent {\bf PROBLEM 5}: Is the variety $\mathbb{DQDS}\rm{t}\mathcal{SH}_1$ (or even $\mathbb{DQDS}\rm{t}\mathbb{H}_1$) generated by its finite members?\\
\noindent{\bf PROBLEM 6}: Is the logic $\mathcal {DS}t$ decidable?  We, however, conjecture that $\mathcal {DS}t\mathcal{H}$ is decidable. \\


We conclude the paper by mentioning a few open-ended problems for future research. \\

{\bf PROBLEM}: Investigate the extensions of the logic $\mathcal{DHMSH}$ in relation to, among others, the following:

(a) Decidability,

(b) Various interpolation properties, 

(c) Definability, 

(d) Disjunction property, 

(e) Finite model property, 

(f) Finite embedability, 

(g) Structural completeness.  

(h) The property ES  (i.e., the epimorphisms are surjective) \\

\medskip
\section*{Acknowledgements} 

 The first author wants to thank the institutional support of CONICET (Consejo Nacional de Investigaciones Cient\'ificas y T\'ecnicas) and Universidad Nacional del Sur.  The authors would like to express their gratitude to Dr. Alex Citkin for his encouraging remarks and insightful comments on an earlier version of this paper.    





\bigskip


\end{document}